\documentclass[10pt]{amsart}

\usepackage{latexsym,exscale,enumerate,amsfonts,amssymb,mathtools}
\usepackage{amsmath,amsthm,amsfonts,amssymb,amscd,stmaryrd,textcomp}
\usepackage[final=true]{hyperref,bookmark}
\usepackage[normalem]{ulem}
\usepackage{young}
\usepackage{thmtools}
\SetSymbolFont{stmry}{bold}{U}{stmry}{m}{n}

\hypersetup{
  colorlinks = false,
  urlcolor = blue,
  pdfauthor = {David E.V. Rose and Daniel Tubbenhauer},
  pdfkeywords = {},
  pdftitle = {Symmetric webs, Jones--Wenzl recursions and \texorpdfstring{$q$}{q}-Howe duality},
  pdfsubject = {},
  pdfpagemode = UseNone,
  bookmarksopen = true,
  bookmarksopenlevel = 3,
  pdfdisplaydoctitle = true,
  linktocpage=true,
}

\usepackage[all]{xy}
\SelectTips{cm}{}

\usepackage{graphicx}

\usepackage{tikz}
\usetikzlibrary{decorations.markings}
\usetikzlibrary{decorations.pathreplacing}
\usetikzlibrary{arrows,shapes,positioning}
\tikzstyle directed=[postaction={decorate,decoration={markings,
    mark=at position #1 with {\arrow{>}}}}]
\tikzstyle rdirected=[postaction={decorate,decoration={markings,
    mark=at position #1 with {\arrow{<}}}}]
    
\def\cal#1{\mathcal{#1}}%
\newcommand{\sln}{\mathfrak{sl}_n}
\newcommand{\slm}{\mathfrak{sl}_m}
\newcommand{\slnn}[1]{\mathfrak{sl}_{#1}}
\newcommand{\gln}{\mathfrak{gl}_n}
\newcommand{\glm}{\mathfrak{gl}_m}
\newcommand{\gli}{\mathfrak{gl}_{\infty}}
\newcommand{\glnn}[1]{\mathfrak{gl}_{#1}}
\newcommand{\spid}[1]{\textbf{Sp}(\mathfrak{sl}_{{#1}})}

\newcommand{\bV}{\raisebox{0.03cm}{\mbox{\footnotesize$\textstyle{\bigwedge}$}}}
\newcommand{\SymSpf}{\textbf{SymSp}^f\!(\slnn{2})}
\newcommand{\SymSp}{\textbf{SymSp}(\slnn{2})}
\newcommand{\GenSp}{\textbf{GenSp}}
\newcommand{\Sym}{\mathrm{Sym}}
\newcommand{\ii}{\underline{i}}
\newcommand{\jj}{\underline{j}}
\newcommand{\ul}[1]{\underline{{#1}}}
\newcommand{\one}{\boldsymbol{1}}

\newcommand{\Uu}{\textbf{U}}
\newcommand{\Uq}{\textbf{U}_q}
\newcommand{\Ud}{\dot{\textbf{U}}_q}

\newcommand{\Rep}[1]{{#1}\text{-}\textbf{Mod}_{\wedge}}
\newcommand{\Repp}{\mathfrak{sl}_2\text{-}\textbf{fdMod}}

\newcommand{\Repn}{\mathfrak{sl}_n\text{-}\textbf{fdMod}}

\newcommand{\Kar}{\textbf{Kar}}
\newcommand{\Mat}{\textbf{Mat}}
\newcommand{\Hom}{\mathrm{Hom}}
\newcommand{\End}{\mathrm{End}}
\newcommand{\Gam}{\Gamma_{\mathbf{sym}}}

\def\C{{\mathbb C}}

\def\Z{{\mathbb Z}}
\def\Zg{{\mathbb{Z}_{\geq 0}}}

\newcommand{\qbin}[2]{{\textstyle\genfrac{[}{]}{0pt}{}{#1}{#2}}}
\newcommand{\qbinn}[2]{\genfrac{[}{]}{0pt}{}{#1}{#2}}

\theoremstyle{definition}
\newtheorem{theorem}{Theorem}[section]
\newtheorem{corollary}[theorem]{Corollary}
\newtheorem{conjecture}[theorem]{Conjecture}
\newtheorem{lemma}[theorem]{Lemma}
\newtheorem{proposition}[theorem]{Proposition}

\declaretheorem[style=definition,name=Example,qed=$\blacktriangle$,numberlike=theorem]{example}
\declaretheorem[style=definition,name=Definition,qed=$\blacktriangle$,numberlike=theorem]{definition}
\declaretheorem[style=definition,name=Remark,qed=$\blacktriangle$,numberlike=theorem]{remark}

\begin{document}
%

\title[Symmetric webs, Jones--Wenzl recursions and $q$-Howe duality]{Symmetric webs, Jones--Wenzl recursions and $q$-Howe duality}

\author{David E.~V.~Rose}
\address{Department of Mathematics, University of Southern California, Los Angeles, CA 90089, USA}
\email{davidero@usc.edu}

\author{Daniel Tubbenhauer}
\address{Centre for Quantum Geometry of Moduli Spaces, Aarhus University, Ny Munkegade 118, building 1530, room 316, 8000 Aarhus C, Denmark}
\email{dtubben@qgm.au.dk}


\vbadness=10001
\hbadness=10001
\begin{abstract}
We define and study the category of symmetric $\slnn{2}$-webs. 
This is a combinatorial description of the category of all finite-dimensional quantum $\slnn{2}$-modules. 
Explicitly, we show that (the additive closure of) the category of symmetric $\slnn{2}$-webs is (braided monoidally) equivalent to the latter. 
Our main tool is a quantum version of symmetric Howe duality. 
As a corollary of our construction, we provide new insight into Jones--Wenzl projectors and colored Jones polynomials.
\end{abstract}

\maketitle

\renewcommand{\thefootnote}{\arabic{footnote}}
\setcounter{footnote}{0}
\section{Introduction}\label{sec-intro}

\subsection{Temperley--Lieb categories and Jones--Wenzl projectors} \label{sub-jwtl}
A classical result of Rumer, Teller, and Weyl~\cite{rtw}, 
modernly interpreted, states that the \textit{Temperley--Lieb category} $\cal{TL}$ describes 
the full subcategory of quantum $\slnn{2}$-modules 
generated by tensor products of the two dimensional vector representation 
$\C_q^2 = \C(q)^2$ of quantum $\slnn{2}$,
which we denote\footnote{The notation $\Rep{\sln}$ more generally is used below to denote the 
full subcategory of quantum $\sln$-modules tensor 
generated by the fundamental representations (which, in the $\sln$ case, are exterior powers of the 
vector representation). Also, throughout the paper, when we refer to $\sln$-weights, $\sln$-modules, etc. we always 
mean their quantum versions. We also note that we only consider type $1$ (in the sense of Section 5.2 in~\cite{ja}) representations of 
quantum groups throughout. Moreover, for the insistent reader, all modules are finite-dimensional, left modules.} 
by $\Rep{\slnn{2}}$.
The former was first introduced in the study of statistical mechanics 
(as an algebra and also in the non-quantum setting) by Temperley and Lieb 
in~\cite{tl} and has played an important role in several areas of mathematics and physics.

Explicitly, the objects in $\cal{TL}$ are non-negative 
integers $k\in\mathbb{Z}_{\geq 0}$, and the morphisms are given graphically by 
$\C(q)$-linear combinations of non-intersecting planar tangle diagrams, 
which we view 
as mapping from the $k_1$ boundary points at the bottom 
of the tangle to the $k_2$ on the top, 
modulo boundary preserving isotopy and the local relation for evaluating a circle, that is,
\begin{equation}\label{eq-circle}
\xy
(0,0)*{
\begin{tikzpicture} [scale=1]
\draw[very thick] (-1,0) circle (0.3cm);
\end{tikzpicture}}
\endxy=-[2].
\end{equation}
Here, and throughout, $[a]$ for $a\in\Z$ denotes the \textit{quantum integer}, given by 
\[
[a]=\frac{q^a-q^{-a}}{q-q^{-1}}=q^{a-1}+q^{a-3}+\dots+q^{-a+3}+q^{-a+1}\in\Z[q,q^{-1}]
\] 
for $q$ a generic parameter. Note that $[0]=0$.

The correspondence between $\cal{TL}$ and the category $\Rep{\slnn{2}}$ associates the $\slnn{2}$-module $(\C_q^2)^{\otimes k}$
to $k \in \Zg$, and the morphisms are locally generated (by taking tensor products $\otimes$ and compositions $\circ$ of diagrams\footnote{Let us fix our diagrammatic conventions now: we read from left to right and bottom to top. 
Tensoring $u\otimes v$ is stacking picture $v$ to the right of $u$ and composition $v\circ u$ is given by stacking picture $v$ on the top of $u$.}) 
by the basic diagrams
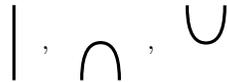
\[
\xy
(0,0)*{
\begin{tikzpicture} [scale=1]
\draw[very thick] (0,0) to (0,1);
\end{tikzpicture}
}
\endxy\quad , \quad
\xy
(0,0)*{
\begin{tikzpicture} [scale=1]
\draw[very thick] (0,0) to [out=90,in=180] (.25,.5) to [out=0,in=90] (.5,0);
\node at (.25,.9) {};
\end{tikzpicture}
}
\endxy\quad , \quad
\xy
(0,0)*{
\begin{tikzpicture} [scale=1]
\draw[very thick] (0,1) to [out=270,in=180] (.25,.5) to [out=0,in=270] (.5,1);
\node at (.25,0.1) {};
\end{tikzpicture}
}
\endxy
\]
where the first diagram corresponds to the identity, and the latter two correspond to the unique 
(up to scalar multiplication) $\slnn{2}$-intertwiners $\C_q^2 \otimes \C_q^2 \twoheadrightarrow\C_q$ and $\C_q\hookrightarrow \C_q^2 \otimes \C_q^2$,
where $\C_q$ denotes the trivial representation.
For example,
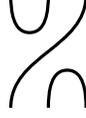
\[
\xy
(0,0)*{
\begin{tikzpicture} [scale=1]
\draw[very thick] (0,1.5) to [out=270,in=180] (.25,1) to [out=0,in=270] (.5,1.5);
\draw[very thick] (0,0) to [out=90,in=270] (1,1.5);
\draw[very thick] (.5,0) to [out=90,in=180] (.75,.5) to [out=0,in=90] (1,0);
\end{tikzpicture}
}
\endxy
\]
corresponds to a morphism $\C_q^2 \otimes \C_q^2 \otimes \C_q^2 \to \C_q^2 \otimes \C_q^2 \otimes \C_q^2$.
It turns out that the isotopy and circle removal~\eqref{eq-circle} relations 
generate all the relations between $\slnn{2}$-intertwiners under this correspondence.
That is, we have the following.

\begin{theorem}\label{thm-equisl2}
The category $\cal{TL}$ and $\Rep{\slnn{2}}$ are equivalent (as pivotal) categories.
\end{theorem}

It is known that every finite-dimensional, irreducible quantum 
$\slnn{2}$-module appears as a direct summand 
of $(\C_q^2)^{\otimes k}$ for some big enough $k\in\Z_{\geq 0}$. Thus, we obtain the 
entire category of finite-dimensional quantum $\slnn{2}$-modules, denoted by $\Repp$, 
by passing to the \textit{Karoubi envelope} $\Kar(\cal{TL})$ of $\cal{TL}$. Recall that the 
Karoubi envelope (sometimes also called 
idempotent completion) is the minimal enlargement 
of a category in which idempotents split; 
objects in the Karoubi envelope are (roughly) idempotent morphisms, 
which should be viewed as corresponding to their images.

It is a striking question if one can give a 
diagrammatic description of $\Kar(\cal{TL})$ as well. 

A solution to this question is known: an 
(in principle) explicit description of the entire 
category $\Repp$ can be given using the
\textit{Jones--Wenzl projectors} (also called Jones--Wenzl idempotents). These were introduced by Jones in~\cite{jones} 
and then further studied by Wenzl in~\cite{wenzl}. 
The Jones--Wenzl projectors are morphisms in $\cal{TL}$ which 
correspond to projecting onto, then including from, 
the highest weight irreducible summand $\Sym_q^k\C_q^2 \subset (\C_q^2)^{\otimes k}$. 
These projectors, which are usually depicted by a box 
with $k$ incoming and outgoing strands at the top and bottom respectively, 
admit a recursive definition describing the $k$-strand 
Jones--Wenzl projector $JW_k$ in terms of $(k-1)$-strand projectors as follows.
\begin{equation}\label{eq-jw}
\xy
(0,0)*{
\begin{tikzpicture} [scale=1]
\draw[very thick] (0,1.25) rectangle (2,.75);
\draw[very thick] (.25,0) to (.25,.75);
\draw[very thick] (.25,1.25) to (.25,2);
\draw[very thick] (1.75,0) to (1.75,.75);
\draw[very thick] (1.75,1.25) to (1.75,2);
\node at (1,1.75) {$\cdots$};
\node at (1,.25) {$\cdots$};
\node at (1,1) {\tiny $JW_k$};
\end{tikzpicture}
}
\endxy
\;\; = \;\;
\xy
(0,0)*{
\begin{tikzpicture} [scale=1]
\draw[very thick] (0,1.25) rectangle (1.5,.75);
\draw[very thick] (.25,0) to (.25,.75);
\draw[very thick] (.25,1.25) to (.25,2);
\draw[very thick] (1.25,0) to (1.25,.75);
\draw[very thick] (1.25,1.25) to (1.25,2);
\draw[very thick] (1.75,0) to (1.75,2);
\node at (.75,1.75) {$\cdots$};
\node at (.75,.25) {$\cdots$};
\node at (.75,1) {\tiny $JW_{k-1}$};
\end{tikzpicture}
}
\endxy 
\;\; + \;\;
\frac{[k-1]}{[k]} \;\;
\xy
(0,0)*{
\begin{tikzpicture} [scale=1]
\draw[very thick] (0,1.6) rectangle (1.5,1.25);
\draw[very thick] (0,.4) rectangle (1.5,.75);
\draw[very thick] (.25,.75) to (.25, 1.25);
\draw[very thick] (.25,0) to (.25,.4);
\draw[very thick] (.25,1.6) to (.25,2);
\draw[very thick] (1.25,0) to (1.25,.4);
\draw[very thick] (1.25,1.6) to (1.25,2);
\draw[very thick] (1.75,0) to [out=90,in=0] (1.5,.9) to [out=180,in=90] (1.25,.75);
\draw[very thick] (1.75,2) to [out=270,in=0] (1.5,1.1) to [out=180,in=270] (1.25,1.25);
\node at (.75,1.75) {$\cdots$};
\node at (.75,.25) {$\cdots$};
\node at (.75,1.4) {\tiny $JW_{k-1}$};
\node at (.75,0.55) {\tiny $JW_{k-1}$};
\end{tikzpicture}
}
\endxy
\end{equation}
We point out that some authors have a different sign convention here. Our convention comes from the fact that a circle evaluates to $-[2]$ instead of to $[2]$, see~\eqref{eq-circle1}.

However, working with such projectors in the 
Karoubi envelope quickly becomes cumbersome 
and computationally unmanageable due to their 
recursive definition. In this article, we provide a new, alternative 
diagrammatic description of the \textit{entire} category $\Repp$ of 
finite-dimensional quantum $\slnn{2}$-modules.

\subsection{A reminder on \texorpdfstring{$\sln$}{sln}-webs}\label{sub-reminder}

In pioneering work, see~\cite{kup}, Kuperberg 
extended the diagrammatic description of $\Rep{\slnn{2}}$ 
to the Lie algebra $\slnn{3}$ (and the other two 
rank $2$ Lie algebras of type $B_2$ and $G_2$ -- but we do not use them in this paper). 
Recall that the question was to find a diagrammatic 
and combinatorial model for $\Rep{\slnn{3}}$, the full subcategory of finite-dimensional 
quantum $\slnn{3}$-modules whose objects 
are finite tensor products of $\bV_q^k\C_q^3$'s, the fundamental $\slnn{3}$-modules\footnote{The notation $\bV_q^k$ means the \textit{quantum} alternating tensors. 
These are the quantum analogs of the classical alternating tensors, see for instance Subsection 4.2 in~\cite{ckm}.}. 
Since every irreducible $\slnn{3}$-module will appear as a summand of tensor products of $\bV^k_q\C_q^3$'s, 
we again have that ``morally'' the study of $\Rep{\slnn{3}}$ suffices to understand the entire category of finite-dimensional $\slnn{3}$-modules.

Kuperberg succeeded: he introduced in Section 4 of~\cite{kup} 
the \textit{$\slnn{3}$-spider}, denoted here by $\spid{3}$. 
This is a category whose morphisms, called \textit{$\slnn{3}$-webs}, 
are freely generated (via tensoring and composition) 
by local pieces of certain oriented, trivalent, planar graphs.
The category $\spid{3}$ is then obtained by taking a 
certain quotient, and the main difficulty is to find the ``correct'' relations 
such that there is an equivalence of (pivotal) categories
$\spid{3} \cong \Rep{\slnn{3}}$. 
Kuperberg gave the relations needed 
to obtain the aforementioned equivalence. 
While in the $\slnn{2}$ case the circle removal 
relation~\eqref{eq-circle} suffices, the $\slnn{3}$ case requires three local relations 
(that we do not need and thus, do not explicitly recall here).

It was long an open problem to extend Kuperberg's 
results to describe $\Rep{\sln}$, the 
full subcategory of all finite-dimensional
$\sln$-modules whose objects are 
finite tensor products of the fundamental 
$\sln$-representations $\bV^k_q\C_q^n$. 
As before, by Karoubi completing, it suffices 
to study $\Rep{\sln}$ to obtain a 
description of the entire category of finite-dimensional 
$\sln$-modules. A description of this 
subcategory in terms of $\sln$-webs was realized by 
Cautis, Kamnitzer and Morrison using the 
novel method of \textit{quantum skew Howe duality} (for short: $q$-skew Howe duality), see~\cite{ckm}. 
Our description of the entire category of finite-dimensional quantum $\slnn{2}$-modules in 
this paper is, surprisingly, 
related to Cautis, Kamnitzer and Morrison's $\sln$-webs, 
which we briefly recall now. 
Much more, of course, can be found in their paper.

In Theorem 3.3.1 from~\cite{ckm}, Cautis, Kamnitzer, and Morrison show that
$\Rep{\sln}$ is (pivotal) equivalent to the category of $\sln$-webs,  
a \textit{combinatorially} defined category in which objects are sequences in the
symbols $1^\pm, \dots, (n-1)^\pm$, and morphisms are given by $\C(q)$-linear combinations of
oriented, trivalent, planar graphs with edges labeled by $1,\dots,n-1$, 
such that the sum of the incoming and outgoing labels agree at each vertex. 
Moreover, by convention, the edges are directed outward at the bottom and
inward at the top if the corresponding boundary number is positive.

The correspondence between this diagrammatic category 
and the category of $\sln$-modules is 
given by associating a tensor 
product of fundamental $\sln$-modules 
and their duals to each sequence, with $k^+$ corresponding 
to $\bV_q^k \C_q^n$ and $k^{-}$ to its dual 
$(\bV_q^k \C_q^n)^*$. The generating $\sln$-webs are
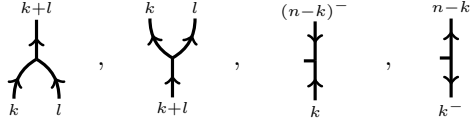
\[
\xy
(0,0)*{
\begin{tikzpicture}[scale=.3]
	\draw [very thick, directed=.55] (0, .75) to (0,2.5);
	\draw [very thick, directed=.45] (1,-1) to [out=90,in=330] (0,.75);
	\draw [very thick, directed=.45] (-1,-1) to [out=90,in=210] (0,.75); 
	\node at (0, 3) {\tiny $k{+}l$};
	\node at (-1,-1.5) {\tiny $k$};
	\node at (1,-1.5) {\tiny $l$};
\end{tikzpicture}
};
\endxy
\quad , \quad
\xy
(0,0)*{
\begin{tikzpicture}[scale=.3]
	\draw [very thick, directed=.55] (0,-1) to (0,.75);
	\draw [very thick, directed=.65] (0,.75) to [out=30,in=270] (1,2.5);
	\draw [very thick, directed=.65] (0,.75) to [out=150,in=270] (-1,2.5); 
	\node at (0, -1.5) {\tiny $k{+}l$};
	\node at (-1,3) {\tiny $k$};
	\node at (1,3) {\tiny $l$};
\end{tikzpicture}
};
\endxy
\quad , \quad
\xy
(0,0)*{
\begin{tikzpicture}[scale=.3]
	\draw [very thick, directed=.55] (0,-1.75) to (0,0);
	\draw[very thick, directed=.55] (0,1.75) to (0,0);
	\draw[very thick] (0,0) to (-.5,0);
	\node at (0,-2.25) {\tiny $k$};
	\node at (0,2.25) {\tiny $(n{-}k)^{-}$};
\end{tikzpicture}
};
\endxy
\quad , \quad
\xy
(0,0)*{
\begin{tikzpicture}[scale=.3]
	\draw [very thick, rdirected=.55] (0,-1.75) to (0,0);
	\draw[very thick, rdirected=.55] (0,1.75) to (0,0);
	\draw[very thick] (0,0) to (-.5,0);
	\node at (0,-2.25) {\tiny $k^{-}$};
	\node at (0,2.25) {\tiny $n{-}k$};
\end{tikzpicture}
};
\endxy
\]
which are called (reading from left to right) \textit{merge}, \textit{split}, 
\textit{tag in} and \textit{tag out}. 
These generators correspond to the unique (up to scalar) quantum $\sln$-intertwiners 
$\bV_q^k \C_q^n \otimes \bV_q^l \C_q^n \twoheadrightarrow \bV_q^{k+l} \C_q^n$, 
$\bV_q^{k+l} \C_q^n \hookrightarrow \bV_q^k \C_q^n \otimes \bV_q^l\C_q^n$, 
$\bV_q^k \C_q^n \xrightarrow{\cong} (\bV_q^{n-k}\C_q^n)^*$, 
and $(\bV_q^k \C_q^n)^* \xrightarrow{\cong} \bV_q^{n-k}\C_q^n$, see Section 3.2 in~\cite{ckm}.

As before, the main difficulty is deducing the correct collection of relations between these generators, 
which Cautis, Kamnitzer and Morrison give in Subsection 2.2 of~\cite{ckm}. 
The subset of their relations consisting of relations between 
``upward'' $\sln$-webs (i.e. those only factoring 
through tensor products of $\bV_q^k \C_q^n$'s, and not their duals) 
is of particular relevance to the current work, 
hence, we recall them now.

The upward relations are the following.
First, we have the \textit{associativity relations}:
\begin{equation}\label{eq-frob}
\xy,
(0,0)*{
\begin{tikzpicture}[scale=.3]
	\draw [very thick, directed=.45] (0,.75) to [out=90,in=220] (1,2.5);
	\draw [very thick, directed=.45] (1,-1) to [out=90,in=330] (0,.75);
	\draw [very thick, directed=.45] (-1,-1) to [out=90,in=210] (0,.75);
	\draw [very thick, directed=.45] (3,-1) to [out=90,in=330] (1,2.5);
	\draw [very thick, directed=.45] (1,2.5) to (1,4.25);
	\node at (-1,-1.5) {\tiny $h$};
	\node at (1,-1.5) {\tiny $k$};
	\node at (-1.375,1.5) {\tiny $h{+}k$};
	\node at (3,-1.5) {\tiny $l$};
	\node at (1,4.75) {\tiny $h{+}k{+}l$};
\end{tikzpicture}
};
\endxy=\xy
(0,0)*{
\begin{tikzpicture}[scale=.3]
	\draw [very thick, directed=.45] (0,.75) to [out=90,in=340] (-1,2.5);
	\draw [very thick, directed=.45] (-1,-1) to [out=90,in=210] (0,.75);
	\draw [very thick, directed=.45] (1,-1) to [out=90,in=330] (0,.75);
	\draw [very thick, directed=.45] (-3,-1) to [out=90,in=220] (-1,2.5);
	\draw [very thick, directed=.45] (-1,2.5) to (-1,4.25);
	\node at (1,-1.5) {\tiny $l$};
	\node at (-1,-1.5) {\tiny $k$};
	\node at (1.25,1.5) {\tiny $k{+}l$};
	\node at (-3,-1.5) {\tiny $h$};
	\node at (-1,4.75) {\tiny $h{+}k{+}l$};
\end{tikzpicture}
};
\endxy
\quad\text{and}\quad
\xy
(0,0)*{\rotatebox{180}{
\begin{tikzpicture}[scale=.3]
	\draw [very thick, rdirected=.55] (0,.75) to [out=90,in=220] (1,2.5);
	\draw [very thick, rdirected=.55] (1,-1) to [out=90,in=330] (0,.75);
	\draw [very thick, rdirected=.55] (-1,-1) to [out=90,in=210] (0,.75);
	\draw [very thick, rdirected=.55] (3,-1) to [out=90,in=330] (1,2.5);
	\draw [very thick, rdirected=.55] (1,2.5) to (1,4.25);
	\node at (-1,-1.5) {\rotatebox{180}{\tiny $l$}};
	\node at (1,-1.5) {\rotatebox{180}{\tiny $k$}};
	\node at (-1.325,1.5) {\rotatebox{180}{\tiny $k\! +\! l$}};
	\node at (3,-1.5) {\rotatebox{180}{\tiny $h$}};
	\node at (1,4.75) {\rotatebox{180}{\tiny $h\! +\! k\! +\! l$}};
\end{tikzpicture}
}};
\endxy=
\xy
(0,0)*{\reflectbox{\rotatebox{180}{
\begin{tikzpicture}[scale=.3]
	\draw [very thick, rdirected=.55] (0,.75) to [out=90,in=220] (1,2.5);
	\draw [very thick, rdirected=.55] (1,-1) to [out=90,in=330] (0,.75);
	\draw [very thick, rdirected=.55] (-1,-1) to [out=90,in=210] (0,.75);
	\draw [very thick, rdirected=.55] (3,-1) to [out=90,in=330] (1,2.5);
	\draw [very thick, rdirected=.55] (1,2.5) to (1,4.25);
	\node at (-1,-1.5) {\reflectbox{\rotatebox{180}{\tiny $h$}}};
	\node at (1,-1.5) {\reflectbox{\rotatebox{180}{\tiny $k$}}};
	\node at (-1.325,1.5) {\reflectbox{\rotatebox{180}{\tiny $h\! +\! k$}}};
	\node at (3,-1.5) {\reflectbox{\rotatebox{180}{\tiny $l$}}};
	\node at (1,4.75) {\reflectbox{\rotatebox{180}{\tiny $h\! +\! k\! +\! l$}}};
\end{tikzpicture}
}}};
\endxy
\end{equation}
To state the remaining relations,
define the so-called \textit{$F^{(j)}$ and $E^{(j)}$-ladders} as
\begin{equation}\label{eq-ladders}
\xy
(0,0)*{
\begin{tikzpicture}[scale=.3]
	\draw [very thick, directed=.55] (-2,-2) to (-2,0);
	\draw [very thick, directed=.55] (-2,0) to (-2,2);
	\draw [very thick, directed=.55] (2,-2) to (2,0);
	\draw [very thick, directed=.55] (2,0) to (2,2);
	\draw [very thick, directed=.55] (-2,0) to (2,0);
	\node at (-2,-2.5) {\tiny $k$};
	\node at (2,-2.5) {\tiny $l$};
	\node at (-2,2.5) {\tiny $k{-}j$};
	\node at (2,2.5) {\tiny $l{+}j$};
	\node at (0,0.75) {\tiny $j$};
	\node at (0,-2) {$F^{(j)}$};
\end{tikzpicture}
};
\endxy=
\xy
(0,0)*{
\begin{tikzpicture}[scale=.3]
	\draw [very thick, directed=.55] (1,.75) to (1,2);
	\draw [very thick] (2,-1) to [out=90,in=320] (1,.75);
	\draw [very thick, rdirected=0.65] (2,-1) to (2,-2);
	\draw [very thick, directed=.55] (-1,-2) to (-1,-.25);
	\draw [very thick] (-1,-0.25) to [out=150,in=270] (-2,1.5);
	\draw [very thick, directed=0.35] (-2,1.5) to (-2,2);
	\draw [very thick, directed=.55] (-1,-.25) to (1,.75);
	\node at (-1,-2.5) {\tiny $k$};
	\node at (2,-2.5) {\tiny $l$};
	\node at (-2,2.5) {\tiny $k{-}j$};
	\node at (1,2.5) {\tiny $l{+}j$};
	\node at (0,.75) {\tiny $j$};
\end{tikzpicture}
};
\endxy\quad\text{and}\quad
\xy
(0,0)*{
\begin{tikzpicture}[scale=.3]
	\draw [very thick, directed=.55] (-2,-2) to (-2,0);
	\draw [very thick, directed=.55] (-2,0) to (-2,2);
	\draw [very thick, directed=.55] (2,-2) to (2,0);
	\draw [very thick, directed=.55] (2,0) to (2,2);
	\draw [very thick, rdirected=.55] (-2,0) to (2,0);
	\node at (-2,-2.5) {\tiny $k$};
	\node at (2,-2.5) {\tiny $l$};
	\node at (-2,2.5) {\tiny $k{+}j$};
	\node at (2,2.5) {\tiny $l{-}j$};
	\node at (0,0.75) {\tiny $j$};
	\node at (0,-2) {$E^{(j)}$};
\end{tikzpicture}
};
\endxy=
\xy
(0,0)*{\reflectbox{
\begin{tikzpicture}[scale=.3]
	\draw [very thick, directed=.55] (1,.75) to (1,2);
	\draw [very thick] (2,-1) to [out=90,in=320] (1,.75);
	\draw [very thick, rdirected=0.65] (2,-1) to (2,-2);
	\draw [very thick, directed=.55] (-1,-2) to (-1,-.25);
	\draw [very thick] (-1,-0.25) to [out=150,in=270] (-2,1.5);
	\draw [very thick, directed=0.35] (-2,1.5) to (-2,2);
	\draw [very thick, directed=.55] (-1,-.25) to (1,.75);
	\node at (-1,-2.5) {\reflectbox{\tiny $l$}};
	\node at (2,-2.5) {\reflectbox{\tiny $k$}};
	\node at (-2,2.5) {\reflectbox{\tiny $l{-}j$}};
	\node at (1,2.5) {\reflectbox{\tiny $k{+}j$}};
	\node at (-0.25,.75) {\reflectbox{\tiny $j$}};
\end{tikzpicture}
}};
\endxy
\end{equation}
The remaining relations are (including a reflection of the right equation):
\begin{equation}\label{eq-simpler1}
\xy
(0,0)*{
\begin{tikzpicture}[scale=.3]
	\draw [very thick, directed=.55] (0,.75) to (0,2.5);
	\draw [very thick, directed=.55] (0,-2.75) to [out=30,in=330] (0,.75);
	\draw [very thick, directed=.55] (0,-2.75) to [out=150,in=210] (0,.75);
	\draw [very thick, directed=.55] (0,-4.5) to (0,-2.75);
	\node at (0,-5) {\tiny $k{+}l$};
	\node at (0,3) {\tiny $k{+}l$};
	\node at (-1.5,-1) {\tiny $k$};
	\node at (1.5,-1) {\tiny $l$};
\end{tikzpicture}
};
\endxy=\qbinn{k+l}{l}\xy
(0,0)*{
\begin{tikzpicture}[scale=.3]
	\draw [very thick, directed=.55] (0,-4.5) to (0,2.5);
	\node at (0,-5) {\tiny $k{+}l$};
	\node at (0,3) {\tiny $k{+}l$};
\end{tikzpicture}
};
\endxy
\quad\text{and}\quad
\xy
(0,0)*{
\begin{tikzpicture}[scale=.3]
	\draw [very thick, directed=.55] (-2,-4) to (-2,-2);
	\draw [very thick, directed=1] (-2,-2) to (-2,0.25);
	\draw [very thick, directed=.55] (2,-4) to (2,-2);
	\draw [very thick, directed=1] (2,-2) to (2,0.25);
	\draw [very thick, directed=.55] (-2,-2) to (2,-2);
	\draw [very thick] (-2,0.25) to (-2,2);
	\draw [very thick, directed=.55] (-2,2) to (-2,4);
	\draw [very thick] (2,0.25) to (2,2);
	\draw [very thick, directed=.55] (2,2) to (2,4);
	\draw [very thick, directed=.55] (-2,2) to (2,2);
	\node at (-2,-4.5) {\tiny $k$};
	\node at (2,-4.5) {\tiny $l$};
	\node at (-2.25,4.5) {\tiny $k{-}j_1{-}j_2$};
	\node at (2.25,4.5) {\tiny $l{+}j_1{+}j_2$};
	\node at (-3.5,0) {\tiny $k{-}j_1$};
	\node at (3.5,0) {\tiny $l{+}j_1$};
	\node at (0,-1.25) {\tiny $j_1$};
	\node at (0,2.75) {\tiny $j_2$};
\end{tikzpicture}
};
\endxy=\qbinn{j_1+j_2}{j_1}\xy
(0,0)*{
\begin{tikzpicture}[scale=.3]
	\draw [very thick, directed=.55] (-2,-2) to (-2,0);
	\draw [very thick, directed=.55] (-2,0) to (-2,2);
	\draw [very thick, directed=.55] (2,-2) to (2,0);
	\draw [very thick, directed=.55] (2,0) to (2,2);
	\draw [very thick, directed=.55] (-2,0) to (2,0);
	\node at (-2,-2.5) {\tiny $k$};
	\node at (2,-2.5) {\tiny $l$};
	\node at (-2.25,2.5) {\tiny $k{-}j_1{-}j_2$};
	\node at (2.25,2.5) {\tiny $l{+}j_1{+}j_2$};
	\node at (0,0.75) {\tiny $j_1{+}j_2$};
\end{tikzpicture}
};
\endxy
\end{equation}
which are called the \textit{digon removal} and 
\textit{square removal} relations. 
In these relations, the \textit{quantum binomial} is given by 
\[
\qbinn{a}{b}=\frac{[a][a-1]\cdots[a-b+2][a-b+1]}{[b]!}
\]
for $a\in\Z$, $b\in\Z_{\geq 0}$, 
where $[b]!=[1]\cdots [b-1][b]$ and by convention $[0]!=1$.
The final relations:
\begin{equation}\label{eq-almostsimpler}
\xy
(0,0)*{
\begin{tikzpicture}[scale=.3]
	\draw [very thick, directed=.55] (-2,-4) to (-2,-2);
	\draw [very thick, directed=1] (-2,-2) to (-2,0.25);
	\draw [very thick, directed=.55] (2,-4) to (2,-2);
	\draw [very thick, directed=1] (2,-2) to (2,0.25);
	\draw [very thick, directed=.55] (-2,-2) to (2,-2);
	\draw [very thick] (-2,0.25) to (-2,2);
	\draw [very thick, directed=.55] (-2,2) to (-2,4);
	\draw [very thick] (2,0.25) to (2,2);
	\draw [very thick, directed=.55] (2,2) to (2,4);
	\draw [very thick, rdirected=.55] (-2,2) to (2,2);
	\node at (-2,-4.5) {\tiny $k$};
	\node at (2,-4.5) {\tiny $l$};
	\node at (-2.25,4.5) {\tiny $k{-}j_1{+}j_2$};
	\node at (2.25,4.5) {\tiny $l{+}j_1{-}j_2$};
	\node at (-3.5,0) {\tiny $k{-}j_1$};
	\node at (3.5,0) {\tiny $l{+}j_1$};
	\node at (0,-1.25) {\tiny $j_1$};
	\node at (0,2.75) {\tiny $j_2$};
\end{tikzpicture}
};
\endxy=\sum_{j^{\prime}\geq 0}\qbinn{k-j_1-l+j_2}{j^{\prime}}\xy
(0,0)*{
\begin{tikzpicture}[scale=.3]
	\draw [very thick, directed=.55] (-2,-4) to (-2,-2);
	\draw [very thick, directed=1] (-2,-2) to (-2,0.25);
	\draw [very thick, directed=.55] (2,-4) to (2,-2);
	\draw [very thick, directed=1] (2,-2) to (2,0.25);
	\draw [very thick, rdirected=.55] (-2,-2) to (2,-2);
	\draw [very thick] (-2,0.25) to (-2,2);
	\draw [very thick, directed=.55] (-2,2) to (-2,4);
	\draw [very thick] (2,0.25) to (2,2);
	\draw [very thick, directed=.55] (2,2) to (2,4);
	\draw [very thick, directed=.55] (-2,2) to (2,2);
	\node at (-2,-4.5) {\tiny $k$};
	\node at (2,-4.5) {\tiny $l$};
	\node at (-2.25,4.5) {\tiny $k{-}j_1{+}j_2$};
	\node at (2.25,4.5) {\tiny $l{+}j_1{-}j_2$};
	\node at (-4.25,0) {\tiny $k{+}j_2{-}j^{\prime}$};
	\node at (4.25,0) {\tiny $l{-}j_2{+}j^{\prime}$};
	\node at (0,-1.25) {\tiny $j_2{-}j^{\prime}$};
	\node at (0,2.75) {\tiny $j_1{-}j^{\prime}$};
\end{tikzpicture}
};
\endxy
\end{equation}
are the \textit{square switch} relations.
For example, if $j_1=j_2=1$, then the only possible $j^{\prime}$ values
are $j^{\prime}=0,1$ and equation~\eqref{eq-almostsimpler}
gives\footnote{Note that we do not draw $\sln$-web edges labeled zero.}: 
\begin{equation}\label{eq-effe}
\xy
(0,0)*{
\begin{tikzpicture}[scale=.3]
	\draw [very thick, directed=.55] (-2,-4) to (-2,-2);
	\draw [very thick, directed=1] (-2,-2) to (-2,0.25);
	\draw [very thick, directed=.55] (2,-4) to (2,-2);
	\draw [very thick, directed=1] (2,-2) to (2,0.25);
	\draw [very thick, directed=.55] (-2,-2) to (2,-2);
	\draw [very thick] (-2,0.25) to (-2,2);
	\draw [very thick, directed=.55] (-2,2) to (-2,4);
	\draw [very thick] (2,0.25) to (2,2);
	\draw [very thick, directed=.55] (2,2) to (2,4);
	\draw [very thick, rdirected=.55] (-2,2) to (2,2);
	\node at (-2,-4.5) {\tiny $k$};
	\node at (2,-4.5) {\tiny $l$};
	\node at (-2,4.5) {\tiny $k$};
	\node at (2,4.5) {\tiny $l$};
	\node at (-3.5,0) {\tiny $k{-}1$};
	\node at (3.5,0) {\tiny $l{+}1$};
	\node at (0,-1.25) {\tiny $1$};
	\node at (0,2.75) {\tiny $1$};
\end{tikzpicture}
};
\endxy=\xy
(0,0)*{
\begin{tikzpicture}[scale=.3]
	\draw [very thick, directed=.55] (-2,-4) to (-2,-2);
	\draw [very thick, directed=1] (-2,-2) to (-2,0.25);
	\draw [very thick, directed=.55] (2,-4) to (2,-2);
	\draw [very thick, directed=1] (2,-2) to (2,0.25);
	\draw [very thick, rdirected=.55] (-2,-2) to (2,-2);
	\draw [very thick] (-2,0.25) to (-2,2);
	\draw [very thick, directed=.55] (-2,2) to (-2,4);
	\draw [very thick] (2,0.25) to (2,2);
	\draw [very thick, directed=.55] (2,2) to (2,4);
	\draw [very thick, directed=.55] (-2,2) to (2,2);
	\node at (-2,-4.5) {\tiny $k$};
	\node at (2,-4.5) {\tiny $l$};
	\node at (-2,4.5) {\tiny $k$};
	\node at (2,4.5) {\tiny $l$};
	\node at (-3.5,0) {\tiny $k{+}1$};
	\node at (3.5,0) {\tiny $l{-}1$};
	\node at (0,-1.25) {\tiny $1$};
	\node at (0,2.75) {\tiny $1$};
\end{tikzpicture}
};
\endxy+[k-l]\xy
(0,0)*{
\begin{tikzpicture}[scale=.3]
	\draw [very thick, directed=.55] (-2,-4) to (-2,4);
	\draw [very thick, directed=.55] (2,-4) to (2,4);
	\node at (-2,-4.5) {\tiny $k$};
	\node at (2,-4.5) {\tiny $l$};
	\node at (-2,4.5) {\tiny $k$};
	\node at (2,4.5) {\tiny $l$};
\end{tikzpicture}
};
\endxy
\end{equation}
The astute reader will recognize the similarity between these final relations and 
the relations 
\[
EF\one_{(k,l)}=FE\one_{(k,l)}+[k-l]\one_{(k,l)}
\] 
in the Beilinson, Lusztig and MacPherson's idempotented quantum group 
$\dot{\Uu}_q(\glnn{m})$ (see~\cite{blm}) recalled in detail below in Subsection~\ref{sub-howe}. 
Of course, this is no coincidence.
One of the main results of~\cite{ckm} is that $q$-skew Howe duality
induces a functor $\Phi^n_m\colon\Ud(\glm) \to \Ud^n(\glm) \to \Rep{\sln}$, 
where $\Ud^n(\glm)$ denotes the quotient of $\Ud(\glm)$ by 
the ideal generated by $\glm$-weights with entries not in $\{0,\dots,n\}$.

They go on to show in Proposition 5.2.2 of~\cite{ckm} that $\Phi^n_m$
factors through $\spid{n}$ and thus, taking the ``limit'' $m\to\infty$,
all the relations in $\spid{n}$ needed for the diagrammatic description of
$\Rep{\sln}$ follow from relations in $\dot{\Uu}_q(\glnn{\infty})$.
Our main idea in this paper is to adapt Cautis, Kamnitzer and Morrison's approach 
to quantum \textit{symmetric} Howe duality (for short, $q$-symmetric Howe duality).

\subsection{Main result}\label{sub-main}

We now introduce our new description of the representation 
theory of quantum $\slnn{2}$, the category of 
\textit{symmetric $\slnn{2}$-webs}. 
Following Kuperberg~\cite{kup}, we call this category \textit{the symmetric $\slnn{2}$-spider}, 
and denote it by $\SymSp$.

Here a symmetric 
$\slnn{2}$-web $u$ is an equivalence class (modulo 
boundary preserving planar isotopies) of edge-labeled, 
trivalent, planar graphs with boundary. The labels for 
the edges of $u$ are numbers from $\Z_{> 0}$ such that, 
at each trivalent vertex, two of the edge labels sum to the third.

We follow Cautis, Kamnitzer and Morrison and first introduce the \textit{free} 
symmetric $\slnn{2}$-spider. Then the 
symmetric $\slnn{2}$-spider $\SymSp$ is a certain quotient of it.

\begin{definition}\label{defn-fspid}(\textbf{The free symmetric $\slnn{2}$-spider}) 
The \textit{free symmetric $\slnn{2}$-spider}, 
which we denote by $\SymSpf$, is the category determined by the following data.
\begin{itemize}
\item The objects of $\SymSpf$ are tuples $\vec{k}\in\Z_{>0}^m$ 
for some $m\in\Z_{\geq 0}$, together with a zero object. 
We display their entries ordered from 
left to right according to their appearance in $\vec{k}$. 
Note that we allow $\varnothing$ as an object (corresponding to the empty 
sequence in $\Z_{>0}^0$), which is not to be confused with the zero object.

\item The morphisms of $\SymSpf$ from $\vec{k}$ to 
$\vec{l}$, denoted by $\Hom_{\SymSpf}(\vec{k},\vec{l})$, 
are diagrams with bottom boundary $\vec{k}$ and 
top boundary $\vec{l}$ freely 
generated as a $\C(q)$-vector space by all symmetric $\slnn{2}$-webs 
that can be obtained by composition $\circ$ 
(vertical gluing) and tensoring $\otimes$ 
(horizontal juxtaposition) of the following basic pieces (including the empty diagram $\varnothing$).
\begin{equation}\label{eq-symgen}
\xy
(0,0)*{
\begin{tikzpicture} [scale=1]
\draw[very thick] (0,0) to (0,1);
\node at (0,-.15) {\tiny $k$};
\node at (0,1.15) {\tiny $k$};
\end{tikzpicture}
}
\endxy
\quad , \quad
\xy
(0,0)*{
\begin{tikzpicture} [scale=1]
\draw[very thick] (0,0) to [out=90,in=180] (.25,.5) to [out=0,in=90] (.5,0);
\node at (.25,1.15) {};
\node at (0,-.15) {\tiny $k$};
\node at (.5,-.15) {\tiny $k$};
\end{tikzpicture}
}
\endxy
\quad , \quad
\xy
(0,0)*{
\begin{tikzpicture} [scale=1]
\draw[very thick] (0,1) to [out=270,in=180] (.25,.5) to [out=0,in=270] (.5,1);
\node at (.25,.-.15) {};
\node at (0,1.15) {\tiny $k$};
\node at (.5,1.15) {\tiny $k$};
\end{tikzpicture}
}
\endxy
\quad , \quad
\xy
(0,0)*{
\begin{tikzpicture}[scale=.3]
	\draw [very thick] (0, .75) to (0,2.5);
	\draw [very thick] (1,-1) to [out=90,in=330] (0,.75);
	\draw [very thick] (-1,-1) to [out=90,in=210] (0,.75); 
	\node at (0, 3) {\tiny $k{+}l$};
	\node at (-1,-1.5) {\tiny $k$};
	\node at (1,-1.5) {\tiny $l$};
\end{tikzpicture}
};
\endxy
\quad , \quad
\xy
(0,0)*{
\begin{tikzpicture}[scale=.3]
	\draw [very thick] (0,-1) to (0,.75);
	\draw [very thick] (0,.75) to [out=30,in=270] (1,2.5);
	\draw [very thick] (0,.75) to [out=150,in=270] (-1,2.5); 
	\node at (0, -1.5) {\tiny $k{+}l$};
	\node at (-1,3) {\tiny $k$};
	\node at (1,3) {\tiny $l$};
\end{tikzpicture}
};
\endxy
\end{equation}
These are called (from left to right) \textit{identity}, \textit{cap}, 
\textit{cup}, \textit{merge} and \textit{split}. \qedhere
\end{itemize}
\end{definition}

\begin{remark}\label{rem-conventions}
Note the following conventions and properties of $\SymSpf$.
\begin{itemize}
\item We consider the (free) symmetric $\slnn{2}$-webs up to boundary preserving isotopies. 
Formally, a (free) symmetric $\slnn{2}$-web is an equivalence class, but we abuse 
language and suppress this technical distinction.

\item The category is \textit{$\C(q)$-linear}, i.e. the spaces 
$\Hom_{\SymSpf}(\vec{k},\vec{l})$ are $\C(q)$-vector spaces 
and the composition $\circ$ is $\C(q)$-bilinear. 
Moreover, the category is \textit{monoidal} by juxtaposition of objects and morphisms, and $\otimes$ is 
similarly $\C(q)$-bilinear on morphism spaces.

\item The reading conventions for all symmetric 
$\slnn{2}$-webs is from \textit{bottom to top} and \textit{left to right}. 
That is, given $u,v\in\Hom_{\SymSpf}(\vec{k},\vec{l})$, then 
$v \circ u$ is obtained by gluing $v$ on top of $u$ and 
$u\otimes v$ is given by putting $v$ to the right of $u$. In pictures, e.g. we have
\[
\xy
(0,0)*{
\begin{tikzpicture} [scale=1]
\draw[very thick] (0,0) to [out=90,in=180] (.25,.5) to [out=0,in=90] (.5,0);
\node at (.25,.75) {};
\node at (0,-.15) {\tiny $k$};
\node at (.5,-.15) {\tiny $k$};
\end{tikzpicture}
}
\endxy
\circ
\xy
(0,0)*{
\begin{tikzpicture} [scale=1]
\draw[very thick] (0,1) to [out=270,in=180] (.25,.5) to [out=0,in=270] (.5,1);
\node at (.25,.25) {};
\node at (0,1.15) {\tiny $k$};
\node at (.5,1.15) {\tiny $k$};
\end{tikzpicture}
}
\endxy
=\xy
(0,0)*{
\begin{tikzpicture} [scale=1]
\draw[very thick] (-1,0) circle (0.3cm);
\node at (-0.5,0) {\tiny $k$};
\end{tikzpicture}}
\endxy\quad,\quad
\xy
(0,0)*{
\begin{tikzpicture} [scale=1]
\draw[very thick] (0,1) to [out=270,in=180] (.25,.5) to [out=0,in=270] (.5,1);
\node at (.25,.25) {};
\node at (0,1.15) {\tiny $k$};
\node at (.5,1.15) {\tiny $k$};
\end{tikzpicture}
}
\endxy
\circ
\xy
(0,0)*{
\begin{tikzpicture} [scale=1]
\draw[very thick] (0,0) to [out=90,in=180] (.25,.5) to [out=0,in=90] (.5,0);
\node at (.25,.75) {};
\node at (0,-.15) {\tiny $k$};
\node at (.5,-.15) {\tiny $k$};
\end{tikzpicture}
}
\endxy
=
\xy
(0,0)*{
\begin{tikzpicture}[scale=.3]
	\draw [very thick] (-1,2.5) to [out=270,in=180] (0,0.5) to [out=0,in=270] (1,2.5);
	\draw [very thick] (-1,-2.75) to [out=90,in=180] (0,-0.75) to [out=0,in=90] (1,-2.75);
	\node at (-1,3) {\tiny $k$};
	\node at (1,3) {\tiny $k$};
	\node at (-1,-3.15) {\tiny $k$};
	\node at (1,-3.15) {\tiny $k$};
\end{tikzpicture}
};
\endxy
\quad\text{and}\quad
\xy
(0,0)*{
\begin{tikzpicture}[scale=.3]
	\draw [very thick] (0,-1) to (0,.75);
	\draw [very thick] (0,.75) to [out=30,in=270] (1,2.5);
	\draw [very thick] (0,.75) to [out=150,in=270] (-1,2.5); 
	\draw [very thick] (1,-2.75) to [out=90,in=330] (0,-1);
	\draw [very thick] (-1,-2.75) to [out=90,in=210] (0,-1);
	\node at (-1,3) {\tiny $l_1$};
	\node at (1,3) {\tiny $l_2$};
	\node at (-1,-3.15) {\tiny $k_1$};
	\node at (1,-3.15) {\tiny $k_2$};
	\node at (-1.625,0) {\tiny $k_1{+}k_2$};
\end{tikzpicture}
};
\endxy
\otimes
\xy
(0,0)*{
\begin{tikzpicture} [scale=.3]
\draw[very thick] (0,-2.75) to (0,2.5);
\node at (0,-3.15) {\tiny $k_3$};
\node at (0,3) {\tiny $k_3$};
\end{tikzpicture}
}
\endxy
=
\xy
(-5,0)*{
\begin{tikzpicture}[scale=.3]
	\draw [very thick] (0,-1) to (0,.75);
	\draw [very thick] (0,.75) to [out=30,in=270] (1,2.5);
	\draw [very thick] (0,.75) to [out=150,in=270] (-1,2.5); 
	\draw [very thick] (1,-2.75) to [out=90,in=330] (0,-1);
	\draw [very thick] (-1,-2.75) to [out=90,in=210] (0,-1);
	\node at (-1,3) {\tiny $l_1$};
	\node at (1,3) {\tiny $l_2$};
	\node at (-1,-3.15) {\tiny $k_1$};
	\node at (1,-3.15) {\tiny $k_2$};
	\node at (-1.625,0) {\tiny $k_1{+}k_2$};
\end{tikzpicture}
};
(5,0)*{
\begin{tikzpicture} [scale=.3]
\draw[very thick] (0,-2.75) to (0,2.5);
\node at (0,-3.15) {\tiny $k_3$};
\node at (0,3) {\tiny $k_3$};
\end{tikzpicture}
}
\endxy
\]
where in the final equation $k_1+k_2=l_1+l_2$.

\item If any of the top boundary labels of the symmetric $\slnn{2}$-web $u$ 
are different from the corresponding bottom boundary component of the symmetric $\slnn{2}$-web 
$v$, then, by convention, $v\circ u$ is zero.\qedhere
\end{itemize}
\end{remark}

\begin{definition}\label{defn-spid}(\textbf{The symmetric $\slnn{2}$-spider}) 
The \textit{symmetric $\slnn{2}$-spider}, denoted by $\SymSp$, 
is the quotient category obtained from $\SymSpf$ by imposing the following local relations.
\begin{itemize}
\item The \textit{standard relations},
without orientations, that is, associativity~\eqref{eq-frob}, digon and square 
removals~\eqref{eq-simpler1} and the square switches~\eqref{eq-almostsimpler}. 
As before, it is convenient to define the $F^{(j)}$ and $E^{(j)}$-ladders 
as in~\eqref{eq-ladders}. 
In order to keep track of which is which, 
we (sometimes) add an orientation to the middle edges as a reminder, that is,
\[
\xy
(0,0)*{
\begin{tikzpicture}[scale=.3]
	\draw [very thick] (-2,-2) to (-2,0);
	\draw [very thick] (-2,0) to (-2,2);
	\draw [very thick] (2,-2) to (2,0);
	\draw [very thick] (2,0) to (2,2);
	\draw [very thick, directed=.55] (-2,0) to (2,0);
	\node at (-2,-2.5) {\tiny $k$};
	\node at (2,-2.5) {\tiny $l$};
	\node at (-2,2.5) {\tiny $k{-}j$};
	\node at (2,2.5) {\tiny $l{+}j$};
	\node at (0,0.75) {\tiny $j$};
	\node at (0,-2) {$F^{(j)}$};
\end{tikzpicture}
};
\endxy=\xy
(0,0)*{
\begin{tikzpicture}[scale=.3]
	\draw [very thick] (1,.75) to (1,2);
	\draw [very thick] (2,-1) to [out=90,in=320] (1,.75);
	\draw [very thick] (2,-1) to (2,-2);
	\draw [very thick] (-1,-2) to (-1,-.25);
	\draw [very thick] (-1,-0.25) to [out=150,in=270] (-2,1.5);
	\draw [very thick] (-2,1.5) to (-2,2);
	\draw [very thick] (-1,-.25) to (1,.75);
	\node at (-1,-2.5) {\tiny $k$};
	\node at (2,-2.5) {\tiny $l$};
	\node at (-2,2.5) {\tiny $k{-}j$};
	\node at (1,2.5) {\tiny $l{+}j$};
	\node at (0,.75) {\tiny $j$};
\end{tikzpicture}
};
\endxy\quad\text{and}\quad
\xy
(0,0)*{
\begin{tikzpicture}[scale=.3]
	\draw [very thick] (-2,-2) to (-2,0);
	\draw [very thick] (-2,0) to (-2,2);
	\draw [very thick] (2,-2) to (2,0);
	\draw [very thick] (2,0) to (2,2);
	\draw [very thick, rdirected=.55] (-2,0) to (2,0);
	\node at (-2,-2.5) {\tiny $k$};
	\node at (2,-2.5) {\tiny $l$};
	\node at (-2,2.5) {\tiny $k{+}j$};
	\node at (2,2.5) {\tiny $l{-}j$};
	\node at (0,0.75) {\tiny $j$};
	\node at (0,-2) {$E^{(j)}$};
\end{tikzpicture}
};
\endxy=
\xy
(0,0)*{\reflectbox{
\begin{tikzpicture}[scale=.3]
	\draw [very thick] (1,.75) to (1,2);
	\draw [very thick] (2,-1) to [out=90,in=320] (1,.75);
	\draw [very thick] (2,-1) to (2,-2);
	\draw [very thick] (-1,-2) to (-1,-.25);
	\draw [very thick] (-1,-0.25) to [out=150,in=270] (-2,1.5);
	\draw [very thick] (-2,1.5) to (-2,2);
	\draw [very thick] (-1,-.25) to (1,.75);
	\node at (-1,-2.5) {\reflectbox{\tiny $l$}};
	\node at (2,-2.5) {\reflectbox{\tiny $k$}};
	\node at (-2,2.5) {\reflectbox{\tiny $l{-}j$}};
	\node at (1,2.5) {\reflectbox{\tiny $k{+}j$}};
	\node at (-0.35,.75) {\reflectbox{\tiny $j$}};
\end{tikzpicture}
}};
\endxy
\]
By convention, if any label appearing in a symmetric $\slnn{2}$-web 
is less than zero, then the corresponding diagram is defined to be the zero morphism.

\item The \textit{symmetric relations}, that is, \textit{circle removal}:
\begin{equation}\label{eq-circle1}
\xy
(0,0)*{
\begin{tikzpicture} [scale=1]
\draw[very thick] (-1,0) circle (0.3cm);
\node at (-0.55,0) {\tiny $1$};
\end{tikzpicture}}
\endxy=-[2],
\end{equation}
and, finally, the \textit{dumbbell relation}:
\begin{equation}\label{eq-dumb}
\xy
(0,0)*{
\begin{tikzpicture}[scale=.3]
	\draw [very thick] (0,-1) to (0,.75);
	\draw [very thick] (0,.75) to [out=30,in=270] (1,2.5);
	\draw [very thick] (0,.75) to [out=150,in=270] (-1,2.5); 
	\draw [very thick] (1,-2.75) to [out=90,in=330] (0,-1);
	\draw [very thick] (-1,-2.75) to [out=90,in=210] (0,-1);
	\node at (-1,3) {\tiny $1$};
	\node at (1,3) {\tiny $1$};
	\node at (-1,-3.15) {\tiny $1$};
	\node at (1,-3.15) {\tiny $1$};
	\node at (-0.5,0) {\tiny $2$};
\end{tikzpicture}
};
\endxy=[2]
\xy
(0,0)*{
\begin{tikzpicture}[scale=.3] 
	\draw [very thick] (1,-2.75) to (1,2.5);
	\draw [very thick] (-1,-2.75) to (-1,2.5);
	\node at (-1,3) {\tiny $1$};
	\node at (1,3) {\tiny $1$};
	\node at (-1,-3.15) {\tiny $1$};
	\node at (1,-3.15) {\tiny $1$};
\end{tikzpicture}
};
\endxy
+
\xy
(0,0)*{
\begin{tikzpicture}[scale=.3]
	\draw [very thick] (-1,2.5) to [out=270,in=180] (0,0.5) to [out=0,in=270] (1,2.5);
	\draw [very thick] (-1,-2.75) to [out=90,in=180] (0,-0.75) to [out=0,in=90] (1,-2.75);
	\node at (-1,3) {\tiny $1$};
	\node at (1,3) {\tiny $1$};
	\node at (-1,-3.15) {\tiny $1$};
	\node at (1,-3.15) {\tiny $1$};
\end{tikzpicture}
};
\endxy
\end{equation}
\end{itemize}
\end{definition}

\begin{example}\label{ex-circle}
These relations, together with Lemma~\ref{lem-calc} 
below, imply that a $k$-labeled circle evaluates to $(-1)^{k}[k+1]$.
Indeed, we inductively compute:
\begin{align*}
\xy
(0,0)*{
\begin{tikzpicture} [scale=1]
\draw[very thick] (-1,0) circle (0.3cm);
\node at (-0.6,0) {\tiny $k$};
\end{tikzpicture}}
\endxy &=\frac{1}{[k]}
\xy
(0,0)*{
\begin{tikzpicture}[scale=.3]
 \draw [very thick] (0,.75) to (0,1.5);
	\draw [very thick] (0,-3.5) to (0,-2.75);
	\draw [very thick] (0,-2.75) to [out=30,in=330] (0,.75);
	\draw [very thick] (0,-2.75) to [out=150,in=210] (0,.75);
	\draw [very thick] (0,-3.5) to [out=270,in=180] (1,-5.5) to [out=0,in=270] (2,-3.5);
	\draw [very thick] (0,1.5) to [out=90,in=180] (1,3.5) to [out=0,in=90] (2,1.5);
	\draw [very thick] (2,-3.5) to (2,1.5);
	\node at (-1.75,0) {\tiny $k{-}1$};
	\node at (1.1,0) {\tiny $1$};
	\node at (2.3,0) {\tiny $k$};
	\end{tikzpicture}
};
\endxy=\frac{1}{[k]}
\xy
(0,0)*{
\begin{tikzpicture}[scale=.3]
	\draw [very thick] (0,-1) to (0,.75);
	\draw [very thick] (0,.75) to [out=30,in=270] (1,2.5);
	\draw [very thick] (0,.75) to [out=150,in=270] (-1,2.5); 
	\draw [very thick] (1,-2.75) to [out=90,in=330] (0,-1);
	\draw [very thick] (-1,-2.75) to [out=90,in=210] (0,-1);
	\draw [very thick] (1,-2.75) to [out=270,in=180] (2,-3.75) to [out=0,in=270] (3,-2.75);
	\draw [very thick] (-1,-2.75) to [out=270,in=180] (2,-4.75) to [out=0,in=270] (5,-2.75);
	\draw [very thick] (1,2.5) to [out=90,in=180] (2,3.5) to [out=0,in=90] (3,2.5);
	\draw [very thick] (-1,2.5) to [out=90,in=180] (2,4.5) to [out=0,in=90] (5,2.5);
	\draw [very thick] (3,-2.75) to (3,2.5);
	\draw [very thick] (5,-2.75) to (5,2.5);
	\node at (-0.5,0) {\tiny $k$};
	\node at (2.625,0) {\tiny $1$};
	\node at (6.1,0) {\tiny $k{-}1$};
\end{tikzpicture}
};
\endxy\\
&\stackrel{\ref{lem-calc}}{=}
\frac{1}{[k]}
\xy
(0,0)*{
\begin{tikzpicture}[scale=.3] 
	\draw[very thick] (1,-2.75) to [out=90,in=330] (0.25,-1);
	\draw[very thick] (-1,-2.75) to [out=30,in=210] (0.25,-1);
	\draw[very thick] (-1,-2.75) to [out=150,in=210] (-1,2.5);
	\draw [very thick] (0.25,-1) to (0.25,.75);
	\draw[very thick] (0.25,.75) to [out=30,in=270] (1,2.5);
	\draw[very thick] (0.25,.75) to [out=150,in=330] (-1,2.5);
	\draw [very thick] (1,-2.75) to [out=270,in=180] (2,-3.75) to [out=0,in=270] (3,-2.75);
	\draw [very thick] (-1,-2.75) to [out=270,in=180] (2,-4.75) to [out=0,in=270] (5,-2.75);
	\draw [very thick] (1,2.5) to [out=90,in=180] (2,3.5) to [out=0,in=90] (3,2.5);
	\draw [very thick] (-1,2.5) to [out=90,in=180] (2,4.5) to [out=0,in=90] (5,2.5);
	\draw [very thick] (3,-2.75) to (3,2.5);
	\draw [very thick] (5,-2.75) to (5,2.5);
	\node at (.75,0) {\tiny $2$};
	\node at (-1.375,0) {\tiny $k{-}2$};
	\node at (2.625,0) {\tiny $1$};
	\node at (6.1,0) {\tiny $k{-}1$};
\end{tikzpicture}
};
\endxy
-
\frac{[k-2]}{[k]}
\xy
(0,0)*{
\begin{tikzpicture}[scale=.3] 
	\draw [very thick] (1,-2.75) to (1,2.5);
	\draw [very thick] (-1,-2.75) to (-1,2.5);
	\draw [very thick] (1,-2.75) to [out=270,in=180] (2,-3.75) to [out=0,in=270] (3,-2.75);
	\draw [very thick] (-1,-2.75) to [out=270,in=180] (2,-4.75) to [out=0,in=270] (5,-2.75);
	\draw [very thick] (1,2.5) to [out=90,in=180] (2,3.5) to [out=0,in=90] (3,2.5);
	\draw [very thick] (-1,2.5) to [out=90,in=180] (2,4.5) to [out=0,in=90] (5,2.5);
	\draw [very thick] (3,-2.75) to (3,2.5);
	\draw [very thick] (5,-2.75) to (5,2.5);
	\node at (2.625,0) {\tiny $1$};
	\node at (6.1,0) {\tiny $k{-}1$};
\end{tikzpicture}
};
\endxy \\
&\stackrel{\eqref{eq-dumb}}{=} 
\frac{[2]}{[k]}
\xy
(0,0)*{
\begin{tikzpicture}[scale=.3] 
	\draw[very thick] (-1,-2) to [out=30,in=330] (-1,1.75);
	\draw[very thick] (-1,-2) to [out=150,in=210] (-1,1.75);
	\draw[very thick] (-1,-2.75) to (-1,-2);
	\draw[very thick] (-1,1.75) to (-1,2.5);
	\draw [very thick] (-1,-2.75) to [out=270,in=180] (2,-4.75) to [out=0,in=270] (5,-2.75);
	\draw [very thick] (-1,2.5) to [out=90,in=180] (2,4.5) to [out=0,in=90] (5,2.5);
	\draw [very thick] (5,-2.75) to (5,2.5);
	\node at (-2.5,1.5) {\tiny $k{-}2$};
	\node at (.25,0) {\tiny $1$};
	\node at (3.25,0) {\tiny $1$};
	\node at (6.1,0) {\tiny $k{-}1$};
	\draw[very thick] (2.25,0) circle (1.5);
\end{tikzpicture}
};
\endxy
+
\frac{1}{[k]}
\xy
(0,0)*{
\begin{tikzpicture}[scale=.3] 
	\draw[very thick] (-1,-2) to [out=30,in=330] (-1,1.75);
	\draw[very thick] (-1,-2) to [out=150,in=210] (-1,1.75);
	\draw[very thick] (-1,-2.75) to (-1,-2);
	\draw[very thick] (-1,1.75) to (-1,2.5);
	\draw [very thick] (-1,-2.75) to [out=270,in=180] (2,-4.75) to [out=0,in=270] (5,-2.75);
	\draw [very thick] (-1,2.5) to [out=90,in=180] (2,4.5) to [out=0,in=90] (5,2.5);
	\draw [very thick] (5,-2.75) to (5,2.5);
	\node at (.25,0) {\tiny $1$};
	\node at (6.1,0) {\tiny $k{-}1$};
	\node at (-2.5,1.5) {\tiny $k{-}2$};
\end{tikzpicture}
};
\endxy
+ (-1)^{k-1} [2] [k-2] \\
&= (-1)^{k} [2]^2[k-1]+ (-1)^{k-1} [2] [k-1]+ (-1)^{k-1} [2] [k-2] = (-1)^{k}[k+1],
\end{align*}
where the last equality follows from $[2][k^{\prime}]=[k^{\prime}+1]+[k^{\prime}-1]$ (for $k^{\prime}\geq 1$).
\end{example}

\begin{remark}\label{rem-TL}
Equation \eqref{eq-circle1} implies there is a functor $\cal{I}\colon\cal{TL} \rightarrow \SymSp$ given by sending objects $k \in \Z_{\geq 0}$ of 
$\cal{TL}$ to a sequence of $1$'s of length $k$, and by viewing morphisms in $\cal{TL}$ as 
symmetric $\slnn{2}$-webs. We will show below in the proof of Theorem~\ref{thm-main} that this functor is in fact an inclusion of a full subcategory.
\end{remark}

\begin{example}\label{ex-lollipop}
The so-called the \textit{lollipop relation}, that is,
\begin{equation}\label{eq-lollipop}
\xy
(0,0)*{
\begin{tikzpicture} [scale=0.3]
	\draw [very thick] (0,-1) to (0,.75);
	\draw [very thick] (0,.75) to [out=30,in=270] (1,2.5);
	\draw [very thick] (0,.75) to [out=150,in=270] (-1,2.5);
	\draw [very thick] (-1,2.5) to [out=90,in=180] (0,3.5) to [out=0,in=90] (1,2.5);  
	\node at (0, -1.5) {\tiny $2$};
	\node at (1.5, 2) {\tiny $1$};
\end{tikzpicture}
}
\endxy=0,
\end{equation}
can be deduced from the relations in the symmetric $\slnn{2}$-spider $\SymSp$:
\begin{align*}
\xy
(0,0)*{
\begin{tikzpicture} [scale=0.3]
	\draw [very thick] (0,-1) to (0,.75);
	\draw [very thick] (0,.75) to [out=30,in=270] (1,2.5);
	\draw [very thick] (0,.75) to [out=150,in=270] (-1,2.5);
	\draw [very thick] (-1,2.5) to [out=90,in=180] (0,3.5) to [out=0,in=90] (1,2.5);  
	\node at (0, -1.5) {\tiny $2$};
	\node at (1.5, 2) {\tiny $1$};
\end{tikzpicture}
}
\endxy
&
\stackrel{\eqref{eq-simpler1}}{=} 
\frac{1}{[2]}\xy
(0,0)*{
\begin{tikzpicture} [scale=0.3]
	\draw [very thick] (0,-3) to (0,-2);
	\draw [very thick] (0,-2) to [out=150,in=270] (-1,-1.115);
	\draw [very thick] (0,-2) to [out=30,in=270] (1,-1.115);
	\draw [very thick] (-1,-1.115) to [out=90,in=210] (0,-.25);
	\draw [very thick] (1,-1.115) to [out=90,in=330] (0,-.25);
	\draw [very thick] (0,-0.25) to (0,.75);
	\draw [very thick] (0,.75) to [out=30,in=270] (1,2.5);
	\draw [very thick] (0,.75) to [out=150,in=270] (-1,2.5);
	\draw [very thick] (-1,2.5) to [out=90,in=180] (0,3.5) to [out=0,in=90] (1,2.5);  
	\node at (0, -3.5) {\tiny $2$};
	\node at (1.5, 2) {\tiny $1$};
\end{tikzpicture}
}
\endxy\\
&
\stackrel{\eqref{eq-dumb}}{=} 
\xy
(0,0)*{
\begin{tikzpicture} [scale=0.3]
	\draw [very thick] (0,-3) to (0,-2);
	\draw [very thick] (0,-2) to [out=150,in=270] (-1,-1.115);
	\draw [very thick] (0,-2) to [out=30,in=270] (1,-1.115);
	\draw [very thick] (-1,-1.115) to(-1,2.5);
	\draw [very thick] (1,-1.115) to (1,2.5);
	\draw [very thick] (-1,2.5) to [out=90,in=180] (0,3.5) to [out=0,in=90] (1,2.5);  
	\node at (0, -3.5) {\tiny $2$};
	\node at (1.5, 2) {\tiny $1$};
\end{tikzpicture}
}
\endxy
+\frac{1}{[2]}\xy
(0,0)*{
\begin{tikzpicture} [scale=0.3]
	\draw [very thick] (0,-3) to (0,-2);
	\draw [very thick] (0,-2) to [out=150,in=270] (-1,-1.115);
	\draw [very thick] (0,-2) to [out=30,in=270] (1,-1.115);
	\draw [very thick] (-1,-1.115) to [out=90,in=180] (0,-.25) to [out=0,in=90] (1,-1.115);
	\draw [very thick] (-1,2.5) [out=270,in=180] to (0,.75) to [out=0,in=270] (1,2.5);
	\draw [very thick] (-1,2.5) to [out=90,in=180] (0,3.5) to [out=0,in=90] (1,2.5);  
	\node at (0, -3.5) {\tiny $2$};
	\node at (1.5, 2) {\tiny $1$};
\end{tikzpicture}
}
\endxy
\stackrel{\eqref{eq-circle1}}{=} 
\xy
(0,0)*{
\begin{tikzpicture} [scale=0.3]
	\draw [very thick] (0,-1) to (0,.75);
	\draw [very thick] (0,.75) to [out=30,in=270] (1,2.5);
	\draw [very thick] (0,.75) to [out=150,in=270] (-1,2.5);
	\draw [very thick] (-1,2.5) to [out=90,in=180] (0,3.5) to [out=0,in=90] (1,2.5);  
	\node at (0, -1.5) {\tiny $2$};
	\node at (1.5, 2) {\tiny $1$};
\end{tikzpicture}
}
\endxy
-
\xy
(0,0)*{
\begin{tikzpicture} [scale=0.3]
	\draw [very thick] (0,-1) to (0,.75);
	\draw [very thick] (0,.75) to [out=30,in=270] (1,2.5);
	\draw [very thick] (0,.75) to [out=150,in=270] (-1,2.5);
	\draw [very thick] (-1,2.5) to [out=90,in=180] (0,3.5) to [out=0,in=90] (1,2.5);  
	\node at (0, -1.5) {\tiny $2$};
	\node at (1.5, 2) {\tiny $1$};
\end{tikzpicture}
}
\endxy
=0.
\end{align*}
Note that Theorem~\ref{thm-main} below implies \textit{a priori} that 
equation~\eqref{eq-lollipop}
must hold, 
since there are no non-trivial $\slnn{2}$-intertwiners from $\Sym^{2}_q\C_q^2$ to the trivial representation $\C_q$.
\end{example}

\begin{remark}\label{rem-extra}
The following ``non-standard'' merge and split $\slnn{2}$-webs can be defined as composites 
of the generating morphisms in $\SymSp$.
\[
\xy
(0,0)*{
\begin{tikzpicture}[scale=.3]
	\draw [very thick] (0, .75) to (0,2.5);
	\draw [very thick] (1,-1) to [out=90,in=330] (0,.75);
	\draw [very thick] (-1,-1) to [out=90,in=210] (0,.75); 
	\node at (0, 3) {\tiny $l$};
	\node at (-1,-1.5) {\tiny $k$};
	\node at (1,-1.5) {\tiny $k{+}l$};
\end{tikzpicture}
};
\endxy
=
\xy
(0,0)*{\reflectbox{
\begin{tikzpicture}[scale=.3]
	\draw [very thick] (0,-1) to (0,.75);
	\draw [very thick] (0,.75) to [out=30,in=270] (1,2.5);
	\draw [very thick] (0,.75) to [out=150,in=270] (-1,2.5);
	\draw [very thick] (1,2.5) to [out=90,in=180] (2,3.5) to [out=0,in=90] (3,2.5);
	\draw [very thick] (3,-1) to (3,2.5);
	\draw [very thick] (-1,2.5) to (-1,3.5);
	\node at (0, -1.5) {\reflectbox{\tiny $k{+}l$}};
	\node at (-1,4) {\reflectbox{\tiny $l$}};
	\node at (3,-1.5) {\reflectbox{\tiny $k$}};
\end{tikzpicture}}
};
\endxy
\quad\text{and}\quad
\xy
(0,0)*{
\begin{tikzpicture}[scale=.3]
	\draw [very thick] (0, .75) to (0,2.5);
	\draw [very thick] (1,-1) to [out=90,in=330] (0,.75);
	\draw [very thick] (-1,-1) to [out=90,in=210] (0,.75); 
	\node at (0, 3) {\tiny $k$};
	\node at (-1,-1.5) {\tiny $k{+}l$};
	\node at (1,-1.5) {\tiny $l$};
\end{tikzpicture}
};
\endxy
=
\xy
(0,0)*{
\begin{tikzpicture}[scale=.3]
	\draw [very thick] (0,-1) to (0,.75);
	\draw [very thick] (0,.75) to [out=30,in=270] (1,2.5);
	\draw [very thick] (0,.75) to [out=150,in=270] (-1,2.5);
	\draw [very thick] (1,2.5) to [out=90,in=180] (2,3.5) to [out=0,in=90] (3,2.5);
	\draw [very thick] (3,-1) to (3,2.5);
	\draw [very thick] (-1,2.5) to (-1,3.5);
	\node at (0, -1.5) {\tiny $k{+}l$};
	\node at (-1,4) {\tiny $k$};
	\node at (3,-1.5) {\tiny $l$};
\end{tikzpicture}
};
\endxy
\]
Similarly for rotated versions.
\end{remark}

\begin{remark}\label{rem-another}
Of course, trivalent, planar graphs have previously appeared in the diagrammatic study of quantum $\slnn{2}$ 
under the guise of quantum spin networks, see~\cite{kl}. 
The difference in the present work is that we view 
trivalent vertices as the generators of our category (and 
deduce all relations between them needed to describe the category of representations)
rather than using trivalent vertices as shorthand for Temperley--Lieb diagrams built from Jones--Wenzl projectors.
\end{remark}

Recall that $\Repp$ denotes the category whose objects are (all) finite-dimensional 
$\slnn{2}$-modules, 
i.e. direct sums 
of the irreducible $\slnn{2}$-modules $\Sym^{k}_q\C_q^2$ (we explain the quantum symmetric tensors in Subsection~\ref{sub-howe} below), 
and whose morphisms are $\slnn{2}$-intertwiners between these tensor products. 
Recall that this is a monoidal category where $\otimes$ is the usual tensor product.

Moreover, recall that the \textit{additive closure} of 
a category $\mathcal C$ consist of finite, formal direct 
sums of objects from $\mathcal C$ with morphisms given by 
matrices whose entries are morphisms from $\mathcal C$.

\begin{theorem}\label{thm-main}
The additive closure\footnote{We must pass to the 
additive closure in order to make sense of direct 
sum decompositions. This is far more satisfying than passing 
to the Karoubi envelope of $\cal{TL}$ since working in 
the additive closure of a category $\mathcal C$ is combinatorially ``the same'' 
as working in $\mathcal C$.} 
of $\SymSp$ is monoidally equivalent to $\Repp$.
\end{theorem}

The functor $\Gam\colon\SymSp \to \Repp$ (see Definition~\ref{defn-mainfunctor}) 
inducing this equivalence is given by assigning the irreducible $\slnn{2}$-module $\Sym_q^{k}\C_q^2$ 
to the label $k$, and sending the generating morphisms in equation~\eqref{eq-symgen} to the (up to scalar) unique $\slnn{2}$-intertwiners 
between the $\slnn{2}$-modules corresponding to their boundaries.

In Section~\ref{sec-proofs}, we will prove Theorem~\ref{thm-main}. Of course, there are essentially two things to check: first, that the 
relations on symmetric $\slnn{2}$-webs are satisfied in the category of $\slnn{2}$-modules, and second, that we describe all morphisms 
(and relations between them) in this category. We accomplish the former task using $q$-symmetric Howe duality, and the latter 
by noticing the surprising result that the square switch relations~\eqref{eq-almostsimpler} gives the Jones--Wenzl recursion formula~\eqref{eq-jw}, 
a result which we think is of independent interest.

Finally, in Section~\ref{sec-colored} we use symmetric $\slnn{2}$-webs to compute the colored Jones polynomial, 
and discuss some further implications of our construction.
To do so, we show that $q$-symmetric Howe duality induces a \textit{braided monoidal} structure 
on our diagrammatic category $\SymSp$ and conclude that the functor $\Gam\colon\SymSp \to \Repp$ 
is an equivalence of braided monoidal categories.

We derive some consequences of this in Section~\ref{sec-colored}. 
For example, in Subsection~\ref{sub-mirror} we observe a connection between 
the $\Sym_q^k\C_q^2$-colored Jones polynomial and the $\bV_q^k\C_q^n$-colored Reshetikhin-Turaev polynomial of a 
colored, oriented link diagram $L_D$. 
For the precise statement see Theorem~\ref{thm-WeakMirrorSymmetry}.
\newline
~\newline
\noindent \textbf{Acknowledgements:}

Special thanks to Sabin Cautis for suggesting the study of a ``symmetric'' presentation of the category of 
$\slnn{2}$-modules.

We also thank Sergei Gukov, Matt Hogancamp, Greg Kuperberg, Aaron Lauda, Gus Lehrer, Gregor Masbaum, Weiwei Pan, Hoel Queffelec, Anonymous Referees,
Peter Samuelson, Antonio Sartori, Marko Sto\v{s}i\'{c}, Catharina Stroppel, 
Roland van der Veen, Pedro Vaz, Paul Wedrich and Geordie Williamson for helpful discussion, comments, and probing questions.

D.R. was supported by the John Templeton Foundation and NSF grant DMS-1255334 during part of this work. 
D.T. was partially supported by the center of excellence grant ``Centre for Quantum Geometry of Moduli Spaces (QGM)'' 
from the ``Danish National Research Foundation (DNRF)'', which also funded D.R.'s research visit, during which this collaboration began.
Additionally, D.R. would like to thank the QGM for their hospitality in hosting him during his research visit. 
D.T. wants to thank the food in India for providing him with enough soulfulness to keep on working on this paper.

\section{The proofs}\label{sec-proofs}

\subsection{\texorpdfstring{$q$}{q}-symmetric Howe duality}\label{sub-howe}

In this subsection, we present the requisite material on quantum groups and $q$-symmetric Howe duality. 
The main objective is to prove Corollary~\ref{cor-functor}, which gives a functor $\Phi_m\colon\Ud(\glm) \to \Repn$.
Along the way, we state $q$-symmetric Howe duality and deduce its consequences for 
any $n\in\Z_{>0}$ before we specialize to $n=2$. 
We use the results in this subsection to demonstrate later in Subsection~\ref{sub-diahowe} how the relations in the symmetric $\slnn{2}$-spider 
$\SymSp$ can be derived from $q$-symmetric Howe duality.

We begin by recalling the quantum general and special linear algebras, 
and their idempotented forms. 
The $\glm$-weight lattice is isomorphic to $\Z^m$. Let 
$\epsilon_i=(0,\ldots,1,\ldots,0)\in \Z^m$, with 
$1$ being in the $i$-th 
coordinate, and $\alpha_i=\epsilon_i-\epsilon_{i+1}
=(0,\ldots,1,-1,\ldots,0)\in\Z^{m}$, for 
$i=1,\ldots,m-1$. Recall that the Euclidean 
inner product on $\Z^m$ is defined by  
$(\epsilon_i,\epsilon_j)=\delta_{i,j}$.

\begin{definition}\label{defn-glm} 
For $m \in \Z_{>1}$, the \textit{quantum general linear algebra} 
$\Uq(\glm)$ is 
the associative, unital $\C(q)$-algebra 
generated by $L_i$ and $L_i^{-1}$, for $i = 1,\ldots, m$, 
and $E_{i}, F_i$, for $i=1,\ldots, m-1$, subject to the relations (for suitable $i,i_1,i_2$)
\begin{gather*}
L_{i_1}L_{i_2}=L_{i_2}L_{i_1},\quad 
L_iL_i^{-1}=L_i^{-1}L_i=1,\quad
L_{i_1}E_{i_2}=q^{(\epsilon_{i_1},\alpha_{i_2})}E_{i_2}L_{i_1},\quad
L_{i_1}F_{i_2}=q^{- (\epsilon_{i_1},\alpha_{i_2})}F_{i_2}L_{i_1},
\\
E_{i_1}F_{i_2} - F_{i_2}E_{i_1} = \delta_{i_1,i_2}\dfrac{L_{i_1}L_{i_1+1}^{-1}-L_{i_1}^{-1}L_{i_1+1}}{q-q^{-1}},
\\
E_{i_1}^2E_{i_2}-[2]E_{i_1}E_{i_2}E_{i_1}+E_{i_2}E_{i_1}^2=0,
\quad\text{if}\quad |i_1-i_2|=1, \qquad
E_{i_1}E_{i_2}-E_{i_2}E_{i_1}=0,\quad\text{else},
\\
F_{i_1}^2F_{i_2}-[2]F_{i_1}F_{i_2}F_{i_1}+F_{i_2}F_{i_1}^2=0,
\quad\text{if}\quad |i_1-i_2|=1, \qquad
F_{i_1}F_{i_2}-F_{i_2}F_{i_1}=0, \quad\text{else}.
\end{gather*}
The leftmost relations in the last two lines are the so-called \textit{(quantum) Serre-relations}.
\end{definition}
\begin{definition}\label{defn-slm} 
For $m \in \Z_{>1}$ the \textit{quantum special linear algebra} 
$\Uq(\slm)$ is the subalgebra of $\Uq(\glm)$ generated by the elements $E_i, F_i, K_i = L_i L_{i+1}^{-1},$ and $K_i^{-1}=L_{i+1}L_{i}^{-1}$ for $i=1,\dots,m-1$.
\end{definition}

To distinguish \textit{dominant integral} $\glm$-weights in $\Z_{\geq 0}^m$ 
(we call these, by abuse of language, just dominant integral $\glm$-weights, although a general 
dominant integral $\glm$-weight can have negative entries)
from general $\glm$-weights, we will denote the former by 
Greek letters as $\lambda$, $\mu$, etc. 
Recall that such $\glm$-weights $\lambda=(\lambda_1,\dots,\lambda_m)$ with $\lambda_i \geq 0$
can be described by partitions of $K$ where $\sum_{i=1}^m \lambda_i = K$. 
We denote the set of all partitions of $K$ of length $m$ by $\Lambda^+(m,K)$. Consequently, these dominant 
integral $\glm$-weights are precisely the elements of $\bigcup_{K\in\Z_{\geq 0}} \Lambda^+(m,K)$. 
We can picture such $\lambda$ as a Young diagram\footnote{We use the English convention 
for Young diagrams.}. For example, if 
$\lambda=(4,3,1,1)\in\Lambda^+(4,9)$, then
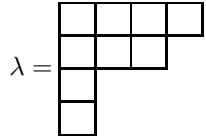
\[
\lambda=\xy(0,0)*{\begin{Young} & & &\cr & &\cr \cr\cr\end{Young}}\endxy
\]
where we abuse notation and denote the Young diagram and the partition by 
the same symbol. Thus, in our notation, dominant integral $\glm$-weights $\lambda$ are in 
bijective correspondence with Young diagram with at most $m$ rows, but 
with any possible (finite) number of columns.

Moreover, recall that $\Uq(\glm)$ has a unique 
highest weight module $V_m(\lambda)$ of highest weight $\lambda$
for each dominant integral $\glm$-weight $\lambda$. 
We point out that, by taking suitable tensors of the form
$V_m(\lambda)\otimes \det^{\otimes -k}$,
one can get any finite-dimensional, 
irreducible $\Uq(\glm)$-module. Here $\det^{\otimes -k}$ denotes 
a tensor product of length $k$ of the dual $\det^*=V_m(-1,\dots,-1)$ of the $1$-dimensional 
$\Uq(\glm)$-module $\det=V_m(1,\dots,1)\cong\bV_q^m\C_q^m$ (which is usually
called the \textit{determinant representation}). 
Thus, it suffices to study the $V_m(\lambda)$'s for most
purposes, including the remainder of this paper.

It is also worth noting that $\Uq(\glm)$ is a Hopf algebra with coproduct $\Delta$ given by
\[
\Delta(E_i)=E_i\otimes L_iL_{i+1}^{-1}+1\otimes E_i,\quad\Delta(F_i)=F_i\otimes 1+L^{-1}_iL_{i+1}\otimes F_i\quad\text{and}\quad\Delta(L_i)=L_i\otimes L_i.
\]
The antipode $S$ and the counit $\varepsilon$ are given by
\[
S(E_i)=-E_iL^{-1}_iL_{i+1},\quad S(F_i)=-L_iL_{i+1}^{-1}F_i,\quad S(L_i)=L_i^{-1},\quad\varepsilon(E_i)=\varepsilon(F_i)=0\quad\text{and}\quad\varepsilon(L_i)=1.
\]
The subalgebra $\Uq(\slm)$ inherits the Hopf algebra structure from $\Uq(\glm)$. 
We point out, since there are variations in different papers, that we use the same conventions as in~\cite{ckm}.
The Hopf algebra structure allows to extend actions to tensor products and duals of representations, 
and gives the existence of a trivial representation (that we simply denote as before by $\C_q$).

Another notion we need in the following is Beilinson, Lusztig and 
MacPherson's \textit{idempotented form}~\cite{blm}, denoted by $\Ud(\glm)$. 
Adjoin an idempotent $\one_{\vec{k}}$ for $\Uq(\glm)$
for each $\vec{k}\in\Z^{m}$ and add the relations
\begin{align*}
\one_{\vec{k}}\one_{\vec{l}} &= \delta_{\vec{k},\vec{l}}\one_{\vec{k}},   
\\
E_{i}\one_{\vec{k}} &= \one_{\vec{k}+\alpha_i}E_{i},\quad\text{with}\;
\alpha_i=(0,\ldots,1,-1,\ldots,0)\text{ as above},
\\
F_{i}\one_{\vec{k}} &= \one_{\vec{k}-\alpha_i}F_{i},\quad\text{with}\;
\alpha_i=(0,\ldots,1,-1,\ldots,0)\text{ as above},
\\
L_i\one_{\vec{k}} &= q^{k_i}\one_{\vec{k}}.
\end{align*}
\begin{definition}\label{defn-blm}The \textit{idempotented} quantum general linear algebra is defined by 
\[
\Ud(\glm)=\bigoplus_{\vec{k},\vec{l}\in\Z^{m}}\one_{\vec{l}}\Uq(\glm)\one_{\vec{k}}.
\]
\end{definition}

\begin{remark}\label{rem-relations-in-blm}
It is convenient to view $\Ud(\glm)$ as generated by the \textit{divided powers}
\[
F^{(j)}_i=\frac{F^j_i}{[j]!}\quad\text{ and }\quad E^{(j)}_i=\frac{E^j_i}{[j]!}
\]
for $i=1,\dots, m-1$. In particular, this point of view is useful if one wishes to work integrally, rather than over a field.
In this case, the integral form of $\Ud(\glm)$ is the $\Z[q,q^{-1}]$-subalgebra generated by divided powers 
and satisfying the following complete list of relations. 
In the following let $\vec{k}\in\Z^m$ and let all the subscripts be in $\{1,\dots,m-1\}$ and all the superscripts be in $\Zg$. 
If some of these indices fall outside of the sets mentioned above, then, by convention, the corresponding element is understood to be zero.

We have \textit{commutation relations} (with the left equations similarly for $E_i^{(j)}$'s)
\[
F^{(j_1)}_{i_1}F^{(j_2)}_{i_2}\one_{\vec{k}}=F^{(j_2)}_{i_2}F^{(j_1)}_{i_1} \one_{\vec{k}},\quad\text{if}\quad|i_1-i_2|>1,\quad F^{(j_1)}_{i_1}E^{(j_2)}_{i_2}\one_{\vec{k}}=E^{(j_2)}_{i_2}F^{(j_1)}_{i_1} \one_{\vec{k}},\quad\text{if}\quad|i_1-i_2|>0
\]
the \textit{Serre} and \textit{divided power relations} (with both equations similarly for $E_i^{(j)}$'s)
\[
F_{i_1}^2F_{i_2}\one_{\vec{k}}-[2]F_{i_1}F_{i_2}F_{i_1}\one_{\vec{k}}+F_{i_2}F_{i_1}^2\one_{\vec{k}}=0,
\quad\text{if}\quad |i_1-i_2|=1,\quad F_i^{(j_1)}F_i^{(j_2)}\one_{\vec{k}}=
\qbinn{j_1+j_2}{j_1}F_i^{(j_1+j_2)}\one_{\vec{k}},
\]
and the \textit{$EF-FE$-relations}
\[
E_i^{(j_2)}F_i^{(j_1)}\one_{\vec{k}}=\sum_{j^{\prime}}
\qbinn{k_i-j_1-k_{i+1}+j_2}{j^{\prime}}F_i^{(j_1-j^{\prime})}E_i^{(j_2-j^{\prime})}\one_{\vec{k}}.
\]
Here we note that the relations in $\Ud(\glm)$ imply that it suffices to specify the idempotents $\one_{\vec{k}}$ only once in an expression.
\end{remark} 

\begin{remark}\label{rem-UdotFun}
We will find it convenient to view $\Ud(\glm)$ as a category. Indeed, this is possible for any algebra containing a system of orthogonal 
idempotents. Explicitly, the objects of $\Ud(\glm)$ are precisely the $\glm$-weights $\vec{k}\in\Z^m$, and $\Hom(\vec{k},\vec{l}) = \one_{\vec{l}}\Uq(\glm)\one_{\vec{k}}$.
\end{remark}

We now discuss $q$-symmetric Howe duality, 
following the approach of Berenstein and Zwicknagl from~\cite{bz1}.
The ``classical'' symmetric Howe duality can be found 
in various sources, see~\cite{ho1} and~\cite{ho2} in the algebraic group setting and 
for example Theorem 5.16 in~\cite{cw1} for the pair $(\Uu(\glnn{m}),\Uu(\glnn{n}))$. 
Note that Cheng and Wang in Theorem 5.19 and Remark 5.20 of~\cite{cw1} also discuss 
\textit{super} Howe duality 
(which is more general and includes symmetric and skew Howe duality as a special case).
A slightly stronger result on super Howe duality which, in the non-quantized setting, 
comes close to what we need can be found in 
Proposition 2.1 of~\cite{sar1}.

Unfortunately, as in the $q$-skew Howe case, the literature about $q$-symmetric Howe duality
is very limited.
We hence adapt Cautis, Kamnitzer and Morrison's results on $q$-skew 
Howe duality to our setting, following closely their notation and exposition.

Denote the standard basis of the $\Uq(\glm)$-module $\C_q^m$ by 
$\{x_1,\dots,x_m\}$ with $\Uq(\glm)$-action given via
\begin{equation}\label{eq-action}
E_i(x_j)=\begin{cases}x_{j-1}, &\text{if }i=j-1,\\0, &\text{else,}\end{cases}\quad\quad F_i(x_j)=\begin{cases}x_{j+1}, &\text{if }i=j,\\0, &\text{else,}\end{cases}\quad\quad
L_i(x_j)=\begin{cases}qx_{j}, &\text{if }i=j,\\ 
x_j, &\text{else.}\end{cases}
\end{equation} 
By our conventions, the action of $\Uq(\slm)$ is almost 
the same as in~\eqref{eq-action}, but the $K_i$ act as 
$q^{+1}$ on $x_i$ and as $q^{-1}$ on $x_{i+1}$.

Now fix $m,n\in\Z_{>0}$. There is an action of $\Uq(\glm) \otimes \Uq(\mathfrak{gl}_n)$ 
on $\C_q^m\otimes \C_q^n$ 
and the latter has a basis given by $z_{ij}=x_i\otimes y_j$ for $x_i\in\C_q^m$ and $y_j\in\C_q^n$. 
The Hopf algebra structures of $\Uq(\glm)$ and $\Uq(\mathfrak{gl}_n)$ induce an action on the tensor algebra 
$\mathcal T(\C_q^m\otimes \C_q^n)$ of $\C_q^m\otimes \C_q^n$.

We now consider the \textit{quantum symmetric algebra}
\[
\mathcal \Sym_q^{\bullet}(\C_q^m\otimes \C_q^n)=\mathcal T(\C_q^m\otimes \C_q^n)/\bV_q^2(\C_q^m\otimes \C_q^n),
\]
where $\bV_q^2(\C_q^m\otimes \C_q^n)$ is the \textit{quantum exterior square} of $\C_q^m\otimes \C_q^n$.
Proposition 2.33 in~\cite{bz1} shows that $\bV_q^2(\C_q^m\otimes \C_q^n)$ is spanned by the elements
\begin{gather*}
z_{ij^{\prime}}\otimes z_{ij}-q z_{ij}\otimes z_{ij^{\prime}},\quad 
z_{ij^{\prime}}\otimes z_{i^{\prime}j}+q z_{i^{\prime}j^{\prime}}\otimes z_{ij}-qz_{ij}\otimes z_{i^{\prime}j^{\prime}}-q^2 z_{i^{\prime}j}\otimes z_{ij^{\prime}},\\
z_{i^{\prime}j}\otimes z_{ij}-q z_{ij}\otimes z_{i^{\prime}j},\quad 
z_{i^{\prime}j}\otimes z_{ij^{\prime}}+q z_{i^{\prime}j^{\prime}}\otimes z_{ij}-qz_{ij}\otimes z_{i^{\prime}j^{\prime}}-q^2 z_{ij^{\prime}}\otimes z_{i^{\prime}j},
\end{gather*}
for all $1\leq i< i^{\prime}\leq m$ and $1\leq j< j^{\prime}\leq n$.
The space $\Sym_q^{\bullet}(\C_q^m\otimes \C_q^n)$ is graded and its 
$k$-homogeneous piece, which we denote by $\Sym_q^{k}(\C_q^m\otimes \C_q^n)$, 
is the $k$-th \textit{quantum symmetric tensor} of $\C_q^m\otimes \C_q^n$. 
By setting $n=1$, we get the $k$-th \textit{quantum symmetric tensor}
of $\C_q^m$ denoted by $\Sym_q^k\C_q^m$. 
Similarly we have the \textit{quantum alternating tensors} $\bV_q^k\C_q^m$, $\bV_q^{\bullet}\C_q^m$, $\bV_q^k(\C_q^m\otimes \C_q^n)$ and $\bV_q^{\bullet}(\C_q^m\otimes \C_q^n)$ 
(we do not need the quantum alternating tensors much in this paper and refer to Subsection 4.2 in~\cite{ckm} for a more detailed treatment of these).

Our next result is a quantum version of symmetric Howe duality. 
We point out one crucial difference to the $q$-skew Howe case is that the direct sum decomposition in (3) of Theorem~\ref{thm-qSymHowe} 
\textit{does not} contain the transpose of $\lambda$.
To this end, we call a dominant integral $\glm$-weight $\lambda$ a \textit{$n$-supported $\glm$-weight} if its Young diagram has at most
$\min(m,n)$ rows, but still any possible (finite) number of columns.

\begin{theorem}(\textbf{$q$-symmetric Howe duality})\label{thm-qSymHowe} 
We have the following.
\begin{enumerate}
\item[(1)] For each $K \in\Zg$, the actions of $\Uq(\glm)$ and $\Uq(\mathfrak{gl}_n)$ on 
$\Sym_q^{K}(\C_q^m\otimes \C_q^n)$ commute and 
each generates the other's commutant.
\item[(2)] There is an isomorphism of $\Uq(\mathfrak{gl}_n)$-modules 
$\Sym_q^{\bullet}(\C_q^m\otimes \C_q^n) \cong (\Sym_q^{\bullet}\C_q^n)^{\otimes m}$ under which
the $\vec{k}$-weight space of $\Sym_q^{\bullet}(\C_q^m\otimes \C_q^n)$ (considered as a $\Uq(\glm)$-module)
is identified with 
$\Sym_q^{k_1}\C_q^n\otimes\cdots\otimes\Sym_q^{k_m}\C_q^n$ (here $\vec{k}=(k_1,\dots,k_m)$).
\item[(3)] As $\Uq(\glm) \otimes \Uq(\mathfrak{gl}_n)$-modules, 
we have a decomposition for each $K\in\Zg$ of the form
\[
\Sym_q^{K}(\C_q^m\otimes \C_q^n) \cong\bigoplus_{\lambda}V_m(\lambda)\otimes V_n(\lambda),
\]
where the $\bigoplus$ runs over all $n$-supported, dominant integral $\glm$-weights $\lambda\in\Lambda^+(m,K)$. 
This induces a $\Uq(\glm) \otimes \Uq(\mathfrak{gl}_n)$-module decomposition
\[
\Sym_q^{\bullet}(\C_q^m\otimes \C_q^n) \cong\bigoplus_{\lambda}V_m(\lambda)\otimes V_n(\lambda),
\]
where the $\bigoplus$ runs over all $n$-supported, dominant integral $\glm$-weights $\lambda$.\qedhere
\end{enumerate}
\end{theorem}

Note that, with the exception of the identification of the $\vec{k}$-weight space in item (2), this is essentially the quantum version of the Theorem in Section 2.1.2 in~\cite{ho2}.

\begin{proof}
The argument is essentially the same as that of Theorem 4.2.2 in~\cite{ckm}, with the exception that our task is easier, since 
from Proposition 2.33 in~\cite{bz1} we already know 
that $\Sym_q^{\bullet}(\C_q^m\otimes \C_q^n)$ is flat, i.e. 
the classical specialization of 
$\Sym_q^{\bullet}(\C_q^m\otimes \C_q^n)$ is $\Sym^{\bullet}(\C^m\otimes \C^n)$ 
and $\dim(\Sym_q^{K}(\C_q^m\otimes \C_q^n))=\dim(\Sym^{K}(\C^m\otimes \C^n))$ for all $K\in\Z_{\geq 0}$.
This then allows us to deduce (1) and (3) above as a consequence of the classical result which can be found, for example, 
in the Theorem of Section 2.1.2 in~\cite{ho2} or in Theorem 5.16 in~\cite{cw1}.

The isomorphism (2) is obtained by piecing together results from~\cite{bz1}. Explicitly, this is precisely their Proposition 4.2, 
using their Lemma 2.32 and Proposition 2.33.
To see that the $\vec{k}$-weight space decomposition holds we have to be
more explicit. Recall that Berenstein and Zwicknagl show that
$\Sym_q^k\C^n_q$ has a basis given by
\[
x_{j_1}\otimes \dots\otimes x_{j_k},\quad\text{for }1\leq j_1\leq\dots\leq  j_k\leq n,
\]
which we denote by $x_{\jj}$ for $\jj=(j_1,\dots,j_k)$.

Consider
\[
T_i\colon \Sym_q^{k}\C_q^n\to \Sym_q^{k}(\C_q^m\otimes \C_q^n),\quad x_{\jj}\mapsto z_{ij_1}\otimes\dots\otimes z_{ij_k}.
\]
for various $i=1,\dots,m$. These can be seen as sections of the $\Uq(\mathfrak{gl}_n)$-isomorphism given by 
Berenstein and Zwicknagl in Proposition 4.2 of~\cite{bz1}. From this, we see that
\[
T\colon\!\!
\displaystyle \bigoplus_{\sum_{i=1}^m\! k_i =K}  \Sym_q^{k_1}\C_q^n\otimes\cdots\otimes\Sym_q^{k_m}\C_q^n\to 
\Sym_q^{K}(\C_q^m\otimes \C_q^n),\quad v_1\otimes\cdots\otimes v_m\mapsto T_1(v_1)\otimes\cdots\otimes T_m(v_m)
\]
is an isomorphism of $\Uq(\mathfrak{gl}_n)$-modules (here $K=k_1+\dots+k_m$).

Since the action of $\Uq(\glm)$ on $\Sym_q^{\bullet}(\C_q^m\otimes \C_q^n)$ is ``row-wise,'' i.e.
\[
L_{i^{\prime}}(z_{ij_1}\otimes\dots\otimes z_{ij_k})=L_{i^{\prime}}(z_{ij_1})\otimes\dots\otimes L_{i^{\prime}}(z_{ij_k})
=
\begin{cases}q^k z_{ij_1}\otimes\dots\otimes z_{ij_k}, &\text{if }i=i^{\prime},\\
z_{ij_1}\otimes\dots\otimes z_{ij_k}, &\text{if }i\neq i^{\prime},\end{cases}
\] 
the $\vec{k}$-weight space identification follows. 
\end{proof}

By Theorem~\ref{thm-qSymHowe} part (2), we get 
linear maps
\begin{equation}\label{eq-symmor}
f_{\vec{k}}^{\vec{l}}\colon \one_{\vec{l}}\Ud(\glm)\one_{\vec{k}}\to\Hom_{\Uq(\mathfrak{gl}_n)}(\Sym_q^{k_1}\C_q^n\otimes\cdots\otimes\Sym_q^{k_m}\C_q^n,
\Sym_q^{l_1}\C_q^n\otimes\cdots\otimes\Sym_q^{l_m}\C_q^n)
\end{equation}
for any two $\vec{k},\vec{l}\in\Z_{>0}^m$ such that $\sum_{i=0}^mk_i=\sum_{i=0}^ml_i$. 
By part (1) of Theorem~\ref{thm-qSymHowe}, the homomorphisms $f_{\vec{k}}^{\vec{l}}$ are all surjective, which immediately implies 
that there exists a functor $\tilde\Phi_m\colon\Ud(\glm) \to \mathfrak{gl}_n\text{-}\textbf{fdMod}$, 
which sends\footnote{It sends all other blue objects $\vec{k}$ to the zero representation.} a $\glm$-weight $\vec{k}=(k_1,\dots,k_m)\in\Z^m_{\geq 0}$ 
to the $\Uq(\mathfrak{gl}_n)$-module $\Sym_q^{k_1}\C_q^n\otimes\cdots\otimes\Sym_q^{k_m}\C_q^n$ 
and morphisms $X\in \one_{\vec{l}}\Ud(\glm)\one_{\vec{k}}$ to $f_{\vec{k}}^{\vec{l}}(X)$, and is surjective on $\Hom$-spaces as in~\eqref{eq-symmor}.
Note that the functors $\tilde\Phi_m$ coming from Howe pairs of the form 
$(\Uq(\mathfrak{gl}_m),\Uq(\mathfrak{gl}_n))$ naturally map to categories of 
$\Uq(\mathfrak{gl}_n)$-modules and \textit{not} to categories of $\Uq(\mathfrak{sl}_n)$-modules
(see also the introduction of~\cite{qs1} or in Remark 1.1 of~\cite{tvw}).
We can, however, compose $\tilde\Phi_m$ with the restriction functor $\mathsf{res} \colon\gln\text{-}\mathbf{fdMod} \to \Repn$ 
to obtain the following.

\phantom{a}

\begin{corollary}\label{cor-functor}
There exists a functor $\Phi_m \!= \!\mathsf{res} \circ\tilde\Phi_m\colon\!\Ud(\glm)\!\to\!\Repn$, 
that we call the \textit{$q$-symmetric Howe functor}, 
which sends a $\glm$-weight $\vec{k}=(k_1,\dots,k_m)\in\Z^m_{\geq 0}$ 
to the $\Uq(\sln)$-module $\Sym_q^{k_1}\C_q^n\otimes\cdots\otimes\Sym_q^{k_m}\C_q^n$ 
and morphisms $X\in \one_{\vec{l}}\Ud(\glm)\one_{\vec{k}}$ to $f_{\vec{k}}^{\vec{l}}(X)$. 
The functor $\Phi_m$ is surjective on $\Hom$-spaces as in~\eqref{eq-symmor}.
\end{corollary}

\begin{proof}
This follows from the decomposition in part (3) of Theorem~\ref{thm-qSymHowe}, since we only 
restrict the action of $\Uq(\mathfrak{gl}_n)$ to $\Uq(\mathfrak{sl}_n)$. Thus, 
the dimensions of the $\Hom$-spaces as in~\eqref{eq-symmor} for both, $\Uq(\mathfrak{gl}_n)$-intertwiners and 
$\Uq(\mathfrak{sl}_n)$-intertwiners, can be read off from the dimensions of 
the irreducible $\Uq(\mathfrak{gl}_m)$-modules which appear in the decomposition from 
(3) of Theorem~\ref{thm-qSymHowe}.
\end{proof}

\begin{remark}\label{rem-new}
Note that the $q$-symmetric Howe functor is only surjective on $\Hom$-spaces as in~\eqref{eq-symmor}. 
For example, there are $\Uq(\mathfrak{sl}_2)$-intertwiners
\begin{gather*}
\mathsf{cap} \colon\C_q^2\otimes \C_q^2\twoheadrightarrow \C_q\quad\text{and}\quad
\mathsf{cup} \colon\C_q\hookrightarrow\C_q^2\otimes \C_q^2
\end{gather*}
given by projection (respectively inclusion) which do not arise under the $q$-symmetric Howe functor. Note that these are not 
$\Uq(\mathfrak{gl}_2)$-intertwiners since $\C_q^2\otimes \C_q^2\cong \bV^2_q\C_q^2\oplus\Sym^2_q\C_q^2$, 
but $\bV^2_q\C_q^2$ is the (non-trivial) determinant representation of $\Uq(\mathfrak{gl}_2)$.
As we will see in the proof of Theorem \ref{thm-main}, adding cap and cup generators suffices to recover 
the $\Uq(\mathfrak{sl}_2)$-intertwiners not in the image of $\Phi_m$.
\end{remark}

Denote by $\Ud^{\infty}(\glm)$ the quotient of $\Ud(\glm)$ by the ideal generated by 
all $\one_{\vec{k}}$ for $\glm$-weights $\vec{k}$ 
with a negative entry $k_i<0$.
By part (3) of Theorem~\ref{thm-qSymHowe}, 
all $\glm$-weights in $\Sym_q^{K}(\C_q^m\otimes \C_q^n)$ 
appear as $\glm$-weights appearing in $V_m(\lambda)$ where $\lambda$ is an $n$-supported, 
dominant integral $\glm$-weight.
Hence, the functors 
$(\Phi_m)_{m\in\Z_{\geq 0}}$ induce functors 
\begin{equation}\label{eq-functorweights}
\Phi^{\infty}_m\colon\Ud^{\infty}(\glm) \to \Repn,\quad \Phi^{\infty}_{\infty}
\colon
\Ud^{\infty}(\gli) :=  \lim_{\longrightarrow}
\Ud^{\infty}(\glm) 
\to \Repn. 
\end{equation}
By part (2) of Theorem~\ref{thm-qSymHowe} (and the restriction to $\Uq(\mathfrak{sl}_n)$), 
these functors are surjective 
on $\Hom$-spaces as in~\eqref{eq-symmor}. 
Since all irreducible $\slnn{2}$-modules are of the 
form $\Sym_q^k\C_q^2$ for some $k\in\Z_{\geq 0}$, we have the 
following more precise statement.
 
\begin{corollary}\label{cor-functor2}
The functor $\Phi^{\infty}_{\infty}\colon\Ud^{\infty}(\gli) \to \Repn$ is surjective 
on $\Hom$-spaces as in~\eqref{eq-symmor}. 
Moreover, for $n=2$ the induced functor from the additive closure (defined before Theorem~\ref{thm-main}) of $\Ud^{\infty}(\gli)$, that is,
\[
\Phi^{\infty}_{\infty}\colon\Mat(\Ud^{\infty}(\gli)) \to \Repp,
\]
is essentially surjective.
\end{corollary}

\begin{remark}\label{rem-doty}
We point out that this is the place where adapting the approach of Cautis, Kamnitzer and Morrison 
to the symmetric setting fails, due to the fact that there will be relations in 
$\mathfrak{gl}_n\text{-}\textbf{fdMod}$ (and hence, in $\Repn$ as well) that 
\textit{do not} come 
from $\Ud^{\infty}(\gli)$.

To this end, recall the dominance order $\unlhd$ for dominant integral
$\glm$-weights, given by setting $\mu\unlhd\lambda$ if and only if $\lambda-\mu$ is a 
$\Z_{\geq 0}$-linear combination of simple roots $\alpha_i$. 
Moreover, a not-necessarily dominant integral 
$\glm$-weight $\vec{k}$ is dominated by $\lambda$, denoted by
$\vec{k}\unlhd\lambda$, if and only if $\vec{k}$ appears in the Weyl group orbit of
a dominant integral $\glm$-weight $\mu$ with $\mu\unlhd\lambda$.

Let $\mathrm{I}_{\lambda}$ denote the ideal of $\Ud(\glm)$ generated by all $\one_{\vec{k}}$ for 
$\glm$-weights $\vec{k}$ that are not dominated by $\lambda$. 
Doty shows in Theorem 4.2 of~\cite{doty} that
\[
\Ud(\glm)/\mathrm{I}_{\lambda}\cong \bigoplus_{\mu\unlhd\lambda}\End_{\C_q}(V_m(\mu)).
\]
Here comes the catch: in part (3) of Theorem~\ref{thm-qSymHowe} we do not 
have all $V_m(\mu)$ appearing, but only those with $n$-supported $\mu$. 
Thus, in order to get faithfulness for the functor $\Phi_m^{\infty}$, one has to kill the endomorphism rings of 
the $V_m(\mu)$'s for non-$n$-supported $\mu$'s. Since this (clearly) depends on $n$, 
this introduces new relations which do not come from killing $\one_{\vec{k}}$ for certain $\gli$-weights $\vec{k}$.

Fortunately, in the $\slnn{2}$ case, it is easy to identify the missing relations, and in the following sections 
we show that they are exactly the symmetric relations from Definition~\ref{defn-spid}.
\end{remark}

\subsection{Jones--Wenzl recursion}\label{sub-jwrecursion}

In this subsection we show how the Jones--Wenzl 
recursion~\eqref{eq-jw} follows from the square switch 
relations~\eqref{eq-almostsimpler} and the dumbbell relation~\eqref{eq-dumb}.

\begin{definition}(\textbf{Symmetric Jones--Wenzl projectors})\label{defn-symJW}
Let $k\in\Z_{> 0}$. The $k$-th \textit{symmetric 
Jones--Wenzl projectors} $\mathcal{JW}_k$ is defined via
\[
\mathcal{JW}_k
=
\frac{1}{[k]!}
\;\;
\xy
(0,0)*{\reflectbox{
\begin{tikzpicture}[scale=.3]
	\draw [very thick] (0,-1) to (0,.75);
	\draw [very thick] (0,.75) to [out=30,in=270] (1,2.5);
	\draw [very thick] (0,.75) to [out=150,in=270] (-1,2.5); 
	\draw [very thick] (1,-2.75) to [out=90,in=330] (0,-1);
	\draw [very thick] (-1,-2.75) to [out=90,in=210] (0,-1);
	\draw [very thick] (1,2.5) to [out=30,in=270] (2,4.25);
	\draw [very thick] (1,2.5) to [out=150,in=270] (0,4.25);
	\draw [very thick] (2,-4.5) to [out=90,in=330] (1,-2.75);
	\draw [very thick] (0,-4.5) to [out=90,in=210] (1,-2.75);
	\draw [very thick] (2,4.25) to [out=30,in=270] (3,6.0);
	\draw [very thick] (2,4.25) to [out=150,in=270] (1,6.0);
	\draw [very thick] (3,-6.25) to [out=90,in=330] (2,-4.5);
	\draw [very thick] (1,-6.25) to [out=90,in=210] (2,-4.5);
	\draw [very thick] (-1,2.5) to (-1,6.0);
	\draw [very thick] (0,4.25) to (0,6.0);
	\draw [very thick] (-1,-2.5) to (-1,-6.25);
	\draw [very thick] (0,-4.25) to (0,-6.25);
	\node at (0,-6.75) {\tiny $\vdots$};
	\node at (0,7.25) {\tiny $\vdots$};
	\node at (-1.25,1.75) {\reflectbox{\tiny $1$}};
	\node at (2,1.75) {\reflectbox{\tiny $k{-}1$}};
	\node at (-1.25,-1.9) {\reflectbox{\tiny $1$}};
	\node at (2,-1.9) {\reflectbox{\tiny $k{-}1$}};
	\node at (0.4,0) {\reflectbox{\tiny $k$}};
	\node at (-0.25,3.5) {\reflectbox{\tiny $1$}};
	\node at (3,3.5) {\reflectbox{\tiny $k{-}2$}};
	\node at (-0.25,-3.65) {\reflectbox{\tiny $1$}};
	\node at (3,-3.65) {\reflectbox{\tiny $k{-}2$}};
	\node at (0.75,5.25) {\reflectbox{\tiny $1$}};
	\node at (4,5.25) {\reflectbox{\tiny $k{-}3$}};
	\node at (0.75,-5.4) {\reflectbox{\tiny $1$}};
	\node at (4,-5.4) {\reflectbox{\tiny $k{-}3$}};
\end{tikzpicture}
}};
\endxy
=
\xy
(0,0)*{
\begin{tikzpicture}[scale=.3]
	\draw [double] (0,-1) to (0,.75);
	\draw [very thick] (0,.75) to [out=30,in=270] (1,2.5);
	\draw [very thick] (0,.75) to [out=150,in=270] (-1,2.5); 
	\draw [very thick] (1,-2.75) to [out=90,in=330] (0,-1);
	\draw [very thick] (-1,-2.75) to [out=90,in=210] (0,-1);
	\node at (-1,3) {\tiny $1$};
	\node at (0.1,3) {$\cdots$};
	\node at (1,3) {\tiny $1$};
	\node at (-1,-3.15) {\tiny $1$};
	\node at (0.1,-3.15) {$\cdots$};
	\node at (1,-3.15) {\tiny $1$};
	\node at (-0.5,0) {\tiny $k$};
\end{tikzpicture}
};
\endxy
\]
where we repeatedly split a $k$-labeled edge until all of 
the top and bottom edges have label $1$. The rightmost picture 
above is a shorthand notation for $\mathcal{JW}_k$ where the ``doubled'' line should 
encode the coefficient $\frac{1}{[k]!}$.
\end{definition}

We need the following lemmata.

\begin{lemma}\label{lem-calc}
We have
\[
\xy
(0,0)*{
\begin{tikzpicture}[scale=.3]
	\draw [very thick] (0,-1) to (0,.75);
	\draw [very thick] (0,.75) to [out=30,in=270] (1,2.5);
	\draw [very thick] (0,.75) to [out=150,in=270] (-1,2.5); 
	\draw [very thick] (1,-2.75) to [out=90,in=330] (0,-1);
	\draw [very thick] (-1,-2.75) to [out=90,in=210] (0,-1);
	\node at (-1,3) {\tiny $k{-}1$};
	\node at (1,3) {\tiny $1$};
	\node at (-1,-3.15) {\tiny $k{-}1$};
	\node at (1,-3.15) {\tiny $1$};
	\node at (-0.5,0) {\tiny $k$};
\end{tikzpicture}
};
\endxy
=
\xy
(0,0)*{\reflectbox{
\begin{tikzpicture}[scale=.3]
	\draw [very thick] (-2,-4) to (-2,-2);
	\draw [very thick] (-2,-2) to (-2,0.25);
	\draw [very thick] (2,-4) to (2,-2);
	\draw [very thick] (2,-2) to (2,0.25);
	\draw [very thick, rdirected=.55] (-2,-2) to (2,-2);
	\draw [very thick] (-2,0.25) to (-2,2);
	\draw [very thick] (-2,2) to (-2,4);
	\draw [very thick] (2,0.25) to (2,2);
	\draw [very thick] (2,2) to (2,4);
	\draw [very thick, directed=.55] (-2,2) to (2,2);
	\node at (-2,-4.5) {\reflectbox{\tiny $1$}};
	\node at (2,-4.5) {\reflectbox{\tiny $k{-}1$}};
	\node at (-2,4.5) {\reflectbox{\tiny $1$}};
	\node at (2,4.5) {\reflectbox{\tiny $k{-}1$}};
	\node at (-2.5,0) {\reflectbox{\tiny $2$}};
	\node at (3.25,0) {\reflectbox{\tiny $k{-}2$}};
	\node at (0,-1.25) {\reflectbox{\tiny $1$}};
	\node at (0,2.75) {\reflectbox{\tiny $1$}};
\end{tikzpicture}
}};
\endxy-[k-2]\xy
(0,0)*{
\begin{tikzpicture}[scale=.3]
	\draw [very thick] (-2,-4) to (-2,4);
	\draw [very thick] (2,-4) to (2,4);
	\node at (-2,-4.5) {\tiny $k{-}1$};
	\node at (2,-4.5) {\tiny $1$};
	\node at (-2,4.5) {\tiny $k{-}1$};
	\node at (2,4.5) {\tiny $1$};
\end{tikzpicture}
};
\endxy
\]
for all $k\in\Z_{>2}$.
\end{lemma}

\begin{proof}
This is an immediate consequence of 
equation~\eqref{eq-effe}.
\end{proof}

\begin{lemma}\label{lem-recursion}
We have
\[
\mathcal{JW}_k=\xy
(0,0)*{
\begin{tikzpicture}[scale=.3]
	\draw [double] (0,-1) to (0,.75);
	\draw [very thick] (0,.75) to [out=30,in=270] (1,2.5);
	\draw [very thick] (0,.75) to [out=150,in=270] (-1,2.5); 
	\draw [very thick] (1,-2.75) to [out=90,in=330] (0,-1);
	\draw [very thick] (-1,-2.75) to [out=90,in=210] (0,-1);
	\node at (-1,3) {\tiny $1$};
	\node at (0.1,3) {$\cdots$};
	\node at (1,3) {\tiny $1$};
	\node at (-1,-3.15) {\tiny $1$};
	\node at (0.1,-3.15) {$\cdots$};
	\node at (1,-3.15) {\tiny $1$};
	\node at (-0.5,0) {\tiny $k$};
\end{tikzpicture}
};
\endxy=\xy
(0,0)*{
\begin{tikzpicture}[scale=.3]
	\draw [double] (0,-1) to (0,.75);
	\draw [very thick] (0,.75) to [out=30,in=270] (1,2.5);
	\draw [very thick] (0,.75) to [out=150,in=270] (-1,2.5); 
	\draw [very thick] (1,-2.75) to [out=90,in=330] (0,-1);
	\draw [very thick] (-1,-2.75) to [out=90,in=210] (0,-1);
	\draw [very thick] (3,-2.75) to (3,2.5);
	\node at (-1,3) {\tiny $1$};
	\node at (0.1,3) {$\cdots$};
	\node at (1,3) {\tiny $1$};
	\node at (-1,-3.15) {\tiny $1$};
	\node at (0.1,-3.15) {$\cdots$};
	\node at (1,-3.15) {\tiny $1$};
	\node at (-1.25,0) {\tiny $k{-}1$};
	\node at (3,-3.15) {\tiny $1$};
	\node at (3,3) {\tiny $1$};
\end{tikzpicture}
};
\endxy
+\frac{[k-1]}{[k]}
\xy
(0,0)*{
\begin{tikzpicture}[scale=.3]
	\draw [double] (0,-1) to (0,.75);
	\draw [very thick] (0,.75) to [out=30,in=270] (1,2.5);
	\draw [very thick] (0,.75) to [out=150,in=270] (-1,2.5); 
	\draw [very thick] (1,-2.75) to [out=90,in=330] (0,-1);
	\draw [very thick] (3,-2.75) to (3,2.5);
	\draw [very thick] (1,-2.75) to [out=270,in=180] (2,-4.5) to [out=0,in=270] (3,-2.75);
	\draw [very thick] (1,-6.75) to [out=90,in=180] (2,-5) to [out=0,in=90] (3,-6.75);
	\draw [double] (0,-10.25) to (0,-8.5);
	\draw [very thick] (0,-8.5) to [out=30,in=270] (1,-6.75);
	\draw [very thick] (0,-1) to [out=210,in=150] (0,-8.5); 
	\draw [very thick] (-0.75,-1.5) to [out=300,in=60] (-0.75,-8); 
	\draw [very thick] (1,-12) to [out=90,in=330] (0,-10.25);
	\draw [very thick] (-1,-12) to [out=90,in=210] (0,-10.25); 
	\draw [very thick] (3,-12) to (3,-6.75);
	\node at (-1,3) {\tiny $1$};
	\node at (0.1,3) {\tiny$\cdots$};
	\node at (1,3) {\tiny $1$};
	\node at (-1,-12.5) {\tiny $1$};
	\node at (0.1,-12.5) {\tiny$\cdots$};
	\node at (1,-12.5) {\tiny $1$};
	\node at (-1,0) {\tiny $k{-}1$};
	\node at (-1,-9.5) {\tiny $k{-}1$};
	\node at (-1.25,-1) {\tiny $k{-}2$};
	\node at (-1.25,-8.5) {\tiny $k{-}2$};
	\node at (-2.25,-4.75) {\tiny $1$};
	\node at (0.5,-4.75) {\tiny $1$};
	\node at (-.75,-4.75) {\tiny $\cdots$};
	\node at (3,-12.5) {\tiny $1$};
	\node at (3,3) {\tiny $1$};
\end{tikzpicture}
};
\endxy
\]
for all $k\in\Z_{>2}$.
\end{lemma}

\begin{proof}
Let $k\in\Z_{>2}$. Using Lemma~\ref{lem-calc}, we find
\[
\frac{1}{[k]!}
\xy
(0,0)*{\reflectbox{
\begin{tikzpicture}[scale=.3]
	\draw [very thick] (0,-1) to (0,.75);
	\draw [very thick] (0,.75) to [out=30,in=270] (1,2.5);
	\draw [very thick] (0,.75) to [out=150,in=270] (-1,2.5); 
	\draw [very thick] (1,-2.75) to [out=90,in=330] (0,-1);
	\draw [very thick] (-1,-2.75) to [out=90,in=210] (0,-1);
	\draw [very thick] (1,2.5) to [out=30,in=270] (2,4.25);
	\draw [very thick] (1,2.5) to [out=150,in=270] (0,4.25);
	\draw [very thick] (2,-4.5) to [out=90,in=330] (1,-2.75);
	\draw [very thick] (0,-4.5) to [out=90,in=210] (1,-2.75);
	\draw [very thick] (2,4.25) to [out=30,in=270] (3,6.0);
	\draw [very thick] (2,4.25) to [out=150,in=270] (1,6.0);
	\draw [very thick] (3,-6.25) to [out=90,in=330] (2,-4.5);
	\draw [very thick] (1,-6.25) to [out=90,in=210] (2,-4.5);
	\draw [very thick] (-1,2.5) to (-1,6.0);
	\draw [very thick] (0,4.25) to (0,6.0);
	\draw [very thick] (-1,-2.5) to (-1,-6.25);
	\draw [very thick] (0,-4.25) to (0,-6.25);
	\node at (0,-6.75) {\tiny $\vdots$};
	\node at (0,7.25) {\tiny $\vdots$};
	\node at (-1.3,1.75) {\reflectbox{\tiny $1$}};
	\node at (2,1.75) {\reflectbox{\tiny $k{-}1$}};
	\node at (-1.3,-1.9) {\reflectbox{\tiny $1$}};
	\node at (2,-1.9) {\reflectbox{\tiny $k{-}1$}};
	\node at (0.4,0) {\reflectbox{\tiny $k$}};
	\node at (-0.25,3.5) {\reflectbox{\tiny $1$}};
	\node at (3,3.5) {\reflectbox{\tiny $k{-}2$}};
	\node at (-0.25,-3.65) {\reflectbox{\tiny $1$}};
	\node at (3,-3.65) {\reflectbox{\tiny $k{-}2$}};
	\node at (0.75,5.25) {\reflectbox{\tiny $1$}};
	\node at (4,5.25) {\reflectbox{\tiny $k{-}3$}};
	\node at (0.75,-5.4) {\reflectbox{\tiny $1$}};
	\node at (4,-5.4) {\reflectbox{\tiny $k{-}3$}};
\end{tikzpicture}
}};
\endxy=
\frac{1}{[k]!}
\xy
(0,0)*{\reflectbox{
\begin{tikzpicture}[scale=.3]
	\draw [very thick] (1,2.5) to [out=30,in=270] (2,4.25);
	\draw [very thick] (1,2.5) to [out=150,in=270] (0,4.25);
	\draw [very thick] (2,-4.5) to [out=90,in=330] (1,-2.75);
	\draw [very thick] (0,-4.5) to [out=90,in=210] (1,-2.75);
	\draw [very thick] (2,4.25) to [out=30,in=270] (3,6.0);
	\draw [very thick] (2,4.25) to [out=150,in=270] (1,6.0);
	\draw [very thick] (3,-6.25) to [out=90,in=330] (2,-4.5);
	\draw [very thick] (1,-6.25) to [out=90,in=210] (2,-4.5);
	\draw [very thick] (-1,2.5) to (-1,6.0);
	\draw [very thick] (0,4.25) to (0,6.0);
	\draw [very thick] (-1,-2.5) to (-1,-6.25);
	\draw [very thick] (0,-4.25) to (0,-6.25);
	\draw [very thick] (1,-2.75) to (1,2.5);
	\draw [very thick] (-1,-2.75) to (-1,2.5);
	\draw [very thick, directed=.55] (-1,1) to (1,1);
	\draw [very thick, rdirected=.55] (-1,-1) to (1,-1);
	\node at (0,-6.75) {\tiny $\vdots$};
	\node at (0,7.25) {\tiny $\vdots$};
	\node at (-1.3,0) {\reflectbox{\tiny $2$}};
	\node at (2,0) {\reflectbox{\tiny $k{-}2$}};
	\node at (-1.3,2) {\reflectbox{\tiny $1$}};
	\node at (2,2) {\reflectbox{\tiny $k{-}1$}};
	\node at (-1.3,-2) {\reflectbox{\tiny $1$}};
	\node at (2,-2) {\reflectbox{\tiny $k{-}1$}};
	\node at (-0.25,3.5) {\reflectbox{\tiny $1$}};
	\node at (3,3.5) {\reflectbox{\tiny $k{-}2$}};
	\node at (-0.25,-3.65) {\reflectbox{\tiny $1$}};
	\node at (3,-3.65) {\reflectbox{\tiny $k{-}2$}};
	\node at (0.75,5.25) {\reflectbox{\tiny $1$}};
	\node at (4,5.25) {\reflectbox{\tiny $k{-}3$}};
	\node at (0.75,-5.4) {\reflectbox{\tiny $1$}};
	\node at (4,-5.4) {\reflectbox{\tiny $k{-}3$}};
\end{tikzpicture}
}};
\endxy - 
\frac{[k-2]}{[k]!}
\xy
(0,0)*{\reflectbox{
\begin{tikzpicture}[scale=.3]
	\draw [very thick] (1,2.5) to [out=30,in=270] (2,4.25);
	\draw [very thick] (1,2.5) to [out=150,in=270] (0,4.25);
	\draw [very thick] (2,-4.5) to [out=90,in=330] (1,-2.75);
	\draw [very thick] (0,-4.5) to [out=90,in=210] (1,-2.75);
	\draw [very thick] (2,4.25) to [out=30,in=270] (3,6.0);
	\draw [very thick] (2,4.25) to [out=150,in=270] (1,6.0);
	\draw [very thick] (3,-6.25) to [out=90,in=330] (2,-4.5);
	\draw [very thick] (1,-6.25) to [out=90,in=210] (2,-4.5);
	\draw [very thick] (-1,2.5) to (-1,6.0);
	\draw [very thick] (0,4.25) to (0,6.0);
	\draw [very thick] (-1,-2.5) to (-1,-6.25);
	\draw [very thick] (0,-4.25) to (0,-6.25);
	\draw [very thick] (1,-2.75) to (1,2.5);
	\draw [very thick] (-1,-2.75) to (-1,2.5);
	\node at (0,-6.75) {\tiny $\vdots$};
	\node at (0,7.25) {\tiny $\vdots$};
	\node at (-1.3,0) {\reflectbox{\tiny $1$}};
	\node at (2,0) {\reflectbox{\tiny $k{-}1$}};
	\node at (-0.25,3.5) {\reflectbox{\tiny $1$}};
	\node at (3,3.5) {\reflectbox{\tiny $k{-}2$}};
	\node at (-0.25,-3.65) {\reflectbox{\tiny $1$}};
	\node at (3,-3.65) {\reflectbox{\tiny $k{-}2$}};
	\node at (0.75,5.25) {\reflectbox{\tiny $1$}};
	\node at (4,5.25) {\reflectbox{\tiny $k{-}3$}};
	\node at (0.75,-5.4) {\reflectbox{\tiny $1$}};
	\node at (4,-5.4) {\reflectbox{\tiny $k{-}3$}};
\end{tikzpicture}
}};
\endxy
\]
There is now a dumbbell (with edge thickness $2$) in the middle picture, and 
we can use 
equation~\eqref{eq-dumb}
to simplify the above to
\[
\frac{1}{[k]!}
\xy
(0,0)*{\reflectbox{
\begin{tikzpicture}[scale=.3]
	\draw [very thick] (0,-1) to (0,.75);
	\draw [very thick] (0,.75) to [out=30,in=270] (1,2.5);
	\draw [very thick] (0,.75) to [out=150,in=270] (-1,2.5); 
	\draw [very thick] (1,-2.75) to [out=90,in=330] (0,-1);
	\draw [very thick] (-1,-2.75) to [out=90,in=210] (0,-1);
	\draw [very thick] (1,2.5) to [out=30,in=270] (2,4.25);
	\draw [very thick] (1,2.5) to [out=150,in=270] (0,4.25);
	\draw [very thick] (2,-4.5) to [out=90,in=330] (1,-2.75);
	\draw [very thick] (0,-4.5) to [out=90,in=210] (1,-2.75);
	\draw [very thick] (2,4.25) to [out=30,in=270] (3,6.0);
	\draw [very thick] (2,4.25) to [out=150,in=270] (1,6.0);
	\draw [very thick] (3,-6.25) to [out=90,in=330] (2,-4.5);
	\draw [very thick] (1,-6.25) to [out=90,in=210] (2,-4.5);
	\draw [very thick] (-1,2.5) to (-1,6.0);
	\draw [very thick] (0,4.25) to (0,6.0);
	\draw [very thick] (-1,-2.5) to (-1,-6.25);
	\draw [very thick] (0,-4.25) to (0,-6.25);
	\node at (0,-6.75) {\tiny $\vdots$};
	\node at (0,7.25) {\tiny $\vdots$};
	\node at (-1.3,1.75) {\reflectbox{\tiny $1$}};
	\node at (2,1.75) {\reflectbox{\tiny $k{-}1$}};
	\node at (-1.3,-1.9) {\reflectbox{\tiny $1$}};
	\node at (2,-1.9) {\reflectbox{\tiny $k{-}1$}};
	\node at (0.4,0) {\reflectbox{\tiny $k$}};
	\node at (-0.25,3.5) {\reflectbox{\tiny $1$}};
	\node at (3,3.5) {\reflectbox{\tiny $k{-}2$}};
	\node at (-0.25,-3.65) {\reflectbox{\tiny $1$}};
	\node at (3,-3.65) {\reflectbox{\tiny $k{-}2$}};
	\node at (0.75,5.25) {\reflectbox{\tiny $1$}};
	\node at (4,5.25) {\reflectbox{\tiny $k{-}3$}};
	\node at (0.75,-5.4) {\reflectbox{\tiny $1$}};
	\node at (4,-5.4) {\reflectbox{\tiny $k{-}3$}};
\end{tikzpicture}
}};
\endxy=
\frac{1}{[k]!}
\xy
(0,0)*{\reflectbox{
\begin{tikzpicture}[scale=.3]
	\draw [very thick] (1,2.5) to [out=30,in=270] (2,4.25);
	\draw [very thick] (1,2.5) to [out=150,in=270] (0,4.25);
	\draw [very thick] (2,-4.5) to [out=90,in=330] (1,-2.75);
	\draw [very thick] (0,-4.5) to [out=90,in=210] (1,-2.75);
	\draw [very thick] (2,4.25) to [out=30,in=270] (3,6.0);
	\draw [very thick] (2,4.25) to [out=150,in=270] (1,6.0);
	\draw [very thick] (3,-6.25) to [out=90,in=330] (2,-4.5);
	\draw [very thick] (1,-6.25) to [out=90,in=210] (2,-4.5);
	\draw [very thick] (-1,2.5) to (-1,6.0);
	\draw [very thick] (0,4.25) to (0,6.0);
	\draw [very thick] (-1,-2.5) to (-1,-6.25);
	\draw [very thick] (0,-4.25) to (0,-6.25);
	\draw [very thick] (1,-2.75) to (1,2.5);
	\draw [very thick] (-1,-2.75) to (-1,-1.5);
	\draw [very thick] (-1,1.5) to (-1,2.5);
	\draw [very thick] (-1,1.5) to [out=270,in=180] (0,0.75) to [out=0,in=270] (1,1.5);
	\draw [very thick] (-1,-1.5) to [out=90,in=180] (0,-0.75) to [out=0,in=90] (1,-1.5);
	\node at (0,-6.75) {\tiny $\vdots$};
	\node at (0,7.25) {\tiny $\vdots$};
	\node at (2.1,0) {\reflectbox{\tiny $k{-}2$}};
	\node at (-1.3,2) {\reflectbox{\tiny $1$}};
	\node at (2,2) {\reflectbox{\tiny $k{-}1$}};
	\node at (-1.3,-2) {\reflectbox{\tiny $1$}};
	\node at (2,-2) {\reflectbox{\tiny $k{-}1$}};
	\node at (-0.25,3.5) {\reflectbox{\tiny $1$}};
	\node at (3,3.5) {\reflectbox{\tiny $k{-}2$}};
	\node at (-0.25,-3.65) {\reflectbox{\tiny $1$}};
	\node at (3,-3.65) {\reflectbox{\tiny $k{-}2$}};
	\node at (0.75,5.25) {\reflectbox{\tiny $1$}};
	\node at (4,5.25) {\reflectbox{\tiny $k{-}3$}};
	\node at (0.75,-5.4) {\reflectbox{\tiny $1$}};
	\node at (4,-5.4) {\reflectbox{\tiny $k{-}3$}};
\end{tikzpicture}
}};
\endxy +
\frac{[2][k-1]-[k-2]}{[k]!}
\xy
(0,0)*{\reflectbox{
\begin{tikzpicture}[scale=.3]
	\draw [very thick] (1,2.5) to [out=30,in=270] (2,4.25);
	\draw [very thick] (1,2.5) to [out=150,in=270] (0,4.25);
	\draw [very thick] (2,-4.5) to [out=90,in=330] (1,-2.75);
	\draw [very thick] (0,-4.5) to [out=90,in=210] (1,-2.75);
	\draw [very thick] (2,4.25) to [out=30,in=270] (3,6.0);
	\draw [very thick] (2,4.25) to [out=150,in=270] (1,6.0);
	\draw [very thick] (3,-6.25) to [out=90,in=330] (2,-4.5);
	\draw [very thick] (1,-6.25) to [out=90,in=210] (2,-4.5);
	\draw [very thick] (-1,2.5) to (-1,6.0);
	\draw [very thick] (0,4.25) to (0,6.0);
	\draw [very thick] (-1,-2.5) to (-1,-6.25);
	\draw [very thick] (0,-4.25) to (0,-6.25);
	\draw [very thick] (1,-2.75) to (1,2.5);
	\draw [very thick] (-1,-2.75) to (-1,2.5);
	\node at (0,-6.75) {\tiny $\vdots$};
	\node at (0,7.25) {\tiny $\vdots$};
	\node at (-1.3,0) {\reflectbox{\tiny $1$}};
	\node at (2,0) {\reflectbox{\tiny $k{-}1$}};
	\node at (-0.25,3.5) {\reflectbox{\tiny $1$}};
	\node at (3,3.5) {\reflectbox{\tiny $k{-}2$}};
	\node at (-0.25,-3.65) {\reflectbox{\tiny $1$}};
	\node at (3,-3.65) {\reflectbox{\tiny $k{-}2$}};
	\node at (0.75,5.25) {\reflectbox{\tiny $1$}};
	\node at (4,5.25) {\reflectbox{\tiny $k{-}3$}};
	\node at (0.75,-5.4) {\reflectbox{\tiny $1$}};
	\node at (4,-5.4) {\reflectbox{\tiny $k{-}3$}};
\end{tikzpicture}
}};
\endxy
\]
where we point out that the additional contribution to the rightmost term above results 
after removing the extra $(k-2,1)$-digon.
A straightforward calculation shows that $[k]=[2][k-1]-[k-2]$ and taking
this into account, the rightmost term 
above is $\mathcal{JW}_{k-1}$ with an extra strand on the right.

To see that the other term works out as well, we 
iteratively ``explode'' the middle edge of thickness $k-2$ 
by using the digon removals~\eqref{eq-simpler1} the other way around, that is
\[
\xy
(0,0)*{
\begin{tikzpicture}[scale=.3]
	\draw [very thick] (0,-4.5) to (0,2.5);
	\node at (0,-5) {\tiny $k{-}2$};
	\node at (0,3) {\tiny $k{-}2$};
\end{tikzpicture}
};
\endxy=\frac{1}{[k-2]}
\xy
(0,0)*{
\begin{tikzpicture}[scale=.3]
	\draw [very thick] (0,.75) to (0,2.5);
	\draw [very thick] (0,-2.75) to [out=30,in=330] (0,.75);
	\draw [very thick] (0,-2.75) to [out=150,in=210] (0,.75);
	\draw [very thick] (0,-4.5) to (0,-2.75);
	\node at (0,-5) {\tiny $k{-}2$};
	\node at (0,3) {\tiny $k{-}2$};
	\node at (-1.5,-1) {\tiny $1$};
	\node at (2,-1) {\tiny $k{-}3$};
\end{tikzpicture}
};
\endxy=\frac{1}{[k-3][k-2]}
\xy
(0,0)*{
\begin{tikzpicture}[scale=.3]
	\draw [very thick] (0,.75) to (0,2.5);
	\draw [very thick] (0,-2.75) to [out=30,in=240] (1,-1.75);
	\draw [very thick] (1,-0.25) to [out=120,in=330] (0,.75);
	\draw [very thick] (1,-1.75) to [out=30,in=330] (1,-0.25);
	\draw [very thick] (1,-1.75) to [out=150,in=210] (1,-0.25);
	\draw [very thick] (0,-2.75) to [out=150,in=210] (0,.75);
	\draw [very thick] (0,-4.5) to (0,-2.75);
	\node at (0,-5) {\tiny $k{-}2$};
	\node at (0,3) {\tiny $k{-}2$};
	\node at (-1.5,-1) {\tiny $1$};
	\node at (1.75,-2.5) {\tiny $k{-}3$};
	\node at (1.75,0.5) {\tiny $k{-}3$};
	\node at (2.375,-1) {\tiny $k{-}4$};
	\node at (0.25,-1) {\tiny $1$};
\end{tikzpicture}
};
\endxy=\cdots=\frac{1}{[k-2]!}\xy
(0,0)*{
\begin{tikzpicture}[scale=.3]
	\draw [very thick] (0,.75) to (0,2.5);
	\draw [very thick] (0,-2.75) to [out=30,in=330] (0,.75);
	\draw [very thick] (0,-2.75) to [out=150,in=210] (0,.75);
	\draw [very thick] (0,-4.5) to (0,-2.75);
	\node at (0,-5) {\tiny $k{-}2$};
	\node at (0,3) {\tiny $k{-}2$};
	\node at (-1.5,-1) {\tiny $1$};
	\node at (0,-1) {\tiny $\cdots$};
	\node at (1.5,-1) {\tiny $1$};
\end{tikzpicture}
};
\endxy
\]
where we continue until all edges are of thickness $1$. 
The diagram now has the desired form in the statement. 
To see that the coefficient works out, note that 
\[
\frac{1}{[k]!}\frac{1}{[k-2]!}=\frac{1}{[k-1]!}\frac{1}{[k-1]!}\frac{[k-1]}{[k]}
\]
and the two factors $\frac{1}{[k-1]!}$ 
give the two symmetric Jones--Wenzl projectors $\mathcal{JW}_{k-1}$.
\end{proof}

Using these lemmata, we now deduce the main result of this subsection.

\begin{proposition}\label{prop-jw}
The symmetric Jones--Wenzl projectors $\mathcal{JW}_k$ are the images of the Jones--Wenzl projectors $JW_k$ in $\cal{TL}$ under the 
functor $\cal{I}\colon\cal{TL} \to\SymSp$, i.e. $\cal{I}(JW_k)=\mathcal{JW}_k$.
\end{proposition}

\begin{proof}
This follows since Lemma~\ref{lem-recursion} and
equation~\eqref{eq-dumb}
show that $\mathcal{JW}_k$ satisfy the Jones--Wenzl recursion~\eqref{eq-jw}, which uniquely determines $JW_k$.
\end{proof}

\begin{remark}
This gives the surprising result that, save for the base case, the Jones--Wenzl recursion exactly corresponds to the $\slnn{2}$-relations $EF\one_{(k,l)}=FE\one_{(k,l)}+[k-l]\one_{(k,l)}$.
\end{remark}

\begin{corollary}\label{cor-idem}
We have
\[
\mathcal{JW}_k^2=\mathcal{JW}_k\quad\text{and}\quad  \xy
(0,0)*{
\begin{tikzpicture}[scale=.3]
	\draw [very thick] (-1,-2.75) to [out=90,in=180] (0,-0.75) to [out=0,in=90] (1,-2.75);
	\draw [very thick] (-3,-2.75) to (-3,2.5);
	\draw [very thick] (-5,-2.75) to (-5,2.5);
	\draw [very thick] (3,-2.75) to (3,2.5);
	\draw [very thick] (5,-2.75) to (5,2.5);
	\node at (-1,-3.15) {\tiny $1$};
	\node at (1,-3.15) {\tiny $1$};
	\node at (-3,-3.15) {\tiny $1$};
	\node at (3,-3.15) {\tiny $1$};
	\node at (-5,-3.15) {\tiny $1$};
	\node at (5,-3.15) {\tiny $1$};
	\node at (-3,3) {\tiny $1$};
	\node at (3,3) {\tiny $1$};
	\node at (-5,3) {\tiny $1$};
	\node at (5,3) {\tiny $1$};
	\node at (-4,0) {\tiny $\dots$};
	\node at (4,0) {\tiny $\dots$};
\end{tikzpicture}
};
\endxy\circ\mathcal{JW}_k=0=\mathcal{JW}_k\circ\xy
(0,0)*{
\begin{tikzpicture}[scale=.3]
	\draw [very thick] (-1,2.5) to [out=270,in=180] (0,0.5) to [out=0,in=270] (1,2.5);
	\draw [very thick] (-3,-2.75) to (-3,2.5);
	\draw [very thick] (-5,-2.75) to (-5,2.5);
	\draw [very thick] (3,-2.75) to (3,2.5);
	\draw [very thick] (5,-2.75) to (5,2.5);
	\node at (-1,3) {\tiny $1$};
	\node at (1,3) {\tiny $1$};
	\node at (-3,-3.15) {\tiny $1$};
	\node at (3,-3.15) {\tiny $1$};
	\node at (-5,-3.15) {\tiny $1$};
	\node at (5,-3.15) {\tiny $1$};
	\node at (-3,3) {\tiny $1$};
	\node at (3,3) {\tiny $1$};
	\node at (-5,3) {\tiny $1$};
	\node at (5,3) {\tiny $1$};
	\node at (-4,0) {\tiny $\dots$};
	\node at (4,0) {\tiny $\dots$};
\end{tikzpicture}
};
\endxy
\]
Thus, $\mathcal{JW}_k$ is an idempotent which is killed 
by all possible cap compositions from the top and all 
possible cup compositions from the bottom.
\end{corollary}

\begin{proof}
Since $\cal{I}$ is a functor and $JW_k$ are idempotents which annihilate caps and cups, this is an immediate consequence of the previous result.
\end{proof}

\subsection{A diagrammatic description of \texorpdfstring{$\Repp$}{sl2-fdMod}}\label{sub-diahowe}

In this subsection, we prove Theorem~\ref{thm-main}. To do so, we must first deduce the existence of a 
functor $\Gam\colon\SymSp \to\Repp$, and then show that $\Gam$ induces the desired equivalence of categories. 
The definition of $\Gam$ is essentially dictated by our desire to have a commutative diagram
\begin{equation}\label{eq-fac}
\begin{gathered}
\xymatrix{
\Ud(\glm) \ar[rr]^{\Phi_m} \ar[dr]_{\Upsilon_m} & & \Repp \\
& \hspace*{0.2cm}\SymSp \ar[ur]_{\Gam} &
}
\end{gathered}
\end{equation}
We will begin by defining the functor $\Upsilon_m$.

\begin{lemma}\label{lem-ladderfunc}
For each $m\in\Z_{\geq 0}$, there exists a functor $\Upsilon_m\colon\Ud(\glm) \to \SymSp$ which sends a $\glm$-weight $\vec{k}\in\Z_{\geq 0}^m$
to the sequence obtained by removing all $0$'s
and all other objects $\vec{k}$ of $\Ud(\glm)$ to the zero object. This functor is 
determined on morphisms by the assignment
\begin{equation}\label{eq-UpOnLadder}
\Upsilon_m(F_i^{(j)} \one_{\vec{k}}) = 
\xy
(0,0)*{
\begin{tikzpicture}[scale=.3]
	\draw [very thick] (-4,-2) to (-4,2);
	\draw [very thick] (-1.5,-2) to (-1.5,0);
	\draw [very thick] (-1.5,0) to (-1.5,2);
	\draw [very thick] (1.5,-2) to (1.5,0);
	\draw [very thick] (1.5,0) to (1.5,2);
	\draw [very thick, directed=.55] (-1.5,0) to (1.5,0);
	\draw [very thick] (4,-2) to (4,2);
	\node at (-4,-2.5) {\tiny $k_1$};
	\node at (-1.5,-2.5) {\tiny $k_i$};
	\node at (1.5,-2.5) {\tiny $k_{i+1}$};
	\node at (4,-2.5) {\tiny $k_m$};
	\node at (-1.65,2.5) {\tiny $k_i{-}j$};
	\node at (1.65,2.5) {\tiny $k_{i+1}{+}j$};
	\node at (-4,2.5) {\tiny $k_1$};
	\node at (4,2.5) {\tiny $k_m$};
	\node at (0,0.75) {\tiny $j$};
	\node at (-2.75,0) {$\cdots$};
	\node at (2.75,0) {$\cdots$};
\end{tikzpicture}
};
\endxy 
\quad , \quad
\Upsilon_m(E_i^{(j)} \one_{\vec{k}}) = 
\xy
(0,0)*{
\begin{tikzpicture}[scale=.3]
	\draw [very thick] (-4,-2) to (-4,2);
	\draw [very thick] (-1.5,-2) to (-1.5,0);
	\draw [very thick] (-1.5,0) to (-1.5,2);
	\draw [very thick] (1.5,-2) to (1.5,0);
	\draw [very thick] (1.5,0) to (1.5,2);
	\draw [very thick, directed=.55] (1.5,0) to (-1.5,0);
	\draw [very thick] (4,-2) to (4,2);
	\node at (-4,-2.5) {\tiny $k_1$};
	\node at (-1.5,-2.5) {\tiny $k_i$};
	\node at (1.5,-2.5) {\tiny $k_{i+1}$};
	\node at (4,-2.5) {\tiny $k_m$};
	\node at (-1.65,2.5) {\tiny $k_i{+}j$};
	\node at (1.65,2.5) {\tiny $k_{i+1}{-}j$};
	\node at (-4,2.5) {\tiny $k_1$};
	\node at (4,2.5) {\tiny $k_m$};
	\node at (0,0.75) {\tiny $j$};
	\node at (-2.75,0) {$\cdots$};
	\node at (2.75,0) {$\cdots$};
\end{tikzpicture}
};
\endxy 
\end{equation}
where we erase any zero labeled edges in the diagrams depicting the images.
\end{lemma}

\begin{proof}
A straightforward check, using arguments found in Lemma 2.2.1 and Proposition 5.2.1 of~\cite{ckm}, shows that the images of relations in $\Ud(\glm)$ are 
consequences of the standard $\sln$-web relations in equations~\eqref{eq-frob},~\eqref{eq-simpler1}, 
and~\eqref{eq-almostsimpler}.
\end{proof}

We now aim to define the functor $\Gam$. We will first define the images of the generating morphisms in $\SymSp$, i.e. define a 
functor from the free symmetric spider $\SymSpf$, and then check that the relations in $\SymSp$ are satisfied.
Given a sequence $\ii=(i_1,\dots, i_m)$ with entries in $\{1,2\}$, 
we write $x^{\ii}$ as shorthand for
\[
x_{i_1}\otimes\cdots\otimes x_{i_m}\in (\C_q^2)^{\otimes m}.
\]
Furthermore, using Lemma 2.32 in~\cite{bz1}, we now fix a basis of $\Sym_q^k\C_q^2$ for all $k$ as the one given by the equivalence classes in $\Sym_q^k\C_q^2$
of all $x^{\ii}$ such that $\ii$ is weakly increasing and of length $k$. We will use the notation $x_{\ii}$ to denote the class of such an 
element in $\Sym_q^k\C_q^2$.

\begin{definition}\label{defn-mainfunctor}
Define a functor $\Gam\colon\SymSpf\to\Repp$ as follows.
\begin{itemize}
\item On objects: the tuples $\vec{k}=(k_1,\dots,k_m)\in\Z_{>0}^m$ are sent 
to the $\Uq(\slnn{2})$-modules 
$\Sym_q^{k_1}\C_q^2\otimes\cdots\otimes\Sym_q^{k_m}\C_q^2$. Moreover, we 
send, by convention, the empty tuple to the trivial $\Uq(\slnn{2})$-module $\C_q$ 
and the zero object to the zero module.
\item On morphisms: we send the generators of $\SymSpf$ to the following 
$\Uq(\slnn{2})$-intertwiners, and extend monoidally. 
We send the thickness $k$ identity strand to 
$\mathrm{id}_k\colon \Sym_q^{k}\C_q^2\to \Sym_q^{k}\C_q^2$, 
and define the functor on $1$-labeled caps and cups via
\[
\Gam \left(
\xy
(0,0)*{
\begin{tikzpicture} [scale=1]
\draw[very thick] (0,0) to [out=90,in=180] (.25,.5) to [out=0,in=90] (.5,0);
\node at (.25,.75) {};
\node at (0,-.15) {\tiny $1$};
\node at (.5,-.15) {\tiny $1$};
\end{tikzpicture}
}
\endxy
\right)
=
\mathsf{cap} \colon\C_q^2\otimes \C_q^2 \twoheadrightarrow \C_q,\quad\begin{cases}
x^{\ul{11}},x^{\ul{22}}\mapsto 0,\\
x^{\ul{12}}\mapsto -q,\\
x^{\ul{21}}\mapsto 1,
\end{cases}
\]
and
\[
\Gam \left(
\xy
(0,0)*{
\begin{tikzpicture} [scale=1]
\draw[very thick] (0,1) to [out=270,in=180] (.25,.5) to [out=0,in=270] (.5,1);
\node at (.25,.25) {};
\node at (0,1.15) {\tiny $1$};
\node at (.5,1.15) {\tiny $1$};
\end{tikzpicture}
}
\endxy
\right)
=
\mathsf{cup} \colon\C_q \hookrightarrow \C_q^2\otimes \C_q^2,\quad 1\mapsto x^{\ul{12}}-q^{-1} x^{\ul{21}}.
\]
On merge and split generators, we define $\Gam$ using the functor $\Phi_2$ from Corollary~\ref{cor-functor}, that is,
\[
\Gam \left(
\xy
(0,0)*{
\begin{tikzpicture}[scale=.3]
	\draw [very thick] (0, .75) to (0,2.5);
	\draw [very thick] (1,-1) to [out=90,in=330] (0,.75);
	\draw [very thick] (-1,-1) to [out=90,in=210] (0,.75); 
	\node at (0, 3) {\tiny $k{+}l$};
	\node at (-1,-1.5) {\tiny $k$};
	\node at (1,-1.5) {\tiny $l$};
\end{tikzpicture}
};
\endxy
\right)
=
\Phi_2(E^{(l)}\one_{(k,l)}), \quad
\Gam \left(
\xy
(0,0)*{
\begin{tikzpicture}[scale=.3]
	\draw [very thick] (0,-1) to (0,.75);
	\draw [very thick] (0,.75) to [out=30,in=270] (1,2.5);
	\draw [very thick] (0,.75) to [out=150,in=270] (-1,2.5); 
	\node at (0, -1.5) {\tiny $k{+}l$};
	\node at (-1,3) {\tiny $k$};
	\node at (1,3) {\tiny $l$};
\end{tikzpicture}
};
\endxy
\right)
=
\Phi_2(F^{(l)}\one_{(k+l,0)}).
\]
Having defined $\Gam$ on these generators, we can extend to $k$-labeled caps via the assignment
\[
\Gam \left(
\xy
(0,0)*{
\begin{tikzpicture} [scale=1]
\draw[very thick] (0,0) to [out=90,in=180] (.25,.5) to [out=0,in=90] (.5,0);
\node at (.25,.75) {};
\node at (0,-.15) {\tiny $k$};
\node at (.5,-.15) {\tiny $k$};
\end{tikzpicture}
}
\endxy
\right)
=
\frac{1}{[k]!}
\Gam \left(
\xy
(0,0)*{
\begin{tikzpicture}[scale=.3]
	\draw [very thick] (0,-1) to (0,.75);
	\draw [very thick] (0,.75) to [out=30,in=270] (1,2.5);
	\draw [very thick] (0,.75) to [out=150,in=270] (-1,2.5); 
	\draw [very thick] (1,2.5) to [out=90,in=180] (2,3.5) to [out=0,in=90] (3,2.5);
	\draw [very thick] (-1,2.5) to [out=90,in=180] (2,5.5) to [out=0,in=90] (5,2.5);
	\draw [very thick] (4,-1) to (4,.75);
	\draw [very thick] (4,-1) to (4,.75);
	\draw [very thick] (4,.75) to [out=30,in=270] (5,2.5);
	\draw [very thick] (4,.75) to [out=150,in=270] (3,2.5);
	\node at (0,-1.35) {\tiny $k$};
	\node at (4,-1.35) {\tiny $k$};
	\node at (2,3) {\tiny $1$};
	\node at (2,6) {\tiny $1$};
	\node at (2,4.75) {\tiny $\vdots$};
	\node at (0,2.25) {\tiny $\cdots$};
	\node at (4,2.25) {\tiny $\cdots$};
\end{tikzpicture}
};
\endxy
\right)
\]
and similarly for $k$-labeled cups.\qedhere
\end{itemize}
\end{definition}

We will denote the images under $\Gam$ of $1$-labeled caps and cups (as above) by $\mathsf{cap}$ and $\mathsf{cup}$, and 
the images of the symmetric $(1,1)$-merge and $(1,1)$-split $\slnn{2}$-webs by $\mathsf{m}$ and $\mathsf{s}$. 
Moreover, for thickened versions we use the notation $\mathsf{cap}_k$, $\mathsf{cup}_k$, 
$\mathsf{m}_{k,l}$ and $\mathsf{s}^{k,l}$ in the evident way.

\begin{remark}\label{rem-ambiguity}
The meticulous reader will note that there is an ambiguity in our 
definition of caps and cups of thickness $k$, in that 
we did not choose a particular choice for the symmetric $\slnn{2}$-web which splits a 
$k$-labeled strand into $k$ strands of thickness $1$. 
Indeed, it follows from the associativity relations~\eqref{eq-frob} in $\Ud(\glm)$ that the corresponding morphisms in $\Repp$ are the same. 
The concerned reader can use their favorite such symmetric $\slnn{2}$-web as the one used in the above definition.

The reader may also be curious about our choices in the definition of $\Gam$ on merge and split morphisms, 
i.e. why not set
\[
\Gam \left(
\xy
(0,0)*{
\begin{tikzpicture}[scale=.3]
	\draw [very thick] (0, .75) to (0,2.5);
	\draw [very thick] (1,-1) to [out=90,in=330] (0,.75);
	\draw [very thick] (-1,-1) to [out=90,in=210] (0,.75); 
	\node at (0, 3) {\tiny $k{+}l$};
	\node at (-1,-1.5) {\tiny $k$};
	\node at (1,-1.5) {\tiny $l$};
\end{tikzpicture}
};
\endxy
\right)
=
\Phi_2(F^{(k)}\one_{(k,l)}), \quad
\Gam \left(
\xy
(0,0)*{
\begin{tikzpicture}[scale=.3]
	\draw [very thick] (0,-1) to (0,.75);
	\draw [very thick] (0,.75) to [out=30,in=270] (1,2.5);
	\draw [very thick] (0,.75) to [out=150,in=270] (-1,2.5); 
	\node at (0, -1.5) {\tiny $k{+}l$};
	\node at (-1,3) {\tiny $k$};
	\node at (1,3) {\tiny $l$};
\end{tikzpicture}
};
\endxy
\right)
=
\Phi_2(E^{(k)}\one_{(0,k+l)})
\quad ?
\]
Indeed, this will lead to the same definition, following from the equalities
\[
E^{(k+l)}F^{(k)}\one_{(k,l)} = E^{(l)}\one_{(k,l)}, \quad
E^{(k)} F^{(k+l)}\one_{(k+l,0)} = F^{(l)}\one_{(k+l,0)}
\]
in $\Ud^{\infty}(\glnn{2})$ and the fact that $\Phi_2(E^{(k)}_{(0,k)})$ and $\Phi_2(F^{(k)}_{(k,0)})$ are 
both the identity morphism of $\Sym_q^k\C_q^2$.
\end{remark}

\begin{example}
Since we will need these explicitly later, we now record the $(1,1)$-merge and the $(1,1)$-split morphisms. They are given by
\[
\Gam \left(
\xy
(0,0)*{
\begin{tikzpicture}[scale=.3]
	\draw [very thick] (0, .75) to (0,2.5);
	\draw [very thick] (1,-1) to [out=90,in=330] (0,.75);
	\draw [very thick] (-1,-1) to [out=90,in=210] (0,.75); 
	\node at (0, 3) {\tiny $2$};
	\node at (-1,-1.5) {\tiny $1$};
	\node at (1,-1.5) {\tiny $1$};
\end{tikzpicture}
};
\endxy
\right)
=
\mathsf{m}\colon \C_q^2\otimes \C_q^2 \twoheadrightarrow \Sym_q^2\C_q^2,\quad
\begin{cases}
x^{\ul{11}}\mapsto x_{\ul{11}},\; x^{\ul{12}}\mapsto x_{\ul{12}},\\
x^{\ul{21}}\mapsto qx_{\ul{12}},\; x^{\ul{22}}\mapsto x_{\ul{22}},
\end{cases}
\]
and
\[
\Gam \left(
\xy
(0,0)*{
\begin{tikzpicture}[scale=.3]
	\draw [very thick] (0,-1) to (0,.75);
	\draw [very thick] (0,.75) to [out=30,in=270] (1,2.5);
	\draw [very thick] (0,.75) to [out=150,in=270] (-1,2.5); 
	\node at (0, -1.5) {\tiny $2$};
	\node at (-1,3) {\tiny $1$};
	\node at (1,3) {\tiny $1$};
\end{tikzpicture}
};
\endxy
\right)
=
\mathsf{s}\colon\Sym_q^2\C_q^2 \hookrightarrow \C_q^2\otimes \C_q^2,\quad
\begin{cases}
x_{\ul{11}}\mapsto [2]x^{\ul{11}},\; x_{\ul{22}}\mapsto [2]x^{\ul{22}},\\
x_{\ul{12}}\mapsto q^{-1}x^{\ul{12}}+x^{\ul{21}}.
\end{cases}
\]
Moreover, the $2$-labeled cap is given by
\begin{equation} \label{eq-cap2}
\Gam \left(
\xy
(0,0)*{
\begin{tikzpicture} [scale=1]
\draw[very thick] (0,0) to [out=90,in=180] (.25,.5) to [out=0,in=90] (.5,0);
\node at (.25,.75) {};
\node at (0,-.15) {\tiny $2$};
\node at (.5,-.15) {\tiny $2$};
\end{tikzpicture}
}
\endxy
\right)
=\mathsf{cap}_2\colon \Sym_q^2\C_q^2\otimes \Sym_q^2\C_q^2 \twoheadrightarrow \C_q,\quad
\begin{cases}
x_{\ul{11}}\otimes x_{\ul{22}}\mapsto q^2[2],\; x_{\ul{12}}\otimes x_{\ul{22}}\mapsto -1,\\
x_{\ul{22}}\otimes x_{\ul{11}}\mapsto [2],\; \text{rest}\mapsto 0.
\end{cases}
\end{equation}
We encourage the reader to work out $\mathsf{cup}_2$, which we will use in algebraic form below as well.
\end{example}

\begin{lemma}\label{lem-mainfunctor}
$\Gam$ descends to give a monoidal functor $\Gam\colon\SymSp\to\Repp$.
\end{lemma}

\begin{proof}
It is clear that, if $\Gam$ is well-defined, then it also preserves the monoidal structure 
(which is given by placing diagrams next to each other).
To check that $\Gam$ is well-defined, it suffices 
to show that the relations of the symmetric $\slnn{2}$-spider $\SymSp$ hold in $\Repp$.

The ``standard'' $\sln$-web relations -- associativity~\eqref{eq-frob}, 
digon and square removals~\eqref{eq-simpler1} and the square switches~\eqref{eq-almostsimpler} -- 
follow from Corollary~\ref{cor-functor}, since these are all induced by relations in $\Ud(\glm)$. 
Here we have to utilize the property that the 
images of the divided powers $F_i^{(j)}\one_{\vec{k}}$ and $E_i^{(j)}\one_{\vec{k}}$ under 
$\Phi_m\colon\Ud(\glm) \to \Repp$ coincide with the images of 
the general symmetric ladders in 
equation~\eqref{eq-UpOnLadder}
under 
$\Gam\colon\SymSp\to\Repp$.
This follows from our definition of $\Gam$ on symmetric merge and split $\slnn{2}$-webs, 
the $\Ud(\glnn{3})$ equalities
\[
E_1^{(j)}\one_{(k,l,0)} = E_2^{(l-j)} E_1^{(j)} F_2^{(l-j)}\one_{(k,l,0)}
\quad , \quad
F_1^{(j)}\one_{(k,l,0)} = E_2^{(l)} F_1^{(j)} F_2^{(l)}\one_{(k,l,0)} 
\]
and the fact that the diagram
\[
\xymatrix{
\Ud(\glm) \ar@{^{(}->}[r] \ar[rd]_{\Phi_m} & \Ud(\glnn{m+1}) \ar[d]^{\Phi_{m+1}} \\
& \Repp}
\]
commutes for any of the standard inclusions $\Ud(\glm) \hookrightarrow \Ud(\glnn{m+1})$.

It now remains to check that the additional symmetric and isotopy relations are satisfied. \\

\noindent \underline{Circle removal.} This relation follows from the computation
\[
(\mathsf{cap}\circ\mathsf{cup})(1)=\mathsf{cap}(x^{\ul{12}})-q^{-1}\mathsf{cap}(x^{\ul{21}})=-q-q^{-1}=-[2],
\]
where we point out the negative sign to the reader. As is known to experts, this is unavoidable if one 
wishes to have isotopy invariance in an unoriented model. \\

\noindent \underline{Dumbbell relation.} This can again be directly verified. For example, we have
\[
(\mathsf{s}\circ\mathsf{m})(x^{\ul{21}})=x^{\ul{12}}+qx^{\ul{21}}
\]
and 
\[
([2]\mathsf{id}+\mathsf{cup}\circ\mathsf{cap})(x^{\ul{21}})=[2]x^{\ul{21}}+x^{\ul{12}}-q^{-1}x^{\ul{21}}=x^{\ul{12}}+qx^{\ul{21}}.
\]
The remainder of the computations follow similarly. \\

\noindent \underline{Isotopy relations.}
The remaining isotopy relations locally reduce to the following relations:
\begin{equation}\label{eq-snake}
\xy
(0,0)*{\reflectbox{
\begin{tikzpicture} [scale=1]
\draw[very thick] (0,0) to [out=90,in=180] (.25,.5) to [out=0,in=90] (.5,0);
\draw[very thick] (-.5,0) to [out=270,in=180] (-.25,-.5) to [out=0,in=270] (0,0);
\draw[very thick] (0.5,0) to (0.5,-.75);
\draw[very thick] (-.5,0) to (-.5,.75);
\node at (0.5,-.9) {\reflectbox{\tiny $k$}};
\node at (-0.5,.9) {\reflectbox{\tiny $k$}};
\end{tikzpicture}
}};
\endxy=
\xy
(0,0)*{
\begin{tikzpicture} [scale=1]
\draw[very thick] (0,0) to (0,1.5);
\node at (0,-.15) {\tiny $k$};
\node at (0,1.65) {\tiny $k$};
\end{tikzpicture}
}
\endxy=\xy
(0,0)*{
\begin{tikzpicture} [scale=1]
\draw[very thick] (0,0) to [out=90,in=180] (.25,.5) to [out=0,in=90] (.5,0);
\draw[very thick] (-.5,0) to [out=270,in=180] (-.25,-.5) to [out=0,in=270] (0,0);
\draw[very thick] (0.5,0) to (0.5,-.75);
\draw[very thick] (-.5,0) to (-.5,.75);
\node at (0.5,-.9) {\tiny $k$};
\node at (-0.5,.9) {\tiny $k$};
\end{tikzpicture}
};
\endxy 
\end{equation}
and
\begin{equation}\label{eq-vertexslide}
\xy
(0,0)*{
\begin{tikzpicture}[scale=.3]
	\draw [very thick] (0,-1) to (0,.75);
	\draw [very thick] (0,.75) to [out=30,in=270] (1,2.5);
	\draw [very thick] (0,.75) to [out=150,in=270] (-1,2.5); 
	\draw [very thick] (1,2.5) to [out=90,in=180] (2,3.5) to [out=0,in=90] (3,2.5);
	\draw [very thick] (-1,2.5) to [out=90,in=180] (2,4.5) to [out=0,in=90] (5,2.5);
	\draw [very thick] (3,-1) to (3,2.5);
	\draw [very thick] (5,-1) to (5,2.5);
	\node at (0,-1.35) {\tiny $k{+}l$};
	\node at (3,-1.35) {\tiny $k$};
	\node at (5,-1.35) {\tiny $l$};
\end{tikzpicture}
};
\endxy=
\xy
(0,0)*{
\begin{tikzpicture}[scale=.3]
	\draw [very thick] (4, .75) to (4,2.5);
	\draw [very thick] (5,-1) to [out=90,in=330] (4,.75);
	\draw [very thick] (3,-1) to [out=90,in=210] (4,.75);
	\draw [very thick] (0, -1) to (0,2.5);
	\draw [very thick] (0,2.5) to [out=90,in=180] (2,4.5) to [out=0,in=90] (4,2.5); 
	\node at (0,-1.35) {\tiny $k{+}l$};
	\node at (3,-1.35) {\tiny $k$};
	\node at (5,-1.35) {\tiny $l$};
\end{tikzpicture}
};
\endxy\quad\text{and}\quad
\xy
(0,0)*{\reflectbox{
\begin{tikzpicture}[scale=.3]
	\draw [very thick] (0,-1) to (0,.75);
	\draw [very thick] (0,.75) to [out=30,in=270] (1,2.5);
	\draw [very thick] (0,.75) to [out=150,in=270] (-1,2.5); 
	\draw [very thick] (1,2.5) to [out=90,in=180] (2,3.5) to [out=0,in=90] (3,2.5);
	\draw [very thick] (-1,2.5) to [out=90,in=180] (2,4.5) to [out=0,in=90] (5,2.5);
	\draw [very thick] (3,-1) to (3,2.5);
	\draw [very thick] (5,-1) to (5,2.5);
	\node at (0,-1.35) {\reflectbox{\tiny $k{+}l$}};
	\node at (3,-1.35) {\reflectbox{\tiny $l$}};
	\node at (5,-1.35) {\reflectbox{\tiny $k$}};
\end{tikzpicture}}
};
\endxy=
\xy
(0,0)*{\reflectbox{
\begin{tikzpicture}[scale=.3]
	\draw [very thick] (4, .75) to (4,2.5);
	\draw [very thick] (5,-1) to [out=90,in=330] (4,.75);
	\draw [very thick] (3,-1) to [out=90,in=210] (4,.75);
	\draw [very thick] (0, -1) to (0,2.5);
	\draw [very thick] (0,2.5) to [out=90,in=180] (2,4.5) to [out=0,in=90] (4,2.5); 
	\node at (0,-1.35) {\reflectbox{\tiny $k{+}l$}};
	\node at (3,-1.35) {\reflectbox{\tiny $l$}};
	\node at (5,-1.35) {\reflectbox{\tiny $k$}};
\end{tikzpicture}
}};
\endxy
\end{equation}
and versions of the ones from~\eqref{eq-vertexslide} involving cups.

We start with~\eqref{eq-vertexslide}, first noting that it suffices to verify the case where either $k=1$ or $l=1$. 
Indeed, assuming that the relation is known in these cases, we can repeatedly use the first relation in
equation~\eqref{eq-simpler1}
to explode the $k$- and $l$-labeled strands into 1-labeled strands. 
We can then pull each of the merge and split vertices (which all have at least one 1-labeled strand) around the cap, and 
use~\eqref{eq-simpler1} to reassemble the exploded strand.

For the remaining cases, we can then use a similar argument to verify the relation (doubly) inductively. 
All together, we see that it suffices to prove the relations explicitly when $k=1=l$.
We hence, compute that the lower part of the left-hand side of the first equation in~\eqref{eq-vertexslide} is given by
\[
\mathsf{s}\otimes\mathrm{id}\otimes\mathrm{id}\colon \Sym_q^2\C_q^2\otimes \C_q^2\otimes \C_q^2\to \C_q^2\otimes \C_q^2\otimes \C_q^2\otimes \C_q^2,\quad\begin{cases}
x_{\ul{11}}\otimes x^{\ul{ij}}\mapsto [2]x^{\ul{11ij}},\\
x_{\ul{12}}\otimes x^{\ul{ij}}\mapsto q^{-1}x^{\ul{12ij}}+x^{\ul{21ij}},\\
x_{\ul{22}}\otimes x^{\ul{ij}}\mapsto [2]x^{\ul{22ij}},
\end{cases}
\]
for all choices of $i,j \in \{1,2\}$.
Most of these terms will be sent to zero after composing with the top, and the only surviving terms are
\[
\mathsf{cap}\circ(\mathrm{id}\otimes\mathsf{cap}\otimes\mathrm{id})\circ(\mathsf{s}\otimes\mathrm{id}
\otimes\mathrm{id})\colon \Sym_q^2\C_q^2\otimes \C_q^2\otimes \C_q^2\to \C_q,\quad\begin{cases}
x_{\ul{11}}\otimes x^{\ul{22}}\mapsto q^2[2],\\
x_{\ul{12}}\otimes x^{\ul{12}}\mapsto -1,\\
x_{\ul{12}}\otimes x^{\ul{21}}\mapsto -q,\\
x_{\ul{22}}\otimes x^{\ul{11}}\mapsto [2].
\end{cases}
\]
The bottom part of the right-hand side is given (for all $1\leq i\leq j\leq 2$) by
\[
\mathrm{id}\otimes\mathsf{m}\colon \Sym_q^2\C_q^2\otimes \C_q^2\otimes \C_q^2\to \Sym_q^2\C_q^2\otimes \Sym_q^2\C_q^2,\quad\begin{cases}
x_{\ul{ij}}\otimes x^{\ul{11}}\mapsto x_{\ul{ij}}\otimes x_{\ul{11}},\\
x_{\ul{ij}}\otimes x^{\ul{12}}\mapsto x_{\ul{ij}}\otimes x_{\ul{12}},\\
x_{\ul{ij}}\otimes x^{\ul{21}}\mapsto q x_{\ul{ij}}\otimes x_{\ul{12}},\\
x_{\ul{ij}}\otimes x^{\ul{22}}\mapsto x_{\ul{ij}}\otimes x_{\ul{22}},
\end{cases}
\]
which composes with the map in
equation~\eqref{eq-cap2}
to give the correct result. 
The check of the second equation in~\eqref{eq-vertexslide} is similar, as are the checks of the 
versions of this relation involving cups.

We can deduce the general form of~\eqref{eq-snake} from the $k=1$ case and 
the relations in~\eqref{eq-vertexslide} (and their analogs) using the following diagrammatic argument:
\[
\xy
(0,0)*{\reflectbox{
\begin{tikzpicture} [scale=1]
\draw[very thick] (0,0) to [out=90,in=180] (.25,.5) to [out=0,in=90] (.5,0);
\draw[very thick] (-.5,0) to [out=270,in=180] (-.25,-.5) to [out=0,in=270] (0,0);
\draw[very thick] (0.5,0) to (0.5,-.75);
\draw[very thick] (-.5,0) to (-.5,.75);
\node at (0.5,-.9) {\reflectbox{\tiny $k$}};
\node at (-0.5,.9) {\reflectbox{\tiny $k$}};
\end{tikzpicture}
}};
\endxy
=\frac{1}{[k]!}
\xy
(0,0)*{\reflectbox{
\begin{tikzpicture} [scale=1]
\draw[very thick] (0,0) to [out=90,in=180] (.25,.5) to [out=0,in=90] (.5,0);
\draw[very thick] (-.5,0) to [out=270,in=180] (-.25,-.5) to [out=0,in=270] (0,0);
\draw[very thick] (0.5,0) to (0.5,-.1);
\draw[very thick] (0.5,-.65) to (0.5,-.75);
\draw[very thick] (-.5,0) to (-.5,.75);
\draw[very thick] (0.35,-.375) to [out=270,in=150] (0.5,-.65);
\draw[very thick] (0.65,-.375) to [out=270,in=30] (0.5,-.65);
\draw[very thick] (0.35,-.375) to [out=90,in=210] (0.5,-.1);
\draw[very thick] (0.65,-.375) to [out=90,in=330] (0.5,-.1);
\node at (0.5,-.375) {\reflectbox{\tiny ...}};
\node at (0.5,-.9) {\reflectbox{\tiny $k$}};
\node at (-0.5,.9) {\reflectbox{\tiny $k$}};
\end{tikzpicture}
}};
\endxy
=
\frac{1}{[k]!}
\xy
(0,0)*{\reflectbox{
\begin{tikzpicture} [scale=1]
\draw[very thick] (-0.65,.375) to (-0.65,0);
\draw[very thick] (0.5,-.65) to (0.5,-.75);
\draw[very thick] (0.65,-.375) to (0.65,0);
\draw[very thick] (0.5,-.65) to (0.5,-.75);
\draw[very thick] (-0.65,0) to (-0.65,-0.18775) to [out=270,in=180] (-0.2,-0.5) to [out=0,in=270] (0.1,-0.18775) to [out=90,in=180] (0.2,0) to [out=0,in=90] (0.35,-.375);
\draw[very thick] (-0.35,.375) to [out=270,in=180] (-0.2,0) to [out=0,in=270] (-.1,0.18775) to [out=90,in=180] (0.2,0.5) to [out=0,in=90] (0.65,0.18775) to (0.65,0);
\draw[very thick] (0.35,-.375) to [out=270,in=150] (0.5,-.65);
\draw[very thick] (0.65,-.375) to [out=270,in=30] (0.5,-.65);
\draw[very thick] (-0.5,.65) to (-0.5,.75);
\draw[very thick] (-0.35,.375) to [out=90,in=330] (-0.5,.65);
\draw[very thick] (-0.65,.375) to [out=90,in=210] (-0.5,.65);
\node at (-0.5,.375) {\reflectbox{\tiny ...}};
\node at (0.5,-.375) {\reflectbox{\tiny ...}};
\node at (-0.25,-.275) {\reflectbox{\tiny ...}};
\node at (0.25,.275) {\reflectbox{\tiny ...}};
\node at (0.5,-.9) {\reflectbox{\tiny $k$}};
\node at (-0.5,.9) {\reflectbox{\tiny $k$}};
\end{tikzpicture}
}};
\endxy
=\frac{1}{[k]!}
\xy
(0,0)*{
\begin{tikzpicture} [scale=1]
\draw[very thick] (0,-.65) to (0,-.75);
\draw[very thick] (0,.65) to (0,.75);
\draw[very thick] (-0.15,-.375) to [out=270,in=150] (0,-.65);
\draw[very thick] (0.15,-.375) to [out=270,in=30] (0,-.65);
\draw[very thick] (0.15,.375) to [out=90,in=330] (0,.65);
\draw[very thick] (-0.15,.375) to [out=90,in=210] (0,.65);
\draw[very thick] (0.15,.375) to (0.15,-.375);
\draw[very thick] (-0.15,.375) to (-0.15,-.375);
\node at (0,0) {\tiny ...};
\node at (0,-.9) {\tiny $k$};
\node at (0,.9) {\tiny $k$};
\end{tikzpicture}
};
\endxy
=\frac{1}{[k]!}
\xy
(0,0)*{
\begin{tikzpicture} [scale=1]
\draw[very thick] (-0.65,.375) to (-0.65,0);
\draw[very thick] (0.5,-.65) to (0.5,-.75);
\draw[very thick] (0.65,-.375) to (0.65,0);
\draw[very thick] (0.5,-.65) to (0.5,-.75);
\draw[very thick] (-0.65,0) to (-0.65,-0.18775) to [out=270,in=180] (-0.2,-0.5) to [out=0,in=270] (0.1,-0.18775) to [out=90,in=180] (0.2,0) to [out=0,in=90] (0.35,-.375);
\draw[very thick] (-0.35,.375) to [out=270,in=180] (-0.2,0) to [out=0,in=270] (-.1,0.18775) to [out=90,in=180] (0.2,0.5) to [out=0,in=90] (0.65,0.18775) to (0.65,0);
\draw[very thick] (0.35,-.375) to [out=270,in=150] (0.5,-.65);
\draw[very thick] (0.65,-.375) to [out=270,in=30] (0.5,-.65);
\draw[very thick] (-0.5,.65) to (-0.5,.75);
\draw[very thick] (-0.35,.375) to [out=90,in=330] (-0.5,.65);
\draw[very thick] (-0.65,.375) to [out=90,in=210] (-0.5,.65);
\node at (-0.5,.375) {\tiny ...};
\node at (0.5,-.375) {\tiny ...};
\node at (-0.25,-.275) {\tiny ...};
\node at (0.25,.275) {\tiny ...};
\node at (0.5,-.9) {\tiny $k$};
\node at (-0.5,.9) {\tiny $k$};
\end{tikzpicture}
};
\endxy
=\frac{1}{[k]!}
\xy
(0,0)*{
\begin{tikzpicture} [scale=1]
\draw[very thick] (0,0) to [out=90,in=180] (.25,.5) to [out=0,in=90] (.5,0);
\draw[very thick] (-.5,0) to [out=270,in=180] (-.25,-.5) to [out=0,in=270] (0,0);
\draw[very thick] (0.5,0) to (0.5,-.1);
\draw[very thick] (0.5,-.65) to (0.5,-.75);
\draw[very thick] (-.5,0) to (-.5,.75);
\draw[very thick] (0.35,-.375) to [out=270,in=150] (0.5,-.65);
\draw[very thick] (0.65,-.375) to [out=270,in=30] (0.5,-.65);
\draw[very thick] (0.35,-.375) to [out=90,in=210] (0.5,-.1);
\draw[very thick] (0.65,-.375) to [out=90,in=330] (0.5,-.1);
\node at (0.5,-.375) {\tiny ...};
\node at (0.5,-.9) {\tiny $k$};
\node at (-0.5,.9) {\tiny $k$};
\end{tikzpicture}
};
\endxy
=
\xy
(0,0)*{
\begin{tikzpicture} [scale=1]
\draw[very thick] (0,0) to [out=90,in=180] (.25,.5) to [out=0,in=90] (.5,0);
\draw[very thick] (-.5,0) to [out=270,in=180] (-.25,-.5) to [out=0,in=270] (0,0);
\draw[very thick] (0.5,0) to (0.5,-.75);
\draw[very thick] (-.5,0) to (-.5,.75);
\node at (0.5,-.9) {\tiny $k$};
\node at (-0.5,.9) {\tiny $k$};
\end{tikzpicture}
};
\endxy
\]
Here the middle equalities follow from repeated application of the case 
$k=1$, and the diagram in the middle 
is, by digon removals, the identity.

The $k=1$ case follows by combining the computation
\[
(\mathrm{id}\otimes\mathsf{cup})(x_i)=x^{\ul{i12}}-q^{-1}x^{\ul{i21}}
\]
with
\[
(\mathsf{cap}\otimes\mathrm{id})(x^{\ul{i12}})=\begin{cases}
0, &\text{if } i=1,\\
+1, &\text{if } i=2,
\end{cases}\quad  \quad
(\mathsf{cap}\otimes\mathrm{id})(x^{\ul{i21}})=\begin{cases}
0, &\text{if } i=2,\\
-q, &\text{if } i=1,
\end{cases}
\]
for the left diagram and 
\[
(\mathsf{cup}\otimes\mathrm{id})(x_i)=x^{\ul{12i}}-q^{-1}x^{\ul{21i}}
\]
with
\[
(\mathrm{id}\otimes\mathsf{cap})(x^{\ul{12i}})=\begin{cases}
0, &\text{if } i=2,\\
+1, &\text{if } i=1,
\end{cases}\quad  \quad
(\mathrm{id}\otimes\mathsf{cap})(x^{\ul{21i}})=\begin{cases}
0, &\text{if } i=1,\\
-q, &\text{if } i=2,
\end{cases}
\]
for the right. We point out that the signs work out as they should.

Finally, all isotopies similar to
\[
\xy
(0,0)*{
\begin{tikzpicture}[scale=.3]
	\draw [very thick] (0, .75) to (0,2.5);
	\draw [very thick] (1,-1) to [out=90,in=330] (0,.75);
	\draw [very thick] (-1,-1) to [out=90,in=210] (0,.75); 
	\node at (0, 3) {\tiny $l$};
	\node at (-1,-1.5) {\tiny $k$};
	\node at (1,-1.5) {\tiny $k{+}l$};
\end{tikzpicture}
};
\endxy
=
\xy
(0,0)*{\reflectbox{
\begin{tikzpicture}[scale=.3]
	\draw [very thick] (0,-1) to (0,.75);
	\draw [very thick] (0,.75) to [out=30,in=270] (1,2.5);
	\draw [very thick] (0,.75) to [out=150,in=270] (-1,2.5);
	\draw [very thick] (1,2.5) to [out=90,in=180] (2,3.5) to [out=0,in=90] (3,2.5);
	\draw [very thick] (3,-1) to (3,2.5);
	\draw [very thick] (-1,2.5) to (-1,3.5);
	\node at (0, -1.5) {\reflectbox{\tiny $k{+}l$}};
	\node at (-1,4) {\reflectbox{\tiny $l$}};
	\node at (3,-1.5) {\reflectbox{\tiny $k$}};
\end{tikzpicture}}
};
\endxy
\quad\text{and}\quad
\xy
(0,0)*{
\begin{tikzpicture}[scale=.3]
	\draw [very thick] (0, .75) to (0,2.5);
	\draw [very thick] (1,-1) to [out=90,in=330] (0,.75);
	\draw [very thick] (-1,-1) to [out=90,in=210] (0,.75); 
	\node at (0, 3) {\tiny $k$};
	\node at (-1,-1.5) {\tiny $k{+}l$};
	\node at (1,-1.5) {\tiny $l$};
\end{tikzpicture}
};
\endxy
=
\xy
(0,0)*{
\begin{tikzpicture}[scale=.3]
	\draw [very thick] (0,-1) to (0,.75);
	\draw [very thick] (0,.75) to [out=30,in=270] (1,2.5);
	\draw [very thick] (0,.75) to [out=150,in=270] (-1,2.5);
	\draw [very thick] (1,2.5) to [out=90,in=180] (2,3.5) to [out=0,in=90] (3,2.5);
	\draw [very thick] (3,-1) to (3,2.5);
	\draw [very thick] (-1,2.5) to (-1,3.5);
	\node at (0, -1.5) {\tiny $k{+}l$};
	\node at (-1,4) {\tiny $k$};
	\node at (3,-1.5) {\tiny $l$};
\end{tikzpicture}
};
\endxy
\]
are not relations, but rather definitions of the elements on the left-hand sides. 
\end{proof}


As a consequence of this proof, we immediately observe the following.

\begin{corollary}\label{cor-lefttothereader}
The diagram from~\eqref{eq-fac} commutes.
\end{corollary}

\begin{remark}\label{rem-additive}
We can extend $\Gam$ additively to a functor
\[
\Gam\colon\Mat(\SymSp)\to\Repp
\]
that we, by abuse of notation, denote using the same symbol. Here $\Mat(\SymSp)$ is the additive closure of the symmetric $\slnn{2}$-spider. 
As we recalled above before Theorem~\ref{thm-main}, this means that objects of $\Mat(\SymSp)$ are finite, formal direct sums of the objects of $\SymSp$ 
and morphisms are matrices (whose entries are morphisms from $\SymSp$) between these sums.
Note that this category is again entirely diagrammatic.
\end{remark}

We are now ready to prove Theorem~\ref{thm-main}.

\begin{proof}[Proof (of Theorem~\ref{thm-main}).]
We have a well-defined functor $\Gam\colon\Mat(\SymSp)\to\Repp$ that preserves the monoidal structure. 
It only remains to show that $\Gam$ is essentially surjective and fully faithful.
\\

\noindent \underline{Essentially surjective.} This follows directly from the definition of the functor $\Gam$, since every finite-dimensional $\Uq(\slnn{2})$-module
is isomorphic to a direct sum of copies of $\Sym_q^k \C_q^2$.
\\

\noindent \underline{Fully faithful.} 
By additivity, we can verify everything on objects of the form $\vec{k}\in\Z_{>0}^m$.
Given $\vec{k}\in\Z^m_{> 0},\vec{l}\in\Z^{m^{\prime}}_{> 0}$, we have to show that
\begin{equation}\label{eq-toverify}
\Hom_{\SymSp}(\vec{k},\vec{l})\cong\Hom_{\Repp}(\Gam(\vec{k}),\Gam(\vec{l}))
\end{equation}
as $\C(q)$-vector spaces.
Surjectivity in~\eqref{eq-toverify} follows\footnote{We point out that this shows that all symmetric $\slnn{2}$-webs that contain cups and caps 
and whose top and bottom labels have the same sum
can be expressed as linear combinations 
of compositions of $F^{(j)}_i$ and $E^{(j)}_i$-ladders.} from Corollary~\ref{cor-functor} and Corollary~\ref{cor-lefttothereader}
in the case that $\sum_{i=1}^m k_i = \sum_{i=1}^{m^{\prime}} l_i$.
To see the case when $\sum_{i=1}^m k_i \neq \sum_{i=1}^{m^{\prime}} l_i$, first note that (all) projections $(\C_q^2)^{\otimes 2r} \twoheadrightarrow \C_q$ and inclusions 
$\C_q \hookrightarrow (\C_q^2)^{\otimes 2r}$ lie in the image of $\Gam$ since they can be built from the images of cap and cup 
morphisms\footnote{This is known as the \textit{first fundamental Theorem of $\mathfrak{sl}_2$-invariant theory} and already appears, in the non-quantum setting, in~\cite{rtw}.}. 
The result now follows since any (non-zero) $\Uq(\mathfrak{sl}_2)$-intertwiner $\Gam(\vec{k}) \to \Gam(\vec{l})$ can be written as the composition 
\[
\Gam(\vec{k}) \to (\C_q^2)^{\otimes 2r} \otimes \Gam(\vec{l}) \twoheadrightarrow \Gam(\vec{l}), \quad \text{if} \quad \sum_{i=1}^m k_i=2r + \sum_{i=1}^{m^{\prime}} l_i,
\]
or as the composition 
\[
\Gam(\vec{k}) \hookrightarrow (\C_q^2)^{\otimes 2r} \otimes \Gam(\vec{k}) \to \Gam(\vec{l}), \quad \text{if} \quad 2r + \sum_{i=1}^m k_i=\sum_{i=1}^{m^{\prime}} l_i.
\]


To see injectivity in~\eqref{eq-toverify}, we start by considering the case when $k_i=1$ and $l_j=1$ for all $1\leq i\leq m$ and $1\leq j\leq m^{\prime}$. 
We claim that in this case, any symmetric $\slnn{2}$-web 
$\Hom_{\SymSp}(\vec{k},\vec{l})$ can be expressed in terms of Temperley--Lieb diagrams. 
Indeed, given such a symmetric $\slnn{2}$-web, we can first use the digon (un)removals from~\eqref{eq-simpler1} to ``explode" the middle of each $k$-labeled edge into $1$-labeled edges. 
Using equation~\eqref{eq-frob}, we see that the symmetric $\slnn{2}$-web is now a multiple of one which is built entirely from symmetric Jones--Wenzl projectors, and we can use 
Proposition~\ref{prop-jw} and the symmetric Jones--Wenzl recursion from Lemma~\ref{lem-recursion} to write this symmetric $\slnn{2}$-web as an element of $\cal{TL}$.
Next, since the diagram 
\[
\begin{gathered}
\xymatrix{
\cal{TL} \ar[rr] \ar[dr]_{\cal{I}} & & \Repp \\
& \hspace*{0.75cm}\SymSp \ar[ur]_{\Gamma_{\mathbf{sym}}} &
}
\end{gathered}
\]
commutes and the top functor is fully faithful (see also Theorem~\ref{thm-equisl2}), 
it follows that $\cal{I}\colon\cal{TL} \rightarrow \SymSp$ is faithful. 
All together, this implies that
\[
\Hom_{\SymSp}(\vec{k},\vec{l}) \cong \Hom_{\cal{TL}}(\vec{k},\vec{l}) \cong \Hom_{\Repp}(\Gam(\vec{k}),\Gam(\vec{l})),
\]
for $\vec{k}=(1,\dots,1)\in\Z_{\geq 0}^k$ and $\vec{l}=(1,\dots,1)\in\Z_{\geq 0}^l$,
which in particular shows injectivity.

The general case follows from this.
Given any symmetric $\slnn{2}$-web $u\in\Hom_{\SymSp}(\vec{k},\vec{l})$ we can compose with merge and split morphisms to obtain
\[
\xy
(0,0)*{
\begin{tikzpicture} [scale=1]
\draw[very thick] (0,1.25) rectangle (2,.75);
\draw[very thick] (.25,0) to (.25,.75);
\draw[very thick] (.25,1.25) to (.25,2);
\draw[very thick] (1.75,0) to (1.75,.75);
\draw[very thick] (1.75,1.25) to (1.75,2);
\node at (1,1.75) {$\cdots$};
\node at (1,.25) {$\cdots$};
\node at (1,1) {$u$};
\node at (2.1,1.675) {\tiny $l_{m^{\prime}}$};
\node at (2.1,0.325) {\tiny $k_m$};
\node at (0,1.675) {\tiny $l_1$};
\node at (0,0.325) {\tiny $k_1$};
\draw [very thick] (0.5,-.5) to [out=90,in=330] (0.25,0);
\draw [very thick] (0,-.5) to [out=90,in=210] (0.25,0);
\draw [very thick] (2,-.5) to [out=90,in=330] (1.75,0);
\draw [very thick] (1.5,-.5) to [out=90,in=210] (1.75,0);
\draw [very thick] (0.25,2) to [out=30,in=270] (0.5,2.5);
\draw [very thick] (0.25,2) to [out=150,in=270] (0,2.5);
\draw [very thick] (1.75,2) to [out=30,in=270] (2,2.5);
\draw [very thick] (1.75,2) to [out=150,in=270] (1.5,2.5);
\node at (0.25,2.65) {\tiny $\cdots$};
\node at (0.25,-.65) {\tiny $\cdots$};
\node at (1.75,2.65) {\tiny $\cdots$};
\node at (1.75,-.65) {\tiny $\cdots$};
\node at (0,2.65) {\tiny $1$};
\node at (0.5,2.65) {\tiny $1$};
\node at (0,-.65) {\tiny $1$};
\node at (0.5,-.65) {\tiny $1$};
\node at (1.5,2.65) {\tiny $1$};
\node at (2,2.65) {\tiny $1$};
\node at (1.5,-.65) {\tiny $1$};
\node at (2,-.65) {\tiny $1$};
\end{tikzpicture}
};
\endxy\in\Hom_{\SymSp}(\underbrace{(1,\dots,1)}_{k_1+\dots+k_m},\underbrace{(1,\dots,1)}_{l_1+\dots+l_{m^{\prime}}})
\]
where we indicate with dots compositions of merge and split morphisms, the order of which do not matter due to the associativity relations~\eqref{eq-frob}.

The above argument, together with the digon removals from~\eqref{eq-simpler1}, 
shows that the images of two distinct symmetric $\slnn{2}$-webs $u,v\in\Hom_{\SymSp}(\vec{k},\vec{l})$ have to be distinct. 
Explicitly, the digon relations show that the splitting procedure is invertible while the argument above shows that the images of their ``enlargements'' are distinct.
\end{proof}

\begin{remark}\label{rem-pivotal}
We do not state and prove Theorem~\ref{thm-main} (following history) in terms of a \textit{pivotal} equivalence between $\SymSp$ and $\Repp$ 
due to an unavoidable sign issue coming from the use of unoriented diagrams. 
In our case, this arises since the vector representation $\C_q^2$ of $\Uq(\slnn{2})$ is \textit{anti}-symmetrically self-dual.
In order to incorporate this, we would have to make the diagrammatic calculus more 
sophisticated 
by introducing extra orientations and tag morphisms (as in~\cite{ckm}). 
Since these issues are usually not relevant to topological applications before categorifying or passing to the $\sln$ case, 
we avoid them for the time being and
stay closer to the ``traditional'' Temperley--Lieb calculus.
\end{remark}

\section{The colored Jones polynomial via symmetric \texorpdfstring{$\slnn{2}$}{sl2}-webs}\label{sec-colored}

\subsection{Braiding via quantum Weyl group elements}\label{sub-braiding}

In this subsection, we extend Theorem~\ref{thm-main} to incorporate the braided structure on $\Repp$. 
We begin by defining the following morphisms in $\Hom_{\SymSp}\big((k,l),(l,k)\big)$:
 \begin{equation}\label{eq-braidedstru}
 \beta_{k,l}^{\Sym} =
 \raisebox{-.05cm}{\xy
(0,0)*{
\begin{tikzpicture}[scale=.3]
	\draw [very thick, ->] (-1,-1) to (1,1);
	\draw [very thick, ->] (-0.25,0.25) to (-1,1);
	\draw [very thick] (0.25,-0.25) to (1,-1);
	\node at (-1,-1.5) {\tiny $k$};
	\node at (1,-1.5) {\tiny $l$};
\end{tikzpicture}
};
\endxy}
=
(-1)^{k}q^{-k-\frac{kl}{2}}\!\!\sum_{\begin{smallmatrix} j_1,j_2\geq 0 \\ j_1-j_2=k-l \end{smallmatrix}}\!\! (-q)^{j_1}
\xy
(0,0)*{
\begin{tikzpicture}[scale=.3]
	\draw [very thick] (-2,-4) to (-2,-2);
	\draw [very thick] (-2,-2) to (-2,0.25);
	\draw [very thick] (2,-4) to (2,-2);
	\draw [very thick] (2,-2) to (2,0.25);
	\draw [very thick, directed=.55] (-2,-2) to (2,-2);
	\draw [very thick] (-2,0.25) to (-2,2);
	\draw [very thick] (-2,2) to (-2,4);
	\draw [very thick] (2,0.25) to (2,2);
	\draw [very thick] (2,2) to (2,4);
	\draw [very thick, rdirected=.55] (-2,2) to (2,2);
	\node at (-2,-4.5) {\tiny $k$};
	\node at (2,-4.5) {\tiny $l$};
	\node at (-2,4.5) {\tiny $l$};
	\node at (2,4.5) {\tiny $k$};
	\node at (-3.25,0) {\tiny $k{-}j_1$};
	\node at (3.25,0) {\tiny $l{+}j_1$};
	\node at (0,-1.25) {\tiny $j_1$};
	\node at (0,2.75) {\tiny $j_2$};
\end{tikzpicture}
};
\endxy
 \end{equation}
which give rise to the braiding. More generally, for any two objects $\vec{k},\vec{l}$ in $\SymSp$ define 
\[
\beta_{\vec{k},\vec{l}}^{\Sym} =
\xy
(0,0)*{
\begin{tikzpicture}[scale=.3]
	\node at (0,0) {\tiny$k_1$};
	\node at (2,0) {$\ldots$};
	\node at (4,0) {\tiny$k_a$};
	\node at (6,0) {\tiny$l_1$};
	\node at (8,0) {$\ldots$};
	\node at (10,0) {\tiny$l_b$};
	\node at (0,7.5) {\tiny$l_1$};
	\node at (2,7.5) {$\ldots$};
	\node at (4,7.5) {\tiny$l_b$};
	\node at (6,7.5) {\tiny$k_1$};
	\node at (8,7.5) {$\ldots$};
	\node at (10,7.5) {\tiny$k_a$};
	\draw[very thick, ->] (0,.75) to (6,6.75);
	\draw[very thick, ->] (4,.75) to (10,6.75);
	\draw[very thick] (6,.75) to (5.25,1.5);
	\draw[very thick] ( 4.75,2) to (3.25,3.5);
	\draw[very thick, ->] (2.75,4) to (0,6.75);
	\draw[very thick] (10,.75) to (7.25,3.5) ;
	\draw[very thick] (6.75,4) to (5.35,5.5);
	\draw[very thick, ->] (4.75,6) to (4,6.75);
\end{tikzpicture}
};
\endxy
\in \Hom_{\SymSp}\big((k_1,\ldots,k_a,l_1,\ldots,l_b),(l_1, \ldots, l_b,k_1,\ldots,k_a)\big)
\]
by taking tensor products of compositions of the morphisms $\beta_{k,l}^{\Sym}$.
We now aim to show the following result. To understand it, recall that $\Repp$ is a 
\textit{braided} monoidal category where the braiding is induced via the 
$\slnn{2}$-$R$-matrix (the explicit construction of the 
braided monoidal structure on the category $\Repp$ can be found in many sources, e.g. 
Chapter XI, Section 2 and Section 7 in~\cite{tur}).

\begin{theorem}\label{thm-main2}
The morphisms $\beta_{\vec{k},\vec{l}}^{\Sym}$ define a braiding on $\SymSp$ and 
the additive closure of $\SymSp$ is braided monoidally equivalent to $\Repp$.
\end{theorem}

In particular, $\beta_{k,l}^{\Sym}$ is invertible, 
with inverse explicitly given by
\begin{equation}\label{eq-braidedinv}
 \left( \beta_{k,l}^{\Sym} \right)^{-1} =
 \raisebox{-.05cm}{\xy
(0,0)*{
\begin{tikzpicture}[scale=.3]
	\draw [very thick, ->] (1,-1) to (-1,1);
	\draw [very thick, ->] (0.25,0.25) to (1,1);
	\draw [very thick] (-0.25,-0.25) to (-1,-1);
	\node at (-1,-1.5) {\tiny $l$};
	\node at (1,-1.5) {\tiny $k$};
\end{tikzpicture}
};
\endxy}
=
(-1)^{k}q^{k+\frac{kl}{2}}\!\!\sum_{\begin{smallmatrix} j_1,j_2\geq 0 \\ j_1-j_2=k-l \end{smallmatrix}}\!\! (-q)^{-j_1}
\xy
(0,0)*{
\begin{tikzpicture}[scale=.3]
	\draw [very thick] (-2,-4) to (-2,-2);
	\draw [very thick] (-2,-2) to (-2,0.25);
	\draw [very thick] (2,-4) to (2,-2);
	\draw [very thick] (2,-2) to (2,0.25);
	\draw [very thick, directed=.55] (2,-2) to (-2,-2);
	\draw [very thick] (-2,0.25) to (-2,2);
	\draw [very thick] (-2,2) to (-2,4);
	\draw [very thick] (2,0.25) to (2,2);
	\draw [very thick] (2,2) to (2,4);
	\draw [very thick, rdirected=.55] (2,2) to (-2,2);
	\node at (-2,-4.5) {\tiny $l$};
	\node at (2,-4.5) {\tiny $k$};
	\node at (-2,4.5) {\tiny $k$};
	\node at (2,4.5) {\tiny $l$};
	\node at (-3.25,0) {\tiny $l{+}j_1$};
	\node at (3.25,0) {\tiny $k{-}j_1$};
	\node at (0,-1.25) {\tiny $j_1$};
	\node at (0,2.75) {\tiny $j_2$};
\end{tikzpicture}
};
\endxy
\end{equation}
as can be verified via a direct computation (compare also to Proposition 5.2.3 in~\cite{lu1}).

To prove Theorem~\ref{thm-main2}, we will again follow Cautis, Kamnitzer and Morrison (who in turn follow Lusztig~\cite{lu1} 
and Chuang and Rouquier~\cite{cr1}) by defining the operator\footnote{Formally, we must work over $\C(q^{\frac{1}{2}})$ to 
define these. Hence, we pass to these coefficients.}
\begin{equation}\label{eq-braiding1}
T_i\one_{\vec{k}} = (-1)^{k_i}q^{-k_i - \frac{k_i k_{i+1}}{2}}
\!\!\sum_{\begin{smallmatrix} j_1,j_2\geq 0 \\ j_1-j_2=k_i-k_{i+1} \end{smallmatrix}}\!\!(-q)^{j_1}E_i^{(j_2)}F_i^{(j_1)}\one_{\vec{k}}
\end{equation}
for any $\glm$-weight $\vec{k}\in\Z^m_{\geq 0}$ and any $i=1,\dots,m-1$, 
called \textit{Lusztig's $i$-th braiding operator}.

\begin{remark}\label{rem-finite}
These operators specify elements in $\Ud^{\infty}(\gli)$, since the 
sum in~\eqref{eq-braiding1} 
truncates to one which is finite. 
This is due to the fact that sufficiently high powers of 
$F_i^{(j_1)}\one_{\vec{k}}$ map to $\glm$-weights 
with negative entries, 
and hence, are zero in $\Ud^{\infty}(\gli)$.

We point out that the elements $T_i\one_{\vec{k}} \in \Ud^{\infty}(\gli)$ 
differ from the corresponding elements of Cautis, Kamnitzer and Morrison, 
both in that we work with (multiples of) Lusztig's $T_{i,+1}^{\prime\prime}$ (instead of $T_{i,-1}^{\prime\prime}$)
and since in their setting they kill all $\one_{\vec{k}}$ for $\glm$-weights $\vec{k}$ whose entries do not lie in $\{0,\dots,n\}$.
Fortunately, most of their calculations follow from those of Lusztig in Subsection 5.1.1 of~\cite{lu1}. 
Thus,
we can adopt most of Cautis, Kamnitzer, and Morrison's 
calculations without further issue.
\end{remark}

\begin{lemma}\label{lem-lustigcalc}
The $T_i\one_{\vec{k}}$ (viewed as elements of $\Ud^{\infty}(\gli)$) are invertible and satisfy the braid relations
\[
T_{i+1}T_{i}T_{i+1}\one_{\vec{k}}=T_{i}T_{i+1}T_{i}\one_{\vec{k}},\quad\text{and}\quad T_{i}T_{i^{\prime}}\one_{\vec{k}}=T_{i^{\prime}}T_{i}\one_{\vec{k}},\quad\text{if }|i-i^{\prime}|,
\]
for all  
$\glm$-weights $\vec{k}\in\Z^m_{\geq 0}$ and all $i,i^{\prime}=1,\dots,m-1$ (and all $m\in\Z_{\geq 0}$).
\end{lemma}

\begin{proof}
Almost word-for-word as in Lemma 6.1.1 and Lemma 6.1.2 from~\cite{ckm}.
\end{proof}

We now proceed with the proof of Theorem~\ref{thm-main2}.

\begin{proof}[Proof (of Theorem~\ref{thm-main2}).]
The one-line explanation is that both $\beta_{\vec{k},\vec{l}}^{\Sym}$ 
and the braiding on $\Repp$ come from Lusztig's 
braiding operator from~\eqref{eq-braiding1} above.

To be more thorough, we first introduce an analog of $\Ud^{\infty}(\gli)$ akin 
to the category studied by Cautis, Kamnitzer and Morrison. 
Let 
\[
\Ud^{\infty}(\glnn{\bullet}) = \bigoplus_{m >0} \Ud^{\infty}(\glm)
\]
which is in fact a monoidal category. For example, the tensor product is given
on objects by concatenating a $\glnn{m_1}$-weight with a $\glnn{m_2}$-weight
to obtain a $\glnn{m_1+m_2}$-weight (see Section 6 of~\cite{ckm} for more details).
Given a $\glnn{m_1}$-weight $\vec{k}$ and a $\glnn{m_2}$-weight $\vec{l}$, define 
the braiding operator
\[
\beta^{\infty}_{\vec{k}\vec{l}} = T_w \one_{\vec{k}\otimes\vec{l}},\quad\text{where $w$ is the permutation } w(i)=\begin{cases}
m_2+i, &\text{if }i\leq m_1,\\
i-m_1, &\text{if }i> m_1,\\
\end{cases}
\]
and $T_w = T_{i_1} \cdots T_{i_r}$ when $w = s_{i_1} \cdots s_{i_r}$ is a reduced 
expression (the choice of reduced expression does not matter by Lemma~\ref{lem-lustigcalc}).
A direct adaptation of Theorem 6.1.4 in~\cite{ckm} shows that these elements endow 
$\Ud^{\infty}(\glnn{\bullet})$ with the structure of a braided monoidal 
category (this uses again Lemma~\ref{lem-lustigcalc} which, as mentioned in 
Remark~\ref{rem-finite}, is based on calculations by Lusztig).

We now claim that the functors in the triangle
\[
\begin{gathered}
\xymatrix{
\Ud^{\infty}(\glnn{\bullet}) \ar[rr]^{\Phi_\bullet} \ar[dr]_{\Upsilon_\bullet} & & \Repp \\
& \SymSp \ar[ur]_{\Gam} &
}
\end{gathered}
\]
induced by the functors in the commuting diagram from~\eqref{eq-fac} are braided, 
which suffices to prove the result. 
The fact that $\SymSp$ is braided and that $\Upsilon_\bullet$ preserves the braiding follows 
directly by comparing equations~\eqref{eq-braidedstru} and~\eqref{eq-braiding1}
(and the fact that this functor is essentially surjective and full).

It finally suffices to show that $\Phi_{\bullet}$ is braided. Explicitly, we must check that 
\[
\Phi_{\bullet}(\beta^{\infty}_{\vec{k},\vec{l}})=
\beta^{R}_{\Phi_{\bullet}(\vec{k}),\Phi_{\bullet}(\vec{l})}=\beta^{R}_{\Sym_q^{k_1}\C_q^n\otimes\cdots\otimes\Sym_q^{k_m}\C_q^n,\Sym_q^{l_1}\C_q^n\otimes\cdots\otimes\Sym_q^{l_{m^\prime}}\C_q^n},
\]
where $\beta^R$ denotes the braiding coming from the $\slnn{2}$-$R$-matrix (as mentioned above).
To see this, we note that all of the steps used to prove Theorem 6.2.1 in~\cite{ckm} carry 
directly over to the symmetric case. 
Their arguments reduce to showing that $\Phi_{\bullet}(\beta_{1,1}^\infty) = \beta^R_{\C_q^2,\C_q^2}$, 
where the latter denotes the standard braiding on $\C_q^2\otimes\C_q^2$ given by the $\slnn{2}$-$R$-matrix.

To check this final equality, it suffices to show that when $k=1=l$,
equation~\eqref{eq-braidedstru} 
maps under $\Gam$ to the braiding on $\C_q^2 \otimes \C_q^2$. As mentioned in Lemma 6.2.2 of~\cite{ckm}, $\beta^R$ is determined on this by the fact that it acts by 
$q^{1/2}$ on $\Sym_q^2\C_q^2$ and by $-q^{-3/2}$ on $\bV_q^2\C_q^2$. 
In this case equation~\eqref{eq-braidedstru} is given by
\[
 \raisebox{-.05cm}{\xy
(0,0)*{
\begin{tikzpicture}[scale=.3]
	\draw [very thick, ->] (-1,-1) to (1,1);
	\draw [very thick, ->] (-0.25,0.25) to (-1,1);
	\draw [very thick] (0.25,-0.25) to (1,-1);
	\node at (-1,-1.5) {\tiny $1$};
	\node at (1,-1.5) {\tiny $1$};
\end{tikzpicture}
};
\endxy} = 
-q^{-3/2} \left(
\xy
(0,0)*{
\begin{tikzpicture}[scale=.2] 
	\draw [very thick] (1,-2.75) to (1,2.5);
	\draw [very thick] (-1,-2.75) to (-1,2.5);
	\node at (-1,3.25) {\tiny $1$};
	\node at (1,3.25) {\tiny $1$};
	\node at (-1,-3.5) {\tiny $1$};
	\node at (1,-3.5) {\tiny $1$};
\end{tikzpicture}
};
\endxy
- q
\xy
(0,0)*{
\begin{tikzpicture}[scale=.2]
	\draw [very thick] (0,-1) to (0,.75);
	\draw [very thick] (0,.75) to [out=30,in=270] (1,2.5);
	\draw [very thick] (0,.75) to [out=150,in=270] (-1,2.5); 
	\draw [very thick] (1,-2.75) to [out=90,in=330] (0,-1);
	\draw [very thick] (-1,-2.75) to [out=90,in=210] (0,-1);
	\node at (-1,3.25) {\tiny $1$};
	\node at (1,3.25) {\tiny $1$};
	\node at (-1,-3.5) {\tiny $1$};
	\node at (1,-3.5) {\tiny $1$};
	\node at (0.6,0) {\tiny $2$};
\end{tikzpicture}
};
\endxy
\right)
\]
and since the second term in 
$\Gam(\beta_{k,l}^{\Sym})$
factors through $\Sym_q^2\C_q^2$, this acts by $-q^{-3/2}$ on $\bV_q^2\C_q^2$. 
Similarly, equation~\eqref{eq-dumb} shows that the dumbbell
acts on $\Sym_q^2\C_q^2$ by multiplying with $[2]$. From this we see that 
$\Gam(\beta_{k,l}^{\Sym})$ acts by $-q^{-3/2} (1-q[2]) = q^{1/2}$ as desired.

Alternatively, we can check graphically that this 
agrees with the standard formula for a (positive) crossing in $\cal{TL}$. 
We compute that
\[
\raisebox{-.05cm}{\xy
(0,0)*{
\begin{tikzpicture}[scale=.3]
	\draw [very thick, ->] (-1,-1) to (1,1);
	\draw [very thick, ->] (-0.25,0.25) to (-1,1);
	\draw [very thick] (0.25,-0.25) to (1,-1);
	\node at (-1,-1.5) {\tiny $1$};
	\node at (1,-1.5) {\tiny $1$};
\end{tikzpicture}
};
\endxy} = 
-q^{-3/2} \left(
\xy
(0,0)*{
\begin{tikzpicture}[scale=.2] 
	\draw [very thick] (1,-2.75) to (1,2.5);
	\draw [very thick] (-1,-2.75) to (-1,2.5);
	\node at (-1,3.25) {\tiny $1$};
	\node at (1,3.25) {\tiny $1$};
	\node at (-1,-3.5) {\tiny $1$};
	\node at (1,-3.5) {\tiny $1$};
\end{tikzpicture}
};
\endxy
- q
\xy
(0,0)*{
\begin{tikzpicture}[scale=.2]
	\draw [very thick] (0,-1) to (0,.75);
	\draw [very thick] (0,.75) to [out=30,in=270] (1,2.5);
	\draw [very thick] (0,.75) to [out=150,in=270] (-1,2.5); 
	\draw [very thick] (1,-2.75) to [out=90,in=330] (0,-1);
	\draw [very thick] (-1,-2.75) to [out=90,in=210] (0,-1);
	\node at (-1,3.25) {\tiny $1$};
	\node at (1,3.25) {\tiny $1$};
	\node at (-1,-3.5) {\tiny $1$};
	\node at (1,-3.5) {\tiny $1$};
	\node at (0.6,0) {\tiny $2$};
\end{tikzpicture}
};
\endxy
\right)
=
q^{1/2} 
\xy
(0,0)*{
\begin{tikzpicture}[scale=.2] 
	\draw [very thick] (1,-2.75) to (1,2.5);
	\draw [very thick] (-1,-2.75) to (-1,2.5);
	\node at (-1,3.25) {\tiny $1$};
	\node at (1,3.25) {\tiny $1$};
	\node at (-1,-3.5) {\tiny $1$};
	\node at (1,-3.5) {\tiny $1$};
\end{tikzpicture}
};
\endxy
+
q^{-1/2}
\xy
(0,0)*{
\begin{tikzpicture}[scale=.2]
	\draw [very thick] (-1,2.5) to [out=270,in=180] (0,0.5) to [out=0,in=270] (1,2.5);
	\draw [very thick] (-1,-2.75) to [out=90,in=180] (0,-0.75) to [out=0,in=90] (1,-2.75);
	\node at (-1,3.25) {\tiny $1$};
	\node at (1,3.25) {\tiny $1$};
	\node at (-1,-3.5) {\tiny $1$};
	\node at (1,-3.5) {\tiny $1$};
\end{tikzpicture}
};
\endxy
\]
Here we remind the reader that the dumbbell can be replaced by $[2]$ times the identity plus a cap-cup.
This is the Kauffman bracket formula for the braiding on $\Repp$ (which is known to give the same result as the one coming from the $\slnn{2}$-$R$-matrix braiding).
\end{proof}

\begin{remark}\label{rem-howtosln}
More generally, the above argument extends without difficulties to show that the braiding in $\Repn$ between tensor products of 
the $\sln$-modules $\Sym_q^k\C_q^n$ coming from the $\sln$-$R$-matrix is given as the image of the braiding $\beta^\infty$ of $\Ud^{\infty}(\glnn{\bullet})$ 
under the functor $\Phi^{n}_{\bullet}\colon\Ud^{\infty}(\glnn{\bullet}) \to \Repn$ induced by the functors from~\eqref{eq-functorweights}.
\end{remark}

\begin{remark}\label{rem-skewsym}
In~\cite{ckm}, they show that the braiding between tensor 
products of fundamental representations in $\Rep{\sln}$
is similarly given by Lusztig's braiding operators 
$T_i \one_{\vec{k}} \in \Ud^{n}(\glnn{\bullet})$, where $\Ud^{n}(\glnn{\bullet})$ is the 
quotient of $\Ud^{\infty}(\glnn{\bullet})$ by all $\one_{\vec{k}}$'s with $\glnn{\bullet}$-weights 
containing an entry strictly lower than $0$ or strictly larger than $n$. 
In addition, they show that $q$-skew Howe duality gives a braided 
monoidal functor $\Ud^{n}(\glnn{\bullet}) \to \Rep{\sln}$. 
Since we have maps $\Ud^{\infty}(\glnn{\bullet}) \to \Ud^{n}(\glnn{\bullet})$, 
this gives the following diagram of 
braided monoidal functors
\[
\xymatrix{
& \Ud^{\infty}(\glnn{\bullet}) \ar[ld]|{q\text{-}\text{skew Howe}} \ar[rd]|{q\text{-}\text{sym. Howe}} & \\
\Rep{\sln} & & \Repp
}
\]
We again point out that there is a slight difference between the 
$q$-symmetric Howe duality and the $q$-skew Howe duality cases coming from the fact that we 
need to use Lusztig's $T_{i,+1}^{\prime\prime}$ instead of $T_{i,-1}^{\prime\prime}$, which is utilized by Cautis, Kamnitzer, and Morrison. 
Since $T_{i,+1}^{\prime\prime}$ and $T_{i,-1}^{\prime\prime}$ only differ by a substitution of $q \leftrightarrow q^{-1}$, this gives the schematic
\[
\xymatrix{
& \beta^{\infty} \ar@{|->}[ld]_{q^{-1}} \ar@{|->}[rd]^{q} & \\
\beta^R_n & & \beta^R
}
\]
where we point out that $\beta^R_n$ is the braiding of $\Rep{\sln}$ coming from the 
$\sln$-$R$-matrix while $\beta^R$ is the braiding of $\Repp$ coming from the 
$\slnn{2}$-$R$-matrix

This observation appears related to the decategorification of the ``mirror symmetry'' between colored 
HOMFLY--PT homology
conjectured in~\cite{gs} (e.g. in (5-17) in their paper). 
See Section~\ref{sub-mirror} below for a more precise discussion.
\end{remark}

\subsection{The colored Jones polynomial via ``MOY''-graphs}\label{sub-myfavourite} 

In this subsection we explore how the braiding from Subsection~\ref{sub-braiding} on the symmetric $\slnn{2}$-spider 
can be used to study the \textit{colored Jones polynomial} of colored, oriented links $L$, which we denote by $\cal{J}_{\vec{c}}(L_D)$.
Here $\vec{c}=(c_1,\dots,c_N)$ denotes the colors of the $N$-component, oriented link $L$
and
$L_D$ is a colored, oriented diagram for $L$. 
In the interest of brevity, we refer the reader to the wide literature on the subject, 
in particular Chapter XI, Section 7 of~\cite{tur}, for the definition of this invariant and a thorough treatment of its properties 
(see~\cite{jonespoly} for Jones' original work on link polynomials).
We only comment that it can be computed by associating a 
morphism between trivial representations in $\Repp$ to any colored, oriented link diagram 
$L_D$ of a colored, oriented link $L$ 
(and rescaling to get an invariant which is not framing-dependent).

This translates to using equations~\eqref{eq-braidedstru} and~\eqref{eq-braidedinv} to view a colored, 
oriented link diagram for $L_D$ as a morphism in $\SymSp$, 
which necessarily evaluates to an element in $\C(q^{\frac{1}{2}})$ (in fact, 
it is clear from our construction that this always gives an element in $\Z[q^{\frac{1}{2}},q^{-\frac{1}{2}}]$), and multiplying by a 
certain normalization factor which can be computed directly from the crossing data of the diagram. 
For the case of a $c$-colored knot $K$ with diagram $K_D$, this factor is 
$(-1)^c q^{- C \omega(K_D)}$, 
where $\omega(K_D)$ is the writhe\footnote{The writhe is the difference between the number of positive $\xy
(0,0)*{
\begin{tikzpicture}[scale=.3]
	\draw [very thick, ->] (-1,-1) to (1,1);
	\draw [very thick, ->] (-0.25,0.25) to (-1,1);
	\draw [very thick] (0.25,-0.25) to (1,-1);
\end{tikzpicture}
};
\endxy$ and negative crossings $\raisebox{.085cm}{\reflectbox{\xy
(0,0)*{
\begin{tikzpicture}[scale=.3]
	\draw [very thick, ->] (-1,-1) to (1,1);
	\draw [very thick, ->] (-0.25,0.25) to (-1,1);
	\draw [very thick] (0.25,-0.25) to (1,-1);
\end{tikzpicture}
};
\endxy}}$ in the diagram.} 
of $K_D$
and 
$C=\frac{c^2+2c}{2}$ is the so-called \textit{quadratic Casimir number} of the color $\Sym_q^c\C_q^{2}$. 
In general, for a colored, oriented link diagram $L_D$ one normalizes by multiplying by the product of the 
normalization factors for each of the colored link components.

We note that this approach is similar in the $1$-colored case to computing the Jones polynomial using the Kauffman bracket, 
but in the colored case avoids the use of cabling and Jones--Wenzl projectors, 
trading them instead for our ``symmetric version''
of the MOY-calculus~\cite{moy} typically used to compute 
the $\bV_q^k\C^n_q$-colored $\sln$-link invariant.

\begin{example}\label{ex-maincalcsmall}
As an example, we compute the ($1$-colored) Jones polynomial 
of the Hopf link using symmetric $\slnn{2}$-webs.
\begin{align*}
\cal{J}_{(1,1)} \left(
\xy
(0,0)*{
\begin{tikzpicture} [scale=.6]
\draw[very thick] (0,0) to [out=90,in=270] (.5,1);
\draw[very thick] (.2,.6) to [out=135,in=270] (0,1);
\draw[very thick] (.5,0) to [out=90,in=315] (.3,.4);
\draw[very thick] (0,-1) to [out=90,in=270] (.5,0);
\draw[very thick] (.2,-.4) to [out=135,in=270] (0,0);
\draw[very thick] (.5,-1) to [out=90,in=315] (.3,-.6);
\draw[very thick, directed=.5] (.5,1) to [out=90,in=180] (.75,1.4) to [out=0,in=90] (1,1) to 
	(1,-1) to [out=270,in=0] (.75,-1.4) to [out=180,in=270] (.5,-1);
\draw[very thick, directed=.5] (0,1) to [out=90,in=180] (.75,1.8) to [out=0,in=90] (1.5,1) to 
	(1.5,-1) to [out=270,in=0] (.75,-1.8) to [out=180,in=270] (0,-1);
	\node at (0.85,-1) {\tiny $1$};
	\node at (1.35,-1) {\tiny $1$};
\end{tikzpicture}
}
\endxy
\right)
&= q^{-3} \;
\xy
(0,0)*{
\begin{tikzpicture} [scale=.6]
\draw[very thick] (0,0) to 
					(0,1);
\draw[very thick] (.5,0) to 
					(.5,1);
\draw[very thick] (0,-1) to 
					(0,0);
\draw[very thick] (.5,-1) to 
					(.5,0);
\draw[very thick] (.5,1) to [out=90,in=180] (.75,1.4) to [out=0,in=90] (1,1) to 
	(1,-1) to [out=270,in=0] (.75,-1.4) to [out=180,in=270] (.5,-1);
\draw[very thick] (0,1) to [out=90,in=180] (.75,1.8) to [out=0,in=90] (1.5,1) to 
	(1.5,-1) to [out=270,in=0] (.75,-1.8) to [out=180,in=270] (0,-1);
	\node at (0.85,-1) {\tiny $1$};
	\node at (1.35,-1) {\tiny $1$};
\end{tikzpicture}
}
\endxy
\; - q^{-2} \; 
\xy
(0,0)*{
\begin{tikzpicture} [scale=.6]
\draw[very thick] (0,0) to [out=90,in=225] (.25,.4) to [out=315,in=90] (.5,0);
\draw[very thick] (0,1) to [out=270,in=135] (.25,.6) to [out=45,in=270] (.5,1);
\draw[very thick] (.25,.4) to (.25,.6);
\draw[very thick] (0,-1) to 
					(0,0);
\draw[very thick] (.5,-1) to 
					(.5,0);
\draw[very thick] (.5,1) to [out=90,in=180] (.75,1.4) to [out=0,in=90] (1,1) to 
	(1,-1) to [out=270,in=0] (.75,-1.4) to [out=180,in=270] (.5,-1);
\draw[very thick] (0,1) to [out=90,in=180] (.75,1.8) to [out=0,in=90] (1.5,1) to 
	(1.5,-1) to [out=270,in=0] (.75,-1.8) to [out=180,in=270] (0,-1);
	\node at (0.85,-1) {\tiny $1$};
	\node at (1.35,-1) {\tiny $1$};
	\node at (0.4,0.5) {\tiny $2$};
\end{tikzpicture}
}
\endxy
\; - q^{-2} \;
\xy
(0,0)*{
\begin{tikzpicture} [scale=.6]
\draw[very thick] (0,0) to 
					(0,1);
\draw[very thick] (.5,0) to 
					(.5,1);
\draw[very thick] (0,-1) to [out=90,in=225] (.25,-.6) to [out=315,in=90] (.5,-1);
\draw[very thick] (0,0) to [out=270,in=135] (.25,-.4) to [out=45,in=270] (.5,0);
\draw[very thick] (.25,-.4) to (.25,-.6);
\draw[very thick] (.5,1) to [out=90,in=180] (.75,1.4) to [out=0,in=90] (1,1) to 
	(1,-1) to [out=270,in=0] (.75,-1.4) to [out=180,in=270] (.5,-1);
\draw[very thick] (0,1) to [out=90,in=180] (.75,1.8) to [out=0,in=90] (1.5,1) to 
	(1.5,-1) to [out=270,in=0] (.75,-1.8) to [out=180,in=270] (0,-1);
	\node at (0.85,-1) {\tiny $1$};
	\node at (1.35,-1) {\tiny $1$};
	\node at (0.4,-0.5) {\tiny $2$};
\end{tikzpicture}
}
\endxy
\; + \;
q^{-1} \;
\xy
(0,0)*{
\begin{tikzpicture} [scale=.6]
\draw[very thick] (0,0) to [out=90,in=225] (.25,.4) to [out=315,in=90] (.5,0);
\draw[very thick] (0,1) to [out=270,in=135] (.25,.6) to [out=45,in=270] (.5,1);
\draw[very thick] (.25,.4) to (.25,.6);
\draw[very thick] (0,-1) to [out=90,in=225] (.25,-.6) to [out=315,in=90] (.5,-1);
\draw[very thick] (0,0) to [out=270,in=135] (.25,-.4) to [out=45,in=270] (.5,0);
\draw[very thick] (.25,-.4) to (.25,-.6);
\draw[very thick] (.5,1) to [out=90,in=180] (.75,1.4) to [out=0,in=90] (1,1) to 
	(1,-1) to [out=270,in=0] (.75,-1.4) to [out=180,in=270] (.5,-1);
\draw[very thick] (0,1) to [out=90,in=180] (.75,1.8) to [out=0,in=90] (1.5,1) to 
	(1.5,-1) to [out=270,in=0] (.75,-1.8) to [out=180,in=270] (0,-1);
	\node at (0.85,-1) {\tiny $1$};
	\node at (1.35,-1) {\tiny $1$};
	\node at (0.4,-0.5) {\tiny $2$};
	\node at (0.4,0.5) {\tiny $2$};
	\node at (0.15,0) {\tiny $1$};
	\node at (0.65,0) {\tiny $1$};
\end{tikzpicture}
}
\endxy \\
&= q^{-3} [2]^2 - 2 q^{-2} [2] [3] + q^{-1} [2]^2 [3] = [2] (q^{2}+q^{-2}) = [4].
\end{align*}
\end{example}

\begin{example}\label{ex-maincalcsmall2}
As another example, we compute the $(2,1)$-colored Jones polynomial 
of the Hopf link using our approach.
\begin{align*}
\cal{J}_{(2,1)} \left(
\xy
(0,0)*{
\begin{tikzpicture} [scale=.6]
\draw[very thick] (0,0) to [out=90,in=270] (.5,1);
\draw[very thick] (.2,.6) to [out=135,in=270] (0,1);
\draw[very thick] (.5,0) to [out=90,in=315] (.3,.4);
\draw[very thick] (0,-1) to [out=90,in=270] (.5,0);
\draw[very thick] (.2,-.4) to [out=135,in=270] (0,0);
\draw[very thick] (.5,-1) to [out=90,in=315] (.3,-.6);
\draw[very thick, directed=.5] (.5,1) to [out=90,in=180] (.75,1.4) to [out=0,in=90] (1,1) to 
	(1,-1) to [out=270,in=0] (.75,-1.4) to [out=180,in=270] (.5,-1);
\draw[very thick, directed=.5] (0,1) to [out=90,in=180] (.75,1.8) to [out=0,in=90] (1.5,1) to 
	(1.5,-1) to [out=270,in=0] (.75,-1.8) to [out=180,in=270] (0,-1);
	\node at (0.85,-1) {\tiny $1$};
	\node at (1.35,-1) {\tiny $2$};
\end{tikzpicture}
}
\endxy
\right)
&= q^{-4} \;
\xy
(0,0)*{
\begin{tikzpicture} [scale=.6]
\draw[very thick] (0,0) to (0,1);
\draw[very thick] (.5,0) to (.5,1);
\draw[very thick, rdirected=.6] (0,.5) to (.5,.5);
\draw[very thick] (0,-1) to (0,0);
\draw[very thick] (.5,-1) to (.5,0);
\draw[very thick, directed=.7] (0,-.5) to (.5,-.5);
\draw[very thick, directed=.5] (.5,1) to [out=90,in=180] (.75,1.4) to [out=0,in=90] (1,1) to 
	(1,-1) to [out=270,in=0] (.75,-1.4) to [out=180,in=270] (.5,-1);
\draw[very thick, directed=.5] (0,1) to [out=90,in=180] (.75,1.8) to [out=0,in=90] (1.5,1) to 
	(1.5,-1) to [out=270,in=0] (.75,-1.8) to [out=180,in=270] (0,-1);
	\node at (0.85,-1) {\tiny $1$};
	\node at (1.35,-1) {\tiny $2$};
	\node at (.25,.8) {\tiny $1$};
	\node at (.25,-.8) {\tiny $1$};
	\node at (0.15,0) {\tiny $1$};
	\node at (0.65,0) {\tiny $2$};
\end{tikzpicture}
}
\endxy
\; - \;
q^{-3} \;
\xy
(0,0)*{
\begin{tikzpicture} [scale=.6]
\draw[very thick] (0,0) to [out=90,in=225] (.25,.4) to [out=315,in=90] (.5,0);
\draw[very thick] (0,1) to [out=270,in=135] (.25,.6) to [out=45,in=270] (.5,1);
\draw[very thick] (.25,.4) to (.25,.6);
\draw[very thick] (0,-1) to (0,0);
\draw[very thick] (.5,-1) to (.5,0);
\draw[very thick, directed=.7] (0,-.5) to (.5,-.5);
\draw[very thick] (.5,1) to [out=90,in=180] (.75,1.4) to [out=0,in=90] (1,1) to 
	(1,-1) to [out=270,in=0] (.75,-1.4) to [out=180,in=270] (.5,-1);
\draw[very thick] (0,1) to [out=90,in=180] (.75,1.8) to [out=0,in=90] (1.5,1) to 
	(1.5,-1) to [out=270,in=0] (.75,-1.8) to [out=180,in=270] (0,-1);
	\node at (0.85,-1) {\tiny $1$};
	\node at (1.35,-1) {\tiny $2$};
	\node at (.25,-.8) {\tiny $1$};
	\node at (0.4,0.5) {\tiny $3$};
	\node at (0.15,0) {\tiny $1$};
	\node at (0.65,0) {\tiny $2$};
\end{tikzpicture}
}
\endxy
\; - \;
q^{-3} \;
\xy
(0,0)*{
\begin{tikzpicture} [scale=.6]
\draw[very thick] (0,0) to (0,1);
\draw[very thick] (.5,0) to (.5,1);
\draw[very thick, rdirected=.6] (0,.5) to (.5,.5);
\draw[very thick] (0,-1) to [out=90,in=225] (.25,-.6) to [out=315,in=90] (.5,-1);
\draw[very thick] (0,0) to [out=270,in=135] (.25,-.4) to [out=45,in=270] (.5,0);
\draw[very thick] (.25,-.4) to (.25,-.6);
\draw[very thick] (.5,1) to [out=90,in=180] (.75,1.4) to [out=0,in=90] (1,1) to 
	(1,-1) to [out=270,in=0] (.75,-1.4) to [out=180,in=270] (.5,-1);
\draw[very thick] (0,1) to [out=90,in=180] (.75,1.8) to [out=0,in=90] (1.5,1) to 
	(1.5,-1) to [out=270,in=0] (.75,-1.8) to [out=180,in=270] (0,-1);
	\node at (0.85,-1) {\tiny $1$};
	\node at (1.35,-1) {\tiny $2$};
	\node at (0.4,-0.5) {\tiny $3$};
	\node at (.25,.8) {\tiny $1$};
	\node at (0.15,0) {\tiny $1$};
	\node at (0.65,0) {\tiny $2$};
\end{tikzpicture}
}
\endxy
\; + \;
q^{-2} \;
\xy
(0,0)*{
\begin{tikzpicture} [scale=.6]
\draw[very thick] (0,0) to [out=90,in=225] (.25,.4) to [out=315,in=90] (.5,0);
\draw[very thick] (0,1) to [out=270,in=135] (.25,.6) to [out=45,in=270] (.5,1);
\draw[very thick] (.25,.4) to (.25,.6);
\draw[very thick] (0,-1) to [out=90,in=225] (.25,-.6) to [out=315,in=90] (.5,-1);
\draw[very thick] (0,0) to [out=270,in=135] (.25,-.4) to [out=45,in=270] (.5,0);
\draw[very thick] (.25,-.4) to (.25,-.6);
\draw[very thick] (.5,1) to [out=90,in=180] (.75,1.4) to [out=0,in=90] (1,1) to 
	(1,-1) to [out=270,in=0] (.75,-1.4) to [out=180,in=270] (.5,-1);
\draw[very thick] (0,1) to [out=90,in=180] (.75,1.8) to [out=0,in=90] (1.5,1) to 
	(1.5,-1) to [out=270,in=0] (.75,-1.8) to [out=180,in=270] (0,-1);
	\node at (0.85,-1) {\tiny $1$};
	\node at (1.35,-1) {\tiny $2$};
	\node at (0.4,-0.5) {\tiny $3$};
	\node at (0.4,0.5) {\tiny $3$};
	\node at (0.15,0) {\tiny $1$};
	\node at (0.65,0) {\tiny $2$};
\end{tikzpicture}
}
\endxy \\
&= -q^{-4} [2][3]^2 + 2q^{-3}[2][3][4]-q^{-2}[3]^2[4]=-[6].
\end{align*} 
\end{example}
In both of these examples, we recover the known formula for the 
colored Jones polynomial of the Hopf link
\[
\cal{J}_{(k,l)} \left(
\xy
(0,0)*{
\begin{tikzpicture} [scale=.6]
\draw[very thick] (0,0) to [out=90,in=270] (.5,1);
\draw[very thick] (.2,.6) to [out=135,in=270] (0,1);
\draw[very thick] (.5,0) to [out=90,in=315] (.3,.4);
\draw[very thick] (0,-1) to [out=90,in=270] (.5,0);
\draw[very thick] (.2,-.4) to [out=135,in=270] (0,0);
\draw[very thick] (.5,-1) to [out=90,in=315] (.3,-.6);
\draw[very thick, directed=.5] (.5,1) to [out=90,in=180] (.75,1.4) to [out=0,in=90] (1,1) to 
	(1,-1) to [out=270,in=0] (.75,-1.4) to [out=180,in=270] (.5,-1);
\draw[very thick, directed=.5] (0,1) to [out=90,in=180] (.75,1.8) to [out=0,in=90] (1.5,1) to 
	(1.5,-1) to [out=270,in=0] (.75,-1.8) to [out=180,in=270] (0,-1);
	\node at (0.85,-1) {\tiny $l$};
	\node at (1.35,-1) {\tiny $k$};
\end{tikzpicture}
}
\endxy
\right)= (-1)^{k+l} [(k+1)(l+1)]
\] 
for $k,l\in\Z_{\geq 0}$ (the $(-1)^{k+l}$ factor comes from our conventions, 
and varies in the literature).

\subsection{A remark on ``mirror symmetry''}\label{sub-mirror}

We now aim to give a slightly more precise formulation of the ``mirror symmetry'' phenomena 
mentioned in Remark~\ref{rem-skewsym}. 
Consider the \textit{generic spider} $\GenSp$, 
the category whose objects are tuples in the symbols $k^{\pm}$ for $k \in \Z_{> 0}$, 
and whose morphisms are $\C(q)$-linear combinations of \textit{generic webs}, that is, 
oriented, trivalent, planar generated by
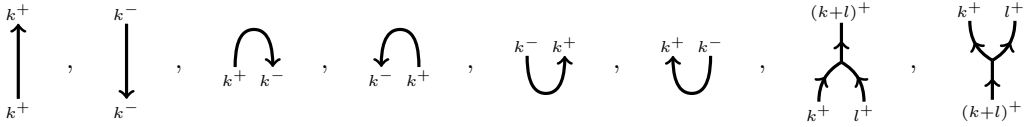
\[
\xy
(0,0)*{
\begin{tikzpicture} [scale=1]
\draw[very thick, ->] (0,0) to (0,1);
\node at (0,-.15) {\tiny $k^+$};
\node at (0,1.15) {\tiny $k^+$};
\end{tikzpicture}
}
\endxy
\quad , \quad
\xy
(0,0)*{
\begin{tikzpicture} [scale=1]
\draw[very thick, <-] (0,0) to (0,1);
\node at (0,-.15) {\tiny $k^-$};
\node at (0,1.15) {\tiny $k^-$};
\end{tikzpicture}
}
\endxy
\quad , \quad
\xy
(0,0)*{
\begin{tikzpicture} [scale=1]
\draw[very thick,->] (0,0) to [out=90,in=180] (.25,.5) to [out=0,in=90] (.5,0);
\node at (0,-.15) {\tiny $k^+$};
\node at (.5,-.15) {\tiny $k^-$};
\end{tikzpicture}
}
\endxy
\quad , \quad
\xy
(0,0)*{
\begin{tikzpicture} [scale=1]
\draw[very thick,<-] (0,0) to [out=90,in=180] (.25,.5) to [out=0,in=90] (.5,0);
\node at (0,-.15) {\tiny $k^-$};
\node at (.5,-.15) {\tiny $k^+$};
\end{tikzpicture}
}
\endxy
\quad , \quad
\xy
(0,0)*{
\begin{tikzpicture} [scale=1]
\draw[very thick,->] (0,1) to [out=270,in=180] (.25,.5) to [out=0,in=270] (.5,1);
\node at (0,1.15) {\tiny $k^-$};
\node at (.5,1.15) {\tiny $k^+$};
\end{tikzpicture}
}
\endxy
\quad , \quad
\xy
(0,0)*{
\begin{tikzpicture} [scale=1]
\draw[very thick,<-] (0,1) to [out=270,in=180] (.25,.5) to [out=0,in=270] (.5,1);
\node at (0,1.15) {\tiny $k^+$};
\node at (.5,1.15) {\tiny $k^-$};
\end{tikzpicture}
}
\endxy
\quad , \quad
\xy
(0,0)*{
\begin{tikzpicture}[scale=.3]
	\draw [very thick,directed=.55] (0, .75) to (0,2.5);
	\draw [very thick,directed=.55] (1,-1) to [out=90,in=330] (0,.75);
	\draw [very thick,directed=.55] (-1,-1) to [out=90,in=210] (0,.75); 
	\node at (0, 3) {\tiny $(k{+}l)^+$};
	\node at (-1,-1.5) {\tiny $k^+$};
	\node at (1,-1.5) {\tiny $l^+$};
\end{tikzpicture}
};
\endxy
\quad , \quad
\xy
(0,0)*{
\begin{tikzpicture}[scale=.3]
	\draw [very thick,directed=.55] (0,-1) to (0,.75);
	\draw [very thick,directed=.55] (0,.75) to [out=30,in=270] (1,2.5);
	\draw [very thick,directed=.55] (0,.75) to [out=150,in=270] (-1,2.5); 
	\node at (0, -1.5) {\tiny $(k{+}l)^+$};
	\node at (-1,3) {\tiny $k^+$};
	\node at (1,3) {\tiny $l^+$};
\end{tikzpicture}
};
\endxy
\]
modulo planar isotopy and the (oriented) standard $\sln$-web relations from equations~\eqref{eq-frob},~\eqref{eq-simpler1} and~\eqref{eq-almostsimpler}.
The oriented version of Lemma~\ref{lem-ladderfunc} gives a functor $\Upsilon_{\mathrm{Gen}}\colon\Ud^{\infty}(\gli)\to\GenSp$.
Since $\GenSp$ clearly admits specialization functors 
to $\spid{n}$ and $\SymSp$, which we denote by $\cal{S}_{\wedge}$ and $\cal{S}_{\Sym}$ respectively, we have the following commuting diagram:
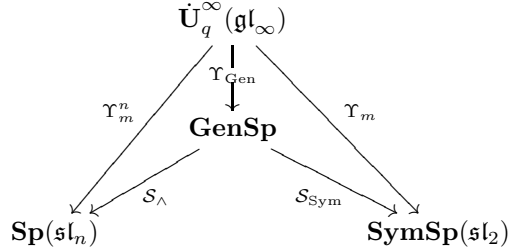
\[
\xymatrix{
& \Ud^{\infty}(\gli) \ar[d]|{\Upsilon_{\mathrm{Gen}}} \ar[ddl]_{\Upsilon^n_{m}} \ar[ddr]^{\Upsilon_{m}}& \\
& \GenSp \ar[dl]^{\cal{S}_{\wedge}} \ar[dr]_{\cal{S}_{\Sym}} \\
\spid{n} & & \SymSp
}
\]
Here $\Upsilon^n_{m}$ is the functor from Subsection 5.2. in~\cite{ckm}.

Given a colored, oriented braid $B$, the non-rescaled crossing formulae
\begin{equation}\label{eq-genbraiding}
 \phantom{aaaa.}\beta_{k,l}^{\mathrm{Gen}} =
 \raisebox{-.05cm}{\xy
(0,0)*{
\begin{tikzpicture}[scale=.3]
	\draw [very thick, ->] (-1,-1) to (1,1);
	\draw [very thick, ->] (-0.25,0.25) to (-1,1);
	\draw [very thick] (0.25,-0.25) to (1,-1);
	\node at (-1,-1.5) {\tiny $k$};
	\node at (1,-1.5) {\tiny $l$};
\end{tikzpicture}
};
\endxy}
=\!\!
\sum_{\begin{smallmatrix} j_1,j_2\geq 0 \\ j_1-j_2=k-l \end{smallmatrix}}\!\! (-q)^{j_1}
\xy
(0,0)*{
\begin{tikzpicture}[scale=.3]
	\draw [very thick] (-2,-4) to (-2,-2);
	\draw [very thick] (-2,-2) to (-2,0.25);
	\draw [very thick] (2,-4) to (2,-2);
	\draw [very thick] (2,-2) to (2,0.25);
	\draw [very thick, directed=.55] (-2,-2) to (2,-2);
	\draw [very thick] (-2,0.25) to (-2,2);
	\draw [very thick] (-2,2) to (-2,4);
	\draw [very thick] (2,0.25) to (2,2);
	\draw [very thick] (2,2) to (2,4);
	\draw [very thick, rdirected=.55] (-2,2) to (2,2);
	\node at (-2,-4.5) {\tiny $k$};
	\node at (2,-4.5) {\tiny $l$};
	\node at (-2,4.5) {\tiny $l$};
	\node at (2,4.5) {\tiny $k$};
	\node at (-3.25,0) {\tiny $k{-}j_1$};
	\node at (3.25,0) {\tiny $l{+}j_1$};
	\node at (0,-1.25) {\tiny $j_1$};
	\node at (0,2.75) {\tiny $j_2$};
\end{tikzpicture}
};
\endxy
\end{equation}
and its inverse
\begin{equation}\label{eq-genbraiding2}
 \left(\beta_{k,l}^{\mathrm{Gen}}\right)^{-1} =
 \raisebox{-.05cm}{\xy
(0,0)*{
\begin{tikzpicture}[scale=.3]
	\draw [very thick, ->] (1,-1) to (-1,1);
	\draw [very thick, ->] (0.25,0.25) to (1,1);
	\draw [very thick] (-0.25,-0.25) to (-1,-1);
	\node at (-1,-1.5) {\tiny $l$};
	\node at (1,-1.5) {\tiny $k$};
\end{tikzpicture}
};
\endxy}
=\!\!
\sum_{\begin{smallmatrix} j_1,j_2\geq 0 \\ j_1-j_2=k-l \end{smallmatrix}}\!\! (-q)^{-j_1}
\xy
(0,0)*{
\begin{tikzpicture}[scale=.3]
	\draw [very thick] (-2,-4) to (-2,-2);
	\draw [very thick] (-2,-2) to (-2,0.25);
	\draw [very thick] (2,-4) to (2,-2);
	\draw [very thick] (2,-2) to (2,0.25);
	\draw [very thick, directed=.55] (2,-2) to (-2,-2);
	\draw [very thick] (-2,0.25) to (-2,2);
	\draw [very thick] (-2,2) to (-2,4);
	\draw [very thick] (2,0.25) to (2,2);
	\draw [very thick] (2,2) to (2,4);
	\draw [very thick, rdirected=.55] (2,2) to (-2,2);
	\node at (-2,-4.5) {\tiny $l$};
	\node at (2,-4.5) {\tiny $k$};
	\node at (-2,4.5) {\tiny $k$};
	\node at (2,4.5) {\tiny $l$};
	\node at (-3.25,0) {\tiny $l{+}j_1$};
	\node at (3.25,0) {\tiny $k{-}j_1$};
	\node at (0,-1.25) {\tiny $j_1$};
	\node at (0,2.75) {\tiny $j_2$};
\end{tikzpicture}
};
\endxy
\end{equation}
assign a morphism in $\GenSp$ to its closure $\mathrm{cl}(B)$, 
which maps to a multiple of the 
colored Jones polynomial of $\mathrm{cl}(B)$ under the functor 
$\cal{S}_{\Sym}\colon\GenSp \to \SymSp$.

Note that  this element also maps to a multiple of the 
$\bV_q^k\C^n_q$-colored $\sln$-link polynomial 
by first making the substitution $q \leftrightarrow q^{-1}$ and 
then applying the functor
$\cal{S}_{\wedge}\colon\GenSp \to \spid{n}$, since the braiding 
in $\spid{n}$ is given by a multiple of the image of 
equations~\eqref{eq-genbraiding} and~\eqref{eq-genbraiding2} after making this substitution.

The relations in $\GenSp$ suffice\footnote{This fact was observed in joint work between the first author 
and Queffelec at the categorical level~\cite{qr2}. 
See also Queffelec's recent preprint with Sartori~\cite{qs} which utilizes and outlines the decategorified statement.} 
to express any closure of a generic web appearing in the morphisms assigned to a 
(colored, oriented) braid in terms of 
colored circles.
Viewing these colored circles as parameters $\{\xi_i\}_{i=1}^{\infty}$, we arrive at the following result.

\begin{theorem}\label{thm-WeakMirrorSymmetry}
There exists an invariant of (colored, oriented) braid conjugacy classes
\[
P(B) \in \Z[q,q^{-1},\xi_1,\xi_2,\ldots] \big/ \mathrm{I},
\] 
where $\mathrm{I}$ is the (possibly empty) ideal of
relations between colored circles in $\GenSp$. 
The specialization $\xi_k = [k+1]$ gives a multiple of the colored Jones polynomial of the closure $\mathrm{cl}(B)$ of $B$, 
and the substitution $q \leftrightarrow q^{-1}$ and subsequent 
specialization $\xi_k = \qbin{n}{k}$ for $k=0,\dots,n$ and $\xi_{>n}=0$ gives 
a multiple of the colored $\sln$-link polynomial of the closure $\mathrm{cl}(B)$ of $B$.
\end{theorem}

\begin{remark}
The method for computing the colored Jones polynomial described before Theorem~\ref{thm-WeakMirrorSymmetry} 
is distinct from the standard method of computation. 
Traditionally, one uses the Jones--Wenzl recursions to express a cabled link with Jones--Wenzl projectors inserted along cabled components 
as a linear combination of $1$-labeled circles, and then evaluates using equation~\eqref{eq-circle}. 
Our method rather expressed a colored braid closure as a linear combination of colored circles, before evaluating these circles as in Example~\ref{ex-circle}.
\end{remark}

Theorem \ref{thm-WeakMirrorSymmetry} is related to properties of the 
HOMFLY--PT polynomial, which were originally defined in~\cite{homfly} and~\cite{pt}.
The following is a slight generalization of the symmetric-skew ``mirror symmetry'' conjecture 
of Gukov and Sto\v{s}i\'{c}, see e.g. (5-17) in~\cite{gs}.

\begin{conjecture}\label{conj-never-conjecture-anything}
There exists a specialization of $P(B)$ which gives a multiple of the $\Sym^k_q$-colored 
HOMFLY--PT polynomial. Applying the substitution $q \leftrightarrow q^{-1}$ yields a multiple of the 
$\bV^k_q$-colored HOMFLY--PT polynomial.
\end{conjecture}

A proof of Conjecture~\ref{conj-never-conjecture-anything} using 
our methods would yield a diagrammatic proof of the 
``mirror symmetry'' between 
colored HOMFLY--PT polynomials. One could hope that our approach is conceptual enough
to give new insights on the categorified level (that is, for colored HOMFLY--PT homology) as well.

\subsection{And the categorified story?}

Khovanov's construction of link homology categorifying the Jones polynomial~\cite{kh} can be viewed as a 
categorification of the Temperley--Lieb category $\cal{TL}$, as made precise in the work of Bar-Natan~\cite{bn2}. 
One hence expects that a categorification of our symmetric $\slnn{2}$-web category will be the natural setting for 
a categorification of the colored Jones polynomial. 
We plan to explore exactly this issue in subsequent work, constructing a certain $2$-category of symmetric $\slnn{2}$-foams, 
akin to previous work by Khovanov~\cite{kh1}, Mackaay, 
Sto\v{s}i\'{c} and Vaz~\cite{msv}, Morrison and Nieh~\cite{mn} and Queffelec and the first author~\cite{qr1}.

Such a categorification should give a colored $\slnn{2}$-link homology theory which avoids the 
use of infinite complexes categorifying Jones--Wenzl projectors as in~\cite{ck1},~\cite{fss1} or~\cite{roz}, 
and will be manifestly finite-dimensional (in contrast to those mentioned above, 
as well as Webster's approach~\cite{web}), but (presumably) different from Khovanov's cabling based approach from~\cite{kh2}. 
We point out that work of Hogancamp~\cite{hog} 
has shown how to extract a finite-dimensional colored $\slnn{2}$-link homology theory 
from these infinite-dimensional theories.

We expect the category of symmetric $\slnn{2}$-foams to be related to categorified quantum groups, 
via a symmetric analog of the categorical skew Howe duality pioneered by Cautis, Kamnitzer, and Licata~\cite{ckl} and 
utilized recently in a large body of work by several researchers 
(including the authors of this paper), see~\cite{cautis},~\cite{lqr1},~\cite{mpt},~\cite{my},~\cite{qr1} and~\cite{tub4}. 
Finally, we suspect that a duality between symmetric and traditional foams 
will lead to a precise formulation of ``mirror symmetry'' between (symmetric or skew) colored $\sln$-link homologies.


\begin{thebibliography}{0}
\bibitem{bn2} D.~Bar-Natan, Khovanov's homology for tangles and cobordisms, Geom. Topol. 9 (2005), 1443-1499, online available arXiv:math/0410495.
\bibitem{blm} A.~Beilinson, G.~Lusztig and R.~MacPherson, A geometric setting for the quantum deformation of $\mathfrak{gl}_n$, Duke Math. J. 61-2 (1990), 655-677.
\bibitem{bz1} A.~Berenstein and S.~Zwicknagl, Braided symmetric and exterior algebras, Trans. Amer. Math. Soc. 360-7 (2008), 3429-3472, online available arXiv:math/0504155.
\bibitem{cautis} S.~Cautis, Clasp technology to knot homology via the affine Grassmannian,  Math. Ann. 363 (2015), no. 3-4, 1053-1115, online available arXiv:1207.2074.
\bibitem{ckl} S.~Cautis, J.~Kamnitzer and A.~Licata, Categorical geometric skew Howe duality, Invent. Math. 180-1 (2009), 111-159, online available arXiv:0902.1795.
\bibitem{ckm} S.~Cautis, J.~Kamnitzer and S.~Morrison, Webs and quantum skew Howe duality, Math. Ann. 360-1-2 (2014), 351-390, online available arXiv:1210.6437.
\bibitem{cw1} S.J.~Cheng and W.~Wang, {\em Dualities and representations of Lie superalgebras}, Graduate Studies in Mathematics 144, American Mathematical Society (2012).
\bibitem{cr1} J.~Chuang and R.~Rouquier, Derived equivalences for symmetric groups and $\mathfrak{sl}_2$-categorification, Ann. of Math. (2) 167-1 (2008), 245-298, online available arXiv:math/0407205.
\bibitem{ck1} B.~Cooper and V.~Krushkal, Categorification of the Jones--Wenzl Projectors, Quantum Topol. 3-2 (2012), 139-180, online available arXiv:1005.5117.
\bibitem{doty} S.~Doty, Presenting generalized $q$-Schur algebras, Represent. Theory 7 (2003), 196-213 (electronic), online available arXiv:math/0305208.
\bibitem{fss1} I.~Frenkel, C.~Stroppel and J.~Sussan, Categorifying fractional Euler characteristics, Jones--Wenzl projector and $3j$-symbols, Quantum Topol. 3-2 (2012), 181-253, online available arXiv:1007.4680.
\bibitem{homfly} P.~Freyd, D.~Yetter, J.~Hoste, W.B.R.~Lickorish, K.~Millett, and A.~Ocneanu, A New Polynomial Invariant of Knots and Links, Bull. Amer. Math. Soc. (N.S.) 12-2 (1985), 239-246.
\bibitem{gs} S.~Gukov and M.~Sto\v{s}i\'{c}, Homological algebra of knots and BPS states, Geom. Topol. Monogr. 18 (2012), 309-367, online available arXiv:1112.0030.
\bibitem{hog} M.~Hogancamp, A polynomial action on colored $\mathfrak{sl}(2)$ link homology, online available arXiv:1405.2574.
\bibitem{ho2} R.~Howe, {\em Perspectives on invariant theory: Schur duality, multiplicity-free actions and beyond}, The Schur lectures, Israel Math. Conf. Proc. 8 (1992), 1-182.
\bibitem{ho1} R.~Howe, Remarks on classical invariant theory, Trans. Amer. Math. Soc. 313-2 (1989), 539-570.
\bibitem{ja} J.C.~Jantzen, {\em Lectures on quantum groups}, Graduate Studies in Mathematics 6, American Mathematical Society (1996).
\bibitem{jonespoly} V.F.R.~Jones, A polynomial invariant for knots via von Neumann algebras,  Bull. Amer. Math. Soc. (N.S.) 12-1 (1985), 103-111.
\bibitem{jones} V.F.R.~Jones, Index for subfactors, Invent. Math. 72-1 (1983), 1-25.
\bibitem{kl} L.~Kauffman and S.~Lins, {\em Temperley--Lieb recoupling theory and invariants of 3-manifolds}, Annals of Math. Studies 134, Princeton University Press (1994).
\bibitem{kh} M.~Khovanov, A categorification of the Jones polynomial, Duke Math. J. 101-3 (2000), 359-426, online available arXiv:math/9908171.
\bibitem{kh2} M.~Khovanov, Categorifications of the colored Jones polynomial, J. Knot Theory Ramifications 14-1 (2005), 111-130, online available arXiv:math/0302060.
\bibitem{kh1} M.~Khovanov, $\mathfrak{sl}_3$ link homology, Algebr. Geom. Topol. 4 (2004), 1045-1081, online available arXiv:math/0304375.
\bibitem{kup} G.~Kuperberg, Spiders for rank $2$ Lie algebras, Comm. Math. Phys. 180-1 (1996), 109-151, online available arXiv:q-alg/9712003.
\bibitem{lqr1} A.D.~Lauda, H.~Queffelec and D.E.V.~Rose, Khovanov homology is a skew Howe $2$-representation of categorified quantum $\mathfrak{sl}(m)$, 
Algebr. Geom. Topol. 15 (2015), no. 5, 2517-2608, online available arXiv:1212.6076.
\bibitem{lu1} G.~Lusztig, {\em Introduction to Quantum Groups}, Reprint of the 1994 edition, Modern Birkh0Š1user Classics (2010).
\bibitem{mpt} M.~Mackaay, W.~Pan and D.~Tubbenhauer, The $\mathfrak{sl}_3$-web algebra, Math. Z. 277-1-2 (2014), 401-479, online available arXiv:1206.2118.
\bibitem{msv} M.~Mackaay, M.~Sto\v{s}i\'{c} and P.~Vaz, $\mathfrak{sl}(N)$ link homology using foams and the Kapustin-Li formula, Geom. Topol. 13-2 (2009), 1075-1128, online available arXiv:0708.2228.
\bibitem{my} M.~Mackaay and Y.~Yonezawa, The $\mathfrak{sl}(N)$ web categories and categorified skew Howe duality, online available arXiv:1306.6242.
\bibitem{mn} S.~Morrison and A.~Nieh, On Khovanov's cobordism theory for $\mathfrak{su}_3$ knot homology, J. Knot Theory Ramifications 17-9 (2008), 1121-1173, online available arXiv:math/0612754.
\bibitem{moy} H.~Murakami, T.~Ohtsuki and S.~Yamada, HOMFLY polynomial via an invariant of colored plane graph, Enseign. Math. (2) 44-3-4 (1998), 325-360.
\bibitem{pt} J.H.~Przytycki and P.~Traczyk, Conway Algebras and Skein Equivalence of Links, Proc. Amer. Math. Soc. 100-4 (1987), 744-748.
\bibitem{qr2} H.~Queffelec and D.E.V.~Rose, Sutured annular Khovanov--Rozansky homology, 
Trans. Amer. Math. Soc. 370 (2018), no. 2, 1285-1319, online available arXiv:1506.08188.
\bibitem{qr1} H.~Queffelec and D.E.V.~Rose, The $\mathfrak{sl}_n$ foam $2$-category: 
A combinatorial formulation of Khovanov--Rozansky homology via categorical skew-Howe duality, 
Adv. Math. 302 (2016), 1251-1339, online available arXiv:1405.5920.
\bibitem{qs} H.~Queffelec and A.~Sartori, HOMFLY--PT and Alexander polynomials from a doubled Schur algebra, Quantum Topol. 9 (2018), no. 2, 323-347, online available arXiv:1412.3824.
\bibitem{qs1} H.~Queffelec and A.~Sartori, Mixed quantum skew Howe duality and link invariants of type A, online available arXiv:1504.01225.
\bibitem{roz} L.~Rozansky, An infinite torus braid yields a categorified Jones--Wenzl projector, Fund. Math. 225-1 (2014), 305-326, online available arXiv:1005.3266.
\bibitem{sar1} A.~Sartori, Categorification of tensor powers of the vector representation of $\Uq(\mathfrak{gl}(1|1))$, Selecta Math. (N.S.) 22 (2016), no. 2, 669-734, online available arXiv:1305.6162.
\bibitem{rtw} G.~Rumer, E.~Teller and H.~Weyl, Eine f\"ur die Valenztheorie geeignete Basis der bin\"aren Vektorinvarianten, Nachrichten von der Ges. der Wiss. Zu G\"ottingen. Math.-Phys. Klasse (1932), 498-504.
\bibitem{tl} H.N.V.~Temperley and E.H.~Lieb, Relations between the ``percolation'' and ``colouring'' problem and other graph-theoretical problems associated with regular planar lattices: some exact results for the ``percolation'' problem,  Proc. Roy. Soc. London Ser. A 322-1549 (1971), 251-280.
\bibitem{tub4} D.~Tubbenhauer, $\mathfrak{sl}_n$-webs, categorification and Khovanov--Rozansky homologies, online available arXiv:1404.5752.
\bibitem{tvw} D.~Tubbenhauer, P.~Vaz and P.~Wedrich, Super $q$-Howe duality and web categories, 
Algebr. Geom. Topol. 17 (2017), no. 6, 3703-3749, online available arXiv:1504.05069.
\bibitem{tur} V.~Turaev, {\em Quantum invariants of knots and $3$-manifolds}, second revised edition, de Gruyter Studies in Mathematics (2010).
\bibitem{web} B.~Webster, Knot invariants and higher representation theory, 
Mem. Amer. Math. Soc. 250 (2017), no. 1191, v+141 pp, online available arXiv:1309.3796.
\bibitem{wenzl} H.~Wenzl, On sequences of projections, C. R. Math. Rep. Acad. Sci. Canada 9-1 (1987), 5-9.

\end{thebibliography}
\end{document}